\documentclass[10pt]{amsart}

\usepackage{amsfonts,amssymb}
\input{diagrams}
\newtheorem{proposition}{Proposition}[section]
\newtheorem{lemma}[proposition]{Lemma}
\newtheorem{corollary}[proposition]{Corollary}
\newtheorem{theorem}[proposition]{Theorem}
\theoremstyle{definition}
\newtheorem{definition}[proposition]{Definition}
\newtheorem{example}[proposition]{Example}

\theoremstyle{remark}
\newtheorem{remark}[proposition]{Remark}
\newtheorem{remarks}[proposition]{Remarks}
\newcommand{\thlabel}[1]{\label{th:#1}}
\newcommand{\thref}[1]{Theorem~\ref{th:#1}}
\newcommand{\selabel}[1]{\label{se:#1}}
\newcommand{\seref}[1]{Section~\ref{se:#1}}
\newcommand{\lelabel}[1]{\label{le:#1}}
\newcommand{\leref}[1]{Lemma~\ref{le:#1}}
\newcommand{\prlabel}[1]{\label{pr:#1}}
\newcommand{\prref}[1]{Proposition~\ref{pr:#1}}
\newcommand{\colabel}[1]{\label{co:#1}}
\newcommand{\coref}[1]{Corollary~\ref{co:#1}}

\newcommand{\relabel}[1]{\label{re:#1}}
\newcommand{\reref}[1]{Remark~\ref{re:#1}}

\newcommand{\reslabel}[1]{\label{res:#1}}
\newcommand{\resref}[1]{Remarks~\ref{res:#1}}

\newcommand{\exlabel}[1]{\label{ex:#1}}
\newcommand{\exref}[1]{Example~\ref{ex:#1}}

\def\hk{\mbox{$\widehat{K}$}}
\def\ra{\rightarrow}
\def\cd{\cdot}

\newcommand\smi{\mbox{$S^{-1}$}}
\def\ot{\otimes}
\def\va{\varepsilon}

\def\un{\underline}

\def\mf{\mathfrak}

\def\mfA{\mf {A}}
\def\mfa{\mf {a}}
\def\le{\langle}
\def\ri{\rangle}
\def\l{\lambda}
\def\bl{\mbox{${\bf l}$}}
\def\br{\mbox{${\bf r}$}}

\def\va{\varepsilon}
\def\v{\varphi}

\def\tr{\triangleright}

\def\rh{\rightharpoonup}
\def\lh{\leftharpoonup}

\def\ra{\rightarrow}

\def\d{\delta}

\def\ov{\overline}
\def\cal{\mathcal}
\def\un{\underline}

\newcommand{\mfb}{\mbox{$\mf {b}$}}

\begin{document}
\title[Braided from CQT Hopf algebras]
{Braided Hopf algebras obtained from coquasitriangular Hopf
algebras}
\author{Margaret Beattie}
\address{Department of Mathematics and Computer Science, Mount Allison University,
Sack-ville, NB E4L 1E6, Canada} \email{mbeattie@mta.ca}\thanks{
Support for the first author's research, and partial support for the
second author, came from an NSERC Discovery Grant.  The second
author held  a postdoctoral fellowship at Mount Allison University
during this project. He would like to thank Mount Allison University
for their warm hospitality. }
\author{Daniel Bulacu}
\address{Faculty of Mathematics and Informatics, University
of Bucharest, Str. Academiei 14, RO-010014 Bucharest 1, Romania}
\email{dbulacu@al.math.unibuc.ro}

\begin{abstract}
Let $(H, \sigma)$ be a coquasitriangular Hopf algebra, not
necessarily finite dimensional. Following methods of Doi and
Takeuchi, which parallel the constructions of Radford in the case of
finite dimensional quasitriangular Hopf algebras, we define
$H_\sigma$, a sub-Hopf algebra of $H^0$, the finite dual of $H$.
Using the generalized quantum double construction and the theory of Hopf algebras
with a projection, we associate to $H$ a braided Hopf algebra structure
in the category of Yetter-Drinfeld modules over   $H_\sigma^{\rm cop}$.
Specializing to $H={\rm SL}_q(N)$, we obtain explicit formulas which endow
${\rm SL}_q(N)$ with a braided Hopf algebra structure within the
category of left Yetter-Drinfeld modules over $U_q^{\rm ext}({\rm sl}_N)^{\rm cop}$.
\end{abstract}
\maketitle

%%%%%%%%%%%%%%%%%%%%%%%%%%%%%%%%%%%%%%%%%%%%%%%%%%%%%
\section{Introduction}
%%%%%%%%%%%%%%%%%%%%%%%%%%%%%%%%%%%%%%%%%%%%%%%%%%%%%
\par The Drinfeld double construction plays
an important role in the theory of minimal quasitriangular Hopf
algebras since the Drinfeld double of a finite  dimensional Hopf
algebra is minimal quasitriangular. Conversely, any minimal
quasitriangular Hopf algebra is a quotient of a Drinfeld double
\cite{rad}. This second statement is a consequence of the following
facts. If $(H, R)$ is a quasitriangular Hopf algebra then from an
expression for $R\in H\ot H$ of minimal length, one can define two
finite dimensional sub-Hopf algebras of $H$, denoted by $R_{(l)}$
and $R_{(r)}$, respectively, and there is a quasitriangular Hopf
algebra morphism between the Drinfeld double of $R_{(l)}$ and $H$
whose image is $H_R=R_{(l)}R_{(r)}$, the minimal quasitriangular
Hopf algebra associated to $(H, R)$.  Note that there is an
isomorphism $R_{(l)}^{*\rm cop}\cong R_{(r)}$ and this isomorphism,
together with the evaluation map, gives a  skew pairing from
$R_{(r)}\ot R_{(l)}$ to $k$.\vspace*{1mm}

\par  Now suppose that $(H, \sigma)$ is a coquasitriangular Hopf
algebra, not necessarily finite dimensional.  The analogues of the
Hopf algebras $R_{(l)}$ and $R_{(r)}$, denoted by $H_l$ and $H_r$,
are, in general, infinite dimensional sub-Hopf algebras of $H^0$,
the finite dual of $H$, and  $H_rH_l = H_lH_r$ is a sub-Hopf algebra
of $H^0$ denoted by $H_\sigma$.  In general, there is no analogue to
the isomorphism $R_{(l)}^{*\rm cop}\cong R_{(r)}$, but the
coquasitriangular map $\sigma$ induces a skew pairing on $H_r \ot
H_l$. Although the classical Drinfeld double is not available in
this setting, there exists a generalized quantum double, in the
sense of Majid \cite{majfoundations} or Doi and Takeuchi
\cite{doitak}, for Hopf algebras in duality, so that the skew
pairing between $H_r$ and $H_l$  gives a generalized quantum double,
denoted $D(H_r, H_l)$. In general, $D(H_r, H_l)$ has no
(co)quasitriangular structure (see \prref{1.8} and \reref{3.4});
however, there is a surjective Hopf algebra morphism from $D(H_r,
H_l)$ to $H_\sigma$ given by multiplication in $H^0$. Note that, if
$(H, \sigma)$ is finite dimensional then $H_{\sigma}=H^*_R$, the
minimal quasitriangular Hopf algebra associated to $(H^*, R)$, where
$R$ is the matrix of $H^*$ obtained from $\sigma$. Thus, in the
finite dimensional case, $(H_{\sigma}, R)$ is quasitriangular.
In the infinite dimensional case, neither
$H_{\sigma}$ nor $H^0$ is guaranteed a
quasitriangular structure.\vspace*{1mm}

\par Drinfeld \cite{dri} observed that for a finite
dimensional quasitriangular Hopf algebra $H$, its double $D(H)$ is a
Hopf algebra with a projection.  Majid \cite{majdoubles} proved that
the converse of Drinfeld's result also holds, and computed on $H^*$
the braided Hopf algebra structure associated to this projection
\cite{r}. Here, for $(H ,\sigma)$ coquasitriangular but not
necessarily finite dimensional, using the evaluation pairing, we
construct the generalized quantum double $D(H_\sigma^{\rm cop}, H)$.
In the finite dimensional case, $H_{\sigma}^{\rm cop}= (H^*_R)^{\rm
cop}$ has a quasitriangular structure and there is a Hopf algebra
morphism from $D(H_\sigma^{\rm cop}, H)$ to $H_\sigma^{\rm cop}$
covering the natural inclusion, so $ D(H_\sigma^{\rm cop}, H)$ is a
Hopf algebra with a projection. We show that, for $H$ not
necessarily finite dimensional, $D(H_\sigma^{\rm cop}, H)$ is still
a Hopf algebra with projection, and we describe explicitly the
induced structure on $H$ of a braided Hopf algebra in the category
of left Yetter-Drinfeld modules over $H_\sigma^{\rm cop}$.\vspace*{1mm}

\par  This paper is organized as follows. In Section
\ref{preliminaries}, Preliminaries,  we first define pairings and
skew pairings on Hopf algebras and give the construction of the
generalized quantum double. Coquasitriangular Hopf algebras give
important examples of Hopf algebras with a skew pairing.  As well,
we describe the structure of a Hopf algebra ${\mf A}$ with
projection (see \cite{r}), so that ${\mf A}$ is isomorphic to a
Radford biproduct $B\times {\cal K}$, where $B$ is a braided Hopf
algebra in the category of left Yetter-Drinfeld modules over ${\cal
K}$.\vspace*{1mm}

\par  In Section \ref{doubleprojections}, we  first work in a rather
general context. For $U, V$ bialgebras such that there is an
invertible pairing from $U \ot V$ to $k$, we show that there is a
projection from $D(U^{\rm cop},V)$ to $U^{\rm cop}$ covering the
natural inclusion if and only if there is a bialgebra morphism
$\gamma: V \rightarrow U^{\rm cop}$ satisfying a relation analogous
to the coquasitriangularity condition for a skew pairing. (See also
\cite{doitak} and \cite{sch}.) Now if $U, V$ are Hopf algebras with
an invertible pairing and $\gamma$ exists as above, then $V$ has a
Hopf algebra structure in the category of left Yetter-Drinfeld
modules over $U^{\rm cop}$. We describe this structure explicitly.
For example, if $(U, R)$ is quasitriangular, then    such a map
$\gamma$ exists.  Then, for $H$ finite dimensional we prove that
there is a projection $\pi: D(H)= D(H^{*\rm cop}, H)\ra H^{* \rm
cop}$ covering the inclusion $H^{*\rm cop}\hookrightarrow D(H)$ if
and only if $H$ is coquasitriangular (\prref{3.9}). A fortiori, $H$
has a braided Hopf algebra structure. Moreover, this structure
allows us to obtain the left version of the transmutation theory for
coquasitriangular Hopf algebras (\reref{4.3}).\vspace*{1mm}

\par The transmutation theory for coquasitriangular Hopf algebras is
due to Majid \cite{majtransm}. Using the dual reconstruction
theorem, he associated to any coquasitriangular Hopf algebra $(H,
\sigma)$, a braided commutative Hopf algebra $\un{\un{H}}$ in the
category of right $H$-comodules, called the function algebra braided
group associated to $H$. In Section \ref{tt} we show  that
$\un{\un{H}}$ can be obtained by writing a generalized quantum
double as a Radford biproduct using the projection $\pi : D(H, H)\ra
H$ given by multiplication in $H$. To this projection corresponds a
braided Hopf algebra $\un{H}$  in the category of left
$H$-comodules.  By considering the co-opposite case, we obtain,
after some identification, that $(\un{H^{\rm cop}})^{\rm
cop}=\un{\un{H}}$, as braided Hopf algebras (\prref{4.2}).\vspace*{1mm}

\par In Section \ref{dc}, we describe the ``dual case'', and show that
$D(H_{\sigma}^{\rm cop}, H)$ is always a Hopf algebra with a
projection. This follows from a more general construction where we
take $A$ and $X$ to be sub-Hopf algebras of $H$, $H_{r_A}$ and
$H_{l_X}$ the corresponding sub-Hopf algebras of $H_r$ and $H_l$
respectively, and $H_{X,A}$ the corresponding sub-Hopf algebra of
$H_\sigma$. Then there is a Hopf algebra projection from $D(H^{\rm
cop}_{X,A}, X)$ to $H^{\rm cop}_{X,A}$ which covers the inclusion.
We stress the fact that this ``dual case" give rise  to a new
braided Hopf algebra structure on $X$ (denoted by $\un{X}$) in the
category of left $H_{X, A}^{\rm cop}$ Yetter-Drinfeld modules. This
new construction cannot be viewed as an example of the transmutation
theory because the transmutation theory associates to an ordinary
coquasitriangular Hopf algebra a braided Hopf algebra in the
category of corepresentations over itself. We show that, in fact,
the two constructions above are related by a non-canonical braided
functor. More exactly, we show that there exists a braided functor
$\mathbf{F}: {\cal M}^X\ra {}_{H_{X, A}^{\rm cop}}^{H_{X, A}^{\rm
cop}}{\cal YD}$ such that $\mathbf{F}({\un{\un X}})=\un{X}$
(\thref{unXtransm}). Nevertheless, we think  it  worthwhile to have
the construction of $\un{X}$. This is first of all  because $\un{X}$
lies in a category of left Yetter-Drinfeld modules over $H_{X,
A}^{\rm cop}$ which, in general, may not have a quasitriangular
structure. (See the example in \seref{slqn}). Secondly, this is
because evaluation gives a duality between $H_{X, A}$ and $X$ and
the associated quantum double, $D(H_{X,A}^{\rm cop}, X)$, is
isomorphic to $\un{X} \times H_{X, A}^{\rm cop}$.\vspace*{1mm}

\par In Section \ref{fdcase} we present the finite
dimensional case in full detail, linking the results of this paper
to those of Radford described above.
Also, we note that starting with a finite dimensional
quasitriangular Hopf algebra $(H, R)$ and taking $A=X=H^*$ then the
corresponding braided Hopf algebra $\un{H^*}$ is precisely the
categorical dual of $\un{\un{H^{\rm cop}}}$, the associated
enveloping algebra braided group of $(H^{\rm cop}, R_{21})$
constructed by Majid in \cite{majqg} (\resref{5.8}).\vspace*{1mm}

\par Finally, in \seref{slqn}, we apply the constructions of Section
\ref{dc}  to the coquasitriangular Hopf algebra $H:={\rm SL}_q(N)$.
By direct computations we show that $H_l$ and $H_r$ are the
Borel-like Hopf algebras $U_q({\mf b}_{+})$ and $U_q({\mf b}_{-})$,
respectively, associated to $U_q^{\rm ext}({\rm sl}_N)$, and obtain
that $H_{\sigma}=U_q^{\rm ext}({\rm sl}_N)$. The general theory
above leads to more conceptual proofs for some well known results.
Namely, there exist dual parings of the pairs of Hopf algebras
$U_q({\mf b}_{+})^{\rm cop}$ and $U_q({\mf b}_{-})$, $U_q^{\rm
ext}({\rm sl}_N)$ and ${\rm SL}_q(N)$, and $U_q({\rm sl}_N)$ and
${\rm SL}_q(N)$, respectively (\coref{4.7} and \coref{4.10}), and
$U_q^{\rm ext}({\rm sl}_N)$ and $U_q({\rm sl}_N)$ are factors of
generalized quantum doubles (\coref{4.8} and \reref{4.9}). Also, the
Hopf algebra structure of $U_q({\rm sl}_N)$ is achieved in a natural
way (\reref{4.5} and \reref{4.9}). The rest of the \seref{slqn} is
dedicated to a description of the braided Hopf algebra structure of
${\rm SL}_q(N)$ in the category of left
$U_q^{\rm ext}({\rm sl}_N)^{\rm cop}$ Yetter-Drinfeld modules.
The explicit formulas for this braided Hopf algebra structure can
be found in \prref{6.14} and \thref{6.15}.\vspace*{1mm}

\par Concluding, the general results in this paper have led not only
to some properties of the couple $(U_q^{\rm ext}({\rm sl}_N), {\rm
SL}_q(N))$ well known in quantum group theory, but also to some new
ones. Namely, ${\rm SL}_q(N)$ has a braided Hopf algebra structure
in the category of left Yetter-Drinfeld modules over $U_q^{\rm
ext}({\rm sl}_N)^{\rm cop}$ which coincides with that of the image
of the transmutation object of ${\rm SL}_q(N)$ through a
non-canonical braided functor, and the corresponding Radford
biproduct is isomorphic to the generalized quantum double associated
to this dual pair.

%%%%%%%%%%%%%%%%%%%%%%%%%%%%%%%%%%%%%%%%%%%%%%
\section{Preliminaries}\label{preliminaries}
%%%%%%%%%%%%%%%%%%%%%%%%%%%%%%%%%%%%%%%%%%%%%%
\par Throughout, we work over a field $k$ and maps are assumed to be
$k$-linear.  Any unexplained definitions or notation may be found in
\cite{dnr},\cite{k}, \cite{majfoundations}, \cite{mo} or \cite{sw}.
\par  For ${\cal B}$ a $k$-bialgebra, we write the
comultiplication in ${\cal B}$ as $\Delta(b) = b_1 \ot b_2$
for $b \in {\cal B}$. For ${\mf M}$ a left ${\cal B}$-comodule, we
write the coaction as $\rho({\mf m}) ={\mf m}_{-1}\ot {\mf m}_0 $.
For ${\mf M}$ a right ${\cal B}$-comodule, we will use the subscript
bracket notation to differentiate subscripts from those in
comultiplication expressions,
i.e. we will write $\rho({\mf m})={\mf m}_{(0)}\ot {\mf m}_{(1)}$
for ${\mf m}\in {\mf M}$. For $k$-spaces ${\mf M}$ and ${\mf N}$,
$tw$ will denote the usual twist map from ${\mf M}\ot {\mf N}$ to
${\mf N}\ot {\mf M}$.

\subsection{Pairings on Hopf algebras  and the generalized quantum
double}\label{defD}

We first recall the definition for two bialgebras or Hopf algebras
to be in duality (see \cite[Section 1]{takGL} or \cite[Section
1.4]{majfoundations}); this notion was first introduced by Takeuchi
and called a Hopf pairing.

\par Throughout this section $U$ and $V$ will denote bialgebras over
$k$.

\begin{definition} A bilinear form   $\le , \ri :\ U\ot V\ra k$
 is called a pairing of $U$ and $V$ if
\begin{eqnarray}
&&\le mn, x\ri =\le m, x_1\ri \le n, x_2\ri , \label{d1}\\
&&\le m, xy\ri =\le m_1, x\ri \le m_2, y\ri ,\label{d2} \\
&&\le 1, x\ri =\va (x),~~\le m, 1\ri =\va (m), \label{d3}
\end{eqnarray}
for all $m, n\in U$ and $x, y\in V$. Then $U$ and $V$ are said to
be in duality.
\end{definition}

\begin{remarks}\reslabel{2.0}
(i)  A bilinear form from $U \ot V$ to $k$ is called a skew pairing
if it is a pairing from $U^{\rm cop} \ot V$ to $k$ or, equivalently,
a pairing from $U \ot V^{op}$ to $k$.

(ii) The form $\le , \ri$ is a pairing of $U$ and $V$ if and only if
there is a bialgebra morphism $\phi :\ U\ra V^0$, defined by
$\phi(u)(v) = \le u,v \ri$, if and only if there exists a bialgebra
morphism $\psi : V\ra U^0$ defined by $\psi(v)(u) = \le u,v\ri $. If
$U$ and $V$ are Hopf algebras, these maps are Hopf algebra maps.

(iii) If $U$ and $V$ are Hopf algebras and $\le , \ri$
is a pairing of $U$ and $V$, then this bilinear form is
invertible in the convolution algebra ${\rm Hom}(U \otimes V, k)$ and its
inverse is the bilinear form which maps $u\ot v$ to $\le S_U(u),
v\ri=\le u, S_V(v)\ri $ .

(iv) If  $\le,\ri$ is a (skew) pairing between $U$ and $V$, then
there is a (skew) pairing between the sub-bialgebras $\phi(U)
\subseteq V^0$ and $\psi(V) \subseteq U^0$ defined by the bilinear
form $B: \phi(U)\ot \psi(V)\ra k$,
\[
B(\phi(u), \psi(v)) = \le u, v\ri.
\]
The form $B$ is well-defined since $\phi(u)= \phi(u')$ if and only if
$\le u, - \ri = \le u', - \ri$ and
$\psi(v)=\psi(v')$ if and only if $\le-, v\ri = \le -, v'\ri$. It is
straightforward to check that $B$ is a (skew) pairing.

(v) If $\le , \ri$ is a pairing of the bialgebras $U$ and $V$ and
${\mf U}\subseteq U$, ${\mf V}\subseteq V$ are sub-bialgebras then
the restriction of $\le , \ri$ to ${\mf U}\ot {\mf V}$ is a pairing
of ${\mf U}$ and ${\mf V}$. As above, the bialgebras $\phi({\mf U})$
and $\psi({\mf V})$ are also in duality.
\end{remarks}

\begin{example}
For $H$ a bialgebra, then the evaluation map provides a duality
between $H^0$, the finite dual of $H$, and $H$, that is, $\le f,
h\ri =f(h)$ for all $f\in H^0$ and $h\in H$.

\end{example}

Another class of examples for bialgebras in duality is provided by
coquasitriangular bialgebras, also called braided bialgebras. This
concept is dual to the idea of quasitriangular bialgebras (see
\cite{k} or \cite{majfoundations} for the defintion). We recall the
notion of coquasitriangularity.

\begin{definition}
A bialgebra $H$ is called coquasitriangular (CQT for short) if there
exists a convolution invertible $k$-bilinear skew pairing $\sigma :
H\ot H\ra k$, i.e., for all $h, h',g \in H$,
\begin{eqnarray}
\sigma(hh', g)&=&\sigma(h, g_1)\sigma(h', g_2)\label{cqt1} \\
\sigma(g, hh')&=&\sigma(g_2, h)\sigma(g_1, h')\label{cqt2}\\
\sigma(1,h)&=&\sigma(h,1)=\va(h), \label{cqt3}
\end{eqnarray}
which also satisfies the coquasitriangular condition,
\begin{equation}
\sigma (h_1, h'_1)h_2h'_2=h'_1h_1\sigma (h_2, h'_2). \label{cqt4}
\end{equation}
\end{definition}

\begin{remarks}\label{cqtremarks}
(i) Doi \cite{doi} showed that for   $(H,\sigma)$ a CQT Hopf algebra,
one may define an invertible element $v \in H^*$ by
\begin{equation}\label{v}
v(h)=\sigma (h_1, S(h_2))~~{\rm with \hspace{1mm}  inverse}~~
v^{-1}(h)=\sigma (S^2(h_1), h_2),
\end{equation}
and
\begin{equation}\label{p2}
S^2(h) = v^{-1}(h_1)h_2 v(h_3),~~\forall~~h\in H.
\end{equation}

Similarly, the element $u\in H^*$ defined by $u(h)=\sigma (h_2,
S(h_1))$,  is invertible with inverse $u^{-1}(h) =
\sigma(S^2(h_2),h_1)$ and also defines the square of the antipode
$S$ of $H$ as a co-inner automorphism of $H$, i.e.,
$S^2(h)=u(h_1)h_2u^{-1}(h_3)$, for all $h \in H$.
In particular, the antipode $S$ is bijective.

(ii)\label{cqtinduality} Let $(H,\sigma)$ be a CQT bialgebra. Then
$\sigma$ is a convolution invertible pairing between $H$ and
$H^{op}$ and also between $H^{\rm cop}$ and $H$.

(iii) Clearly any sub-bialgebra of a CQT bialgebra $(H,\sigma)$ is CQT.

(iv) If $(H, \sigma)$ is CQT, then so are $(H^{\rm op},
\sigma^{-1})$ and $(H^{\rm cop}, \sigma^{-1})$ where, if $H$ is a
Hopf algebra, $\sigma^{-1}=\sigma \circ (S_H \ot Id_H)=\sigma \circ
(Id_H\ot S_H^{-1})$ is the convolution inverse of $\sigma$
\cite[1.2]{doi}. Moreover, $(H, \sigma ^{-1}_{21}=\sigma ^{-1}\circ
{\it tw})$ is another CQT structure for $H$, so that  $(H^{\rm
op},\sigma _{21}:=\sigma \circ {\it tw})$ and  $(H^{\rm cop},\sigma
_{21} )$ are CQT also.
\end{remarks}

Let ${\mf B}$, ${\mf H}$ be bialgebras and $\wp$ an invertible skew pairing
$\wp: {\mf B}\ot {\mf H}\rightarrow k$.  We now define the bialgebra (Hopf
algebra) structure on ${\mf B}\otimes {\mf H}$ to be studied in this paper. We
follow the presentation in \cite{doitak}.

For $A$ a bialgebra, an invertible bilinear form $\tau$ on $A$ is called a
unital 2-cocycle if for all $a,b,c \in A$,
\begin{equation}\nonumber
\tau(a_1,b_1)\tau(a_2b_2,c) = \tau(b_1,c_1)\tau(a,b_2c_2) ~~{\rm
and}~~ \tau(a,1) = \tau(1,a) = \va(a).
\end{equation}

If $\tau$ is a 2-cocycle on $A$, then we may form  $A^\tau$, the
bialgebra which has coalgebra structure from $A$ but with the new
multiplication
\[
a\bullet b = \tau(a_1,b_1) a_2b_2\tau^{-1}(a_3,b_3).
\]
It is shown in \cite{doitak} that since  $\wp$ is an invertible skew
pairing on ${\mf B}\ot {\mf H}$, then the bilinear form on ${\mf
B}\ot {\mf H}$ defined by $\tau(\mfb \ot {\mf h}, \mfb' \ot {\mf
h}')=\va(\mfb) \wp(\mfb', {\mf h}) \va({\mf h}')$ is a unital
2-cocycle. Then we may form the bialgebra $({\mf B}\ot {\mf
H})^{\tau}$. As a  coalgebra, $({\mf B}\ot {\mf H})^{\tau}={\mf
B}\ot {\mf H}$, the unit is $1\ot 1$, and the multiplication is
defined by
\begin{equation}\label{mult}
(\mfb \otimes {\mf h})(\mfb'\otimes {\mf h}') =
\wp(\mfb'_1, {\mf h}_1)\wp^{-1}(\mfb'_3, {\mf h}_3)\mfb \mfb'_2\otimes {\mf h}_2{\mf h}'.
\end{equation}

For ${\mf B}, {\mf H}$ Hopf algebras with bijective antipodes, the inverse of
$\wp$ as a skew-pairing on ${\mf B}\ot {\mf H}$ is given by
$\wp^{-1}(\mfb, {\mf h})=\wp(S_{\mf B}(\mfb), {\mf h})=\wp(\mfb, S_{\mf H}^{-1}({\mf h}))$.
Then $({\mf B}\ot {\mf H})^\tau$, also often denoted by ${\mf B}\bowtie_\tau{\mf H}$,
has antipode given by
\begin{eqnarray}\label{antipode}
S(\mfb \otimes {\mf h})&=&(1\otimes S_{\mf H}({\mf h}))
\bullet(S_{\mf B}(\mfb)\otimes 1)\nonumber \\
&=&\wp(S_{\mf B}(\mfb_3), S_{\mf H}({\mf h}_3))
\wp^{-1}(S_{\mf B}(\mfb_1) , S_{\mf H}({\mf h}_1))S_{\mf B}(\mfb_2)\ot S_{\mf H}({\mf h}_2).
\end{eqnarray}

This construction is a special case of  the generalized quantum
double (see \cite[Chapter 7]{majfoundations}) for matched pairs of
Hopf algebras.

\begin{definition}\label{gendouble}
We denote $({\mf B}\ot {\mf H})^\tau$ by $D({\mf B}, {\mf H})$ and will refer to
this bialgebra (Hopf algebra) with multiplication as in (\ref{mult}),
coalgebra structure from the tensor product ${\mf B}\ot {\mf H}$, and antipode
as in (\ref{antipode}) as a generalized quantum double.
\end{definition}

In fact, in many of our constructions, we will begin with a pairing
$\le, \ri: U \ot V \rightarrow k$ which is then a skew pairing
$\rho$ from $U^{\rm cop} \ot V \rightarrow k$.  Then, continuing to
write subscripts in $U$, and using the fact that $\le S_U^{-1}(m),
S_V(x)\ri =\le m, x\ri$, we have that the formulas for the
multiplication, comultiplication and antipode in $D(U^{\rm cop}, V)$
are given by
\begin{eqnarray}
\label{Dmult}
(m\otimes x)(n\otimes y)
&=&\le n _3 , x_1\ri \le S_U^{-1}(n_1), x_3\ri mn_2\otimes x_2y\\
\label{Dcomult}
\Delta(m \ot x)&=&(m_2 \ot x_1) \ot (m_1 \ot y_2)\\
\label{Dantipode}
S(m \otimes x)&=&\le m_1, x_3\ri \le m_3 , S^{-1}_V(x_1)\ri S_U^{-1}(m_2)\ot S_V(x_2).
\end{eqnarray}

\begin{example}(cf. \cite[7.2.5]{majfoundations})
If $H$ is a finite dimensional Hopf algebra, then $H^*$ and $H$ are
in duality via the evaluation map as mentioned above and the double
$D(H^{* {\rm cop}}, H)$ is the usual Drinfeld double, denoted in
this case by $D(H)$.
\end{example}

It is shown in \cite{majfoundations} that if $(H, \sigma)$ is CQT,
so is the double $(H \ot H)^\tau = D(H , H)$, where $\tau$ is the unital
$2$-cocycle on $H\ot H$ defined by the skew-pairing $\sigma$.

The next proposition shows that this generalizes to
$D({\mf B}, {\mf H})$ where ${\mf B}, {\mf H}$ are CQT Hopf algebras.

\begin{proposition}\prlabel{1.8}
Let ${\mf B}, {\mf H}$ be Hopf algebras,
$\wp: {\mf B}\ot {\mf H}\rightarrow k$ a skew pairing, and $D=D({\mf B}, {\mf H})=
({\mf B}\ot {\mf H})^\tau$ the generalized quantum double of Definition
\ref{gendouble}.  Then $D$ is CQT if and only if ${\mf B}$ and ${\mf H}$ are.
\end{proposition}
\begin{proof}
Since sub-Hopf algebras of a CQT Hopf algebra are CQT,
then if $D$ is CQT, so are ${\mf B}$ and ${\mf H}$.

Now suppose that $({\mf B}, \sigma_{\mf B})$ and $({\mf H},
\sigma_{\mf H})$ are CQT. Then we form the CQT Hopf algebra $(A,
\sigma)=\left({\mf B}\ot {\mf H}, (\sigma_{\mf B}\otimes \sigma_{\mf
H})\circ (Id_{\mf B}\otimes{\it tw}\otimes Id_{\mf H})\right)$. From
the above discussion, $\tau : A\ot A\rightarrow k$ is a unital
2-cocycle, where $\tau (\mfb \ot {\mf h}, \mfb'\ot {\mf
h}')=\va(\mfb) \wp(\mfb', {\mf h})\va({\mf h}')$. Then by \cite[p.61
(2.24)]{majfoundations}, $A^\tau$ is CQT via the bilinear form
$\omega$ defined by $\omega(\mf a, {\mf a }')=\tau({\mf a}'_1, \mf a
_1)\sigma(\mf a_2, {\mf a'}_2) \tau^{-1}(\mf a_3, {\mf a'}_3)$.
Specifically,  $(D= A^\tau, \omega)$ is coquasitriangular  where
$$\omega(\mfb \ot {\mf h}, \mfb' \ot {\mf h}')= \wp(\mfb_1, {\mf
h}'_1)\sigma_{\mf B}(\mfb_2, \mfb'_1) \sigma_{\mf H}({\mf h}_1,{\mf
h}'_2) \wp(S_{\mf B}(\mfb'_2), {\mf h}_2).$$ (By   using
(\ref{cqt4}) for $\wp$ generously, one can even check the
coquasitriangularity conditions (\ref{cqt1}) to (\ref{cqt4}) for
$\omega$ directly.)
\end{proof}

Next we show how the doubles $D(U^{\rm cop}, V)$ and $D(\phi(U^{\rm
cop}), \psi(V))$ are related.

\begin{lemma}
Let $U,V$ be Hopf algebras and $\le, \ri$ a pairing of $U$ and $V$.
Then $\phi(U) \subseteq V^0$ and $\psi(V) \subseteq U^0$ are also
Hopf algebras with a pairing $B$ on $\phi(U) \ot \psi(V)$ defined by
$B(\phi(u), \psi(v)) = \le u,v \ri$. Then $\phi \ot \psi: D(U^{\rm
cop}, V) \rightarrow D(\phi(U^{\rm cop}), \psi(V))$ is a surjection
of Hopf algebras.
\end{lemma}
\begin{proof} From \resref{2.0} and the fact that $\phi \ot \psi:
U^{\rm cop} \ot V \rightarrow \phi(U^{\rm cop}) \ot \psi(V)$ is a
Hopf algebra surjection, it remains only to show that $\phi \ot
\psi$ respects the multiplication in the double.  This can be
checked by a straightforward computation, or by noting that for
$D(U^{\rm cop}, V) = (U^{\rm cop} \ot V)^\tau$ and $D(\phi(U^{\rm
cop}), \psi(V))=( \phi(U^{\rm cop}) \ot \psi(V) )^{\tau'}$ for
cocycles $\tau, \tau'$ as above, then $\tau(m\ot x, n \ot y)=
\tau'(\phi(m)\ot \psi(x), \phi(n)\ot \psi(y))$.
\end{proof}

%%%%%%%%%%%%%%%%%%%%%%%%%%%%%%%%%%%%%%%%%%%%%%%%%%%%%%%%%%%%%%%%%%%
\subsection{Hopf algebras with projection}\label{projection}
%%%%%%%%%%%%%%%%%%%%%%%%%%%%%%%%%%%%%%%%%%%%%%%%%%%%%%%%%%%%%%%%%%%%%%
Let ${\cal K}$ be a bialgebra. Recall that a left Yetter-Drinfeld
module over ${\cal K}$ is a left ${\cal K}$-module ${\mf M}$ which
is also a left ${\cal K}$-comodule,  such that the following
compatibility relation  holds. For all $\kappa\in {\cal K}$ and
${\mf m}\in M$:
\begin{equation}\label{lyd1}
\kappa_1{\mf m}_{-1}\ot \kappa_2\cd {\mf m}_{0}=(\kappa_1\cd {\mf m})_{-1}\kappa_2\ot
(\kappa_1\cd {\mf m})_{0},
\end{equation}
where ${\cal K}\ot {\mf M}\ni \kappa\ot {\mf m}\mapsto \kappa\cd
{\mf m}\in {\mf M}$ is the left ${\cal K}$-action.   The category of
left Yetter-Drinfeld modules over ${\cal K}$ and $k$-linear maps
that preserve the ${\cal K}$-action and ${\cal K}$-coaction is
denoted by ${}_{\cal K}^{\cal K}{\cal YD}$.

\par The category ${}_{\cal K}^{\cal K}{\cal YD}$ is pre-braided.
If ${\mf M}, {\mf N}\in {}_{\cal K}^{\cal K}{\cal YD}$ then ${\mf M}\ot {\mf N}$ is a left
Yetter-Drinfeld module over ${\cal K}$ via the structures defined by
\begin{equation}\label{lyd2}
\kappa\cd ({\mf m}\ot {\mf n})=\kappa_1\cd {\mf m}\ot \kappa_2\cd {\mf n}~~{\rm and}~~
{\mf m}\ot {\mf n}\mapsto {\mf m}_{-1}{\mf n}_{-1}\ot {\mf m}_0\ot {\mf n}_0,
\end{equation}
for all $\kappa\in {\cal K}$, ${\mf m}\in {\mf M}$ and ${\mf n}\in {\mf N}$.
The pre-braiding is given by
\begin{equation}\label{lyd3}
{\bf c}_{{\mf M}, {\mf N}}({\mf m}\ot {\mf n})={\mf m}_{-1}\cd {\mf n}\ot {\mf m}_{0}.
\end{equation}
If ${\cal K}$ is a Hopf algebra then ${\bf c}$ is invertible, so ${}_{\cal K}^{\cal K}{\cal YD}$
is a braided monoidal category.

\par The structure of a Hopf algebra with projection was given in
\cite{r}. More precisely, if ${\cal K}$ and ${\mf A}$ are   Hopf
algebras with  Hopf algebra maps $ {\cal K} \pile{\rTo^{i}\\
\lTo_{\pi}}{\mf A}$ such that $\pi \circ i=Id_{\cal K}$, then there
exists a braided Hopf algebra $B$ in the category of left
Yetter-Drinfeld modules ${}_{\cal K}^{\cal K}{\cal YD}$ such that
${\mf A}\cong B\times {\cal K}$ as Hopf algebras, where $B\times
{\cal K}$ denotes Radford's biproduct between $B$ and ${\cal K}$
(for more details see \cite{r}).

As $k$-vector space $B=\{\mfa\in {\mf A}\mid \mfa_1\ot \pi
(\mfa_2)=\mfa\ot 1\}$. Now, $B$ is a ${\cal K}$-module subalgebra of
${\mf A}$ where ${\mf A}$ is a left ${\cal K}$-module algebra via
the left adjoint action induced by $i$, that is $\kappa \tr
_i\mfa=i({\kappa}_1)\mfa i(S({\kappa}_2))$, for all $\kappa\in {\cal
K}$ and $\mfa\in {\mf A}$. Moreover, $B$ is an algebra in the
braided category ${}_{\cal K}^{\cal K}{\cal YD}$ where the left
coaction of ${\cal K}$ on $B$ is given for all $b\in B$ by
\begin{equation} \label{lbi}
\l _{B}(b)=\pi (b_1)\ot b_2.
\end{equation}

Also, as $k$-vector space, $B$ is the image of the $k$-linear map
$\Pi :\ {\mf A}\ra {\mf A}$ defined for all $\mfa\in {\mf A}$ by
\begin{equation} \label{Pi}
\Pi (\mfa)=\mfa_1i(S(\pi (\mfa_2))).
\end{equation}

For all $\mfa\in \mfA$, we define
\begin{equation} \label{dbi}
\un {\Delta }(\Pi (\mfa))=\Pi (\mfa_1)\ot \Pi (\mfa_2).
\end{equation}

This makes $B$ into a coalgebra in ${}_{\cal K}^{\cal K}{\cal YD}$
and a bialgebra in ${}_{\cal K}^{\cal K}{\cal YD}$. The counit of
$B$ is $\un {\va }= \va \mid _B$. Moreover, we have that $B$ is a
braided Hopf algebra in ${}_{\cal K}^{\cal K}{\cal YD}$ with
antipode $\un {S}$ given by
\begin{equation} \label{ant}
\un {S}(b)=i(\pi (b_1))S_{\mfA}(b_2),
\end{equation}
where $S_{\mfA}$ is the antipode of $\mfA$.

The Hopf algebra isomorphism $\chi :\ B\times {\cal K}\ra \mfA$ is
given by
\begin{equation}\label{th}
\chi (b\times \kappa)=bi(\kappa),
\end{equation}
for all $b\in B$ and $\kappa\in {\cal K}$.

\par Note that the description of a Hopf with a projection in terms of a braided
Hopf algebra is due to Majid \cite{majqg}.
%%%%%%%%%%%%%%%%%%%%%%%%%%%%%%%%%%%%%%%%%%%%%%
\section{Generalized quantum doubles which are Radford
biproducts}\label{doubleprojections}
%%%%%%%%%%%%%%%%%%%%%%%%%%%%%%%%%%%%%%%%%%%%%%
\setcounter{equation}{0}
It is well-known (Majid \cite{majdoubles}) that if $H$ is a finite
dimensional Hopf algebra, then the Drinfeld double $D(H)$ is a
Radford biproduct.
 In this section, we give necessary and sufficient conditions for a
generalized quantum double $D=D(U^{\rm cop}, V)$ to be a Radford
biproduct $B \times U^{\rm cop} $, and determine the structure of
$B$ as a Hopf algebra in $_{U^{\rm cop}}^{U^{\rm cop}}\cal{YD}$.

Suppose first that $U$ and $V$ are bialgebras in duality with $\le,
\ri : U\otimes V \rightarrow k$ an invertible pairing, so that $\rho
=\le, \ri$ is an invertible skew pairing on $U^{\rm cop} \ot V$. We
form the generalized quantum double $D = D(U^{\rm cop}, V)$ as in
Subsection \ref{defD}. There are bialgebra morphisms $i: U^{\rm
cop}\rightarrow D(U^{\rm cop}, V)$ given by $i(m)=m\otimes 1$ and
$j: V \rightarrow D(U^{\rm cop}, V)$ given by $j(x)=1\ot x$.

\begin{proposition}\label{piiffgamma}
Let $U, V$ be as above.

(i) There exists a bialgebra projection $\pi$ from $D=D(U^{\rm cop},
V)$ to $U^{\rm cop}$ that splits $i$ if and only if there is a
bialgebra map $\gamma: V \rightarrow U^{\rm cop}$ such that for all
$y\in V$, $m \in U$, we have
\begin{equation}\label{gamma}
\gamma(y)m = \rho^{-1}(m_1,  y_3)\rho( m_3, y_1)m_2\gamma(y_2).
\end{equation}

(ii) Similarly, there exists a bialgebra projection from $D$ to $V$
that splits $j$ if and only if there is a bialgebra morphism $\mu:
U^{\rm cop} \rightarrow V$ such that  for all $y\in V, m\in U$, we
have
\begin{equation}
y \mu(m)=\rho^{-1}(m_1, y_3)\rho( m_3,y_1 )\mu(m_2)y_2.
\end{equation}
If $U, V$ are Hopf algebras, then these maps are Hopf algebra
morphisms.
\end{proposition}

\begin{proof}
(i) Suppose that there is a bialgebra morphism $\gamma: V
\rightarrow U^{\rm cop}$ satisfying (\ref{gamma}). Define $\pi$ to
be $\pi (m\ot x)=m\gamma(x)$, for all $m\in U$ and $x\in V$. Then
for $m\in U$, we have that $\pi \circ i(m) = \pi(m\ot 1)=
m\gamma(1)=m$. Furthermore, by \cite[2.4]{doitak} with $B=J=U^{\rm
cop}$, $H=V$, $\alpha = Id_U$, $\beta = \gamma$, we have that $\pi$
is a bialgebra morphism.

Conversely, given $\pi$, define $\gamma$ by $\gamma = \pi \circ j$
where $j: V \rightarrow D(U^{\rm cop}, V)$ is defined by $j(x) = 1
\otimes x$. Then, since $\pi,j$ are bialgebra maps, so is $\gamma$.
To verify that (\ref{gamma}) holds, we compute
\begin{eqnarray*}
\gamma(y)m&=&(\pi \circ j)(y)m\\
&=&\pi(1 \otimes y)\pi(m \otimes 1)
=\pi((1 \otimes y)(m\otimes 1))\\
&=& \rho^{-1}(m_1, y_3 )\rho(m_3, y_1)m_2\gamma(y_2).
\end{eqnarray*}

(ii) The proof of (ii) is analogous.
\end{proof}

\begin{example}\label{dte}(cf. \cite[3.1]{doitak})
Let $\sigma$ be an invertible
skew pairing on a bialgebra $H$ and form $D(H,H)$.  Then the
identity map $Id_H$ satisfies $(\ref{gamma})$ if and only if
$(H, \sigma)$ is CQT if and only if the multiplication map
$\pi : D(H,H)\ra H$, $\pi(h \ot l) = hl$ is a bialgebra map.
\end{example}

\begin{example}\exlabel{2.3}
In the setting of Proposition \ref{piiffgamma}, if $V$ is
quasitriangular via $R = R^1 \ot R^2 \in V \ot V$, then the map
$\pi: D(U^{\rm cop}, V) \ra V$ defined by $\pi(m\ot x) = \le m, R^1
\ri R^2x$ is a bialgebra projection. Here the map $\mu: U^{\rm
cop}\ra V$ is given by $\mu(m)=\le m, R^1\ri R^2$. (The details can
be found in \cite[2.5]{doitak}.)

Likewise, if $U$ is quasitriangular with
the $R$-matrix $R=R^1\ot R^2\in U\ot U$,
then $\pi: D(U^{\rm cop}, V)\ra U^{\rm cop}$
defined by $\pi (m\ot x)=\le R^2, x\ri mS(R^1)$
is a Hopf algebra projection. In this case the map
$\gamma: V\ra U^{\rm cop}$ is given by
$\gamma(y)=\le R^2, y\ri S(R^1)$.
\end{example}

\begin{remarks}\label{equivgamma}
(i) In    later   computations we will use  that
  (\ref{gamma}) is equivalent to
\begin{equation}\label{gamma'}
\rho( m_1, y_2) \gamma(y_1)m_2=\rho( m_2, y_1) m_1\gamma(y_2).
\end{equation}
If $m=\gamma(x)$, $x\in V$, then (\ref{gamma}) becomes
\begin{equation}\label{gammagamma}
\rho ( \gamma(x_2), y_2) \gamma(y_1)\gamma(x_1)=\rho ( \gamma(x_1),
y_1)\gamma(x_2)\gamma(y_2).
\end{equation}

(ii) Note that if the above map $\gamma: V \rightarrow U^{\rm cop}$
is injective then the map $\sigma: V \ot V \ra k$ defined by
$\sigma(x,y) = \rho(\gamma(x),y) = \le \gamma(x),y \ri$ gives $V$ a
CQT structure. The relations (\ref{cqt1}), (\ref{cqt2}),
(\ref{cqt3}) are easy to check and (\ref{cqt4}) is equivalent to
\[
\le\gamma(x_1),y_1\ri x_2y_2 = \le\gamma(x_2),y_2\ri y_1x_1 .
\]
This equation holds if and only if it holds when the injective map
$\gamma$ is applied to both sides, i.e., when (\ref{gammagamma})
holds.
\end{remarks}

If $U, V$ are Hopf algebras with bijective antipodes, then for
$\rho$ a skew pairing from $U^{\rm cop} \ot V $ to $k$, we have
$\rho^{-1}(m, x)=\le S^{-1}_U(m), x\ri = \le m, S^{-1}_V(x)\ri$. In
this case, we have the identities below which are useful in the
following computations and also provide generalizations of  the
equations describing the square of the antipode in Remarks
\ref{cqtremarks}(i).

\begin{proposition}
Let $U, V$ be Hopf algebras in duality and assume that there exists a map $\gamma$ as
in Proposition \ref{piiffgamma}. Then:
\begin{itemize}
\item[(i)]
The map $\vartheta\in V^*$ defined by $\vartheta (x)=\le \gamma(x_1), S_V(x_2)\ri$,
for all $x\in V$, is convolution invertible with
$\vartheta^{-1}(x)=\le \gamma(S_V^2(x_1)), x_2\ri$. Moreover, for any $x\in V$,
\begin{equation}\label{vgamma}
\gamma(S_V^2(x))=\vartheta^{-1}(x_1)\gamma(x_2)\vartheta(x_3).
\end{equation}
\item[(ii)]
Similarly, the map $\upsilon\in V^*$ defined by $\upsilon(x)=\le \gamma(x_2), S_V(x_1)\ri$,
for all $x\in V$, is convolution invertible with
$\upsilon ^{-1}(x)=\le \gamma(S_V^2(x_2)), x_1\ri$. In addition, for all $x\in V$,
\begin{equation}\label{vSgamma}
\gamma(S_V^2(x))=\upsilon(x_1)\gamma(x_2)\upsilon^{-1}(x_3).
\end{equation}
\end{itemize}
\end{proposition}
\begin{proof}
We only sketch the proof for (i); the rest of the details are left
to the reader.

For all $x\in V$ we have
\begin{eqnarray*}
\vartheta (x_1)\gamma(S_V^2(x_2))
&=&\le \gamma(x_3), S_V(x_4)\ri \gamma(x_1)\gamma(S_V(x_2))\gamma(S_V^2(x_5))\\
&=&\le \gamma(x_3), S_V(x_4)\ri \gamma(x_1)S_U^{-1}
\left(\gamma(S_V(x_5))\gamma(x_2)\right)\\
&{{(\ref{gammagamma})}\atop =}&\le \gamma(x_2), S_V(x_5)\ri
\gamma(x_1)S_U^{-1}\left(\gamma(x_3)\gamma(S_V(x_4))\right)\\
&=&\le \gamma(x_2), S_V(x_3)\ri \gamma(x_1)=\gamma(x_1)\vartheta(x_2),
\end{eqnarray*}
and, in a similar manner, one can prove that
\begin{equation}\label{vt2}
\gamma(S^2_V(x_1))\vartheta^{-1}(x_2)=\vartheta^{-1}(x_1)\gamma(x_2).
\end{equation}
Now, for  $x\in V$, using the fact that $\gamma: V\ra U^{\rm cop}$
is a Hopf algebra map, we have
\begin{eqnarray*}
&&\hspace*{-1cm}
\vartheta(x_1)\vartheta^{-1}(x_2)=\vartheta(x_1)\le \gamma(S_V^2(x_2)), x_3\ri
=\vartheta(x_2)\le \gamma(x_1), x_3\ri \\
&&\hspace*{3mm} =\le \gamma(x_2), S_V(x_3)\ri \le \gamma(x_1),
x_4\ri = \le \gamma(x_1), S_V(x_2)x_3\ri =\va(x).
\end{eqnarray*}
Similarly, using (\ref{vt2}) we can show
that $\vartheta^{-1}(x_1)\vartheta(x_2)=\va(x)$, so we are done.
\end{proof}

Now suppose that $U, V$ are Hopf algebras with bijective antipodes
and we have a Hopf algebra projection $\pi$ from $D(U^{\rm cop}, V)$
to $U^{\rm cop}$ that splits $i$. Then there exists a Hopf algebra
$B$ in the category of Yetter-Drinfeld modules $^{U^{\rm
cop}}_{U^{\rm cop}}{\cal YD}$ such that $D = D(U^{\rm cop},V) \cong
B \times U^{\rm cop}$, a Radford biproduct. From Subsection
\ref{projection}, we know that
\[
B = \{a \in D\mid a_1 \otimes \pi(a_2) = a \otimes 1 \}.
\]

\begin{proposition}\label{thetabijection}
The map $\theta: V \rightarrow B$ given by $\theta(y) =
\gamma(S_V^{-1}(y_2))\otimes y_1 = S_U(\gamma(y_2))\otimes y_1$ is
a bijection.
\end{proposition}
\begin{proof}
Clearly the map $\theta$ is injective since $\va \circ \gamma=\va$.
It remains to show that $Im(\theta) = B$. For $y \in V$, we have
$\theta(y) \in B$ since
\begin{eqnarray*}
\theta(y)_1 \otimes \pi(\theta(y)_2)&=&(S_U(\gamma(y_4))\otimes y_1)
\otimes \pi(S_{U}(\gamma(y_3))\otimes y_2)\\
&=&(S_U(\gamma(y_4))\otimes y_1)\otimes \gamma(S_V^{-1}(y_3)y_2)\\
&=&(S_U(\gamma (y_2))\otimes y_1)\otimes 1 =\theta(y) \otimes 1.
\end{eqnarray*}
Conversely, suppose that $m\otimes y\in B$, i.e., $(m_2\otimes
y_1)\otimes m_1\gamma(y_2)=(m\otimes y)\otimes 1\in D\otimes U^{\rm
cop}$. Then we have
\begin{eqnarray*}
m\otimes y &=&S_U(1)m\otimes y\\
&=&S_U(m_1 \gamma(y_2))m_2\otimes y_1  \\
&=&S_U(\gamma(y_2))S_U(m_1)m_2\otimes y_1\\
&=&\va_U(m)\theta(y).
\end{eqnarray*}
Similarly, if $z=\sum\limits_i m_i\otimes y_i \in B$, then
$z=\theta(\sum\limits_i\va(m_i) y_i)$.
\end{proof}

Now we denote by $\un{V}$ the vector space $V$ with the structure
of a Yetter-Drinfeld module induced by that of $B$.

\begin{proposition}\label{leftleftYD}
The structure of the left $U^{\rm cop}$ Yetter-Drinfeld module
$\underline{V}$  is given by the left action and left coaction
\begin{eqnarray}
&&m\rhd y=\le m_1, S^{-1}_V(y_1)\ri y_2\le m_2,
S^{-2}_V(y_3)\ri =\le m, S_V^{-1}(S_V^{-1}(y_3)y_1)\ri y_2;\\
&&\lambda_{\underline{V}}(y)=\gamma(S^{-1}_V(y_3)y_1)\ot y_2.
\end{eqnarray}
\end{proposition}
\begin{proof}
For $m\in U$ and $y\in \underline{V}$, we have that
$m\rhd y=\theta^{-1}(m\rhd_i \theta(y))$ by Section \ref{projection}
and using (\ref{Dmult}), (\ref{gamma}) and the fact that
$\le m, v\ri =\le S_U^{-1}(m), S_V(v)\ri$, we compute
\begin{eqnarray*}
&&\hspace*{-1.5cm}
i(m_2)\theta(y)i(S_U^{-1}(m_1))\\
&=&(m_2\otimes 1)(\gamma(S_V^{-1}(y_2))\otimes y_1)(S_U^{-1}(m_1)\otimes 1)\\
&=&m_4\gamma(S^{-1}_V(y_4))S_U^{-1}(m_2)\le S^{-1}_U(m_3),
S^{-1}_V(y_3)\ri \otimes y_2\le S^{-1}_U(m_1),y_1\ri\\
&{{(\ref{gamma'})}\atop =}&
m_4S^{-1}_U(m_3)\gamma(S^{-1}_V(y_3))\le S_U^{-1}(m_2), S^{-1}_V(y_4)\ri
\ot y_2\le S_U^{-1}(m_1), y_1\ri\\
&=&\le S^{-1}_U(m_1), y_1\ri \theta(y_2)
\le S_U^{-1}(m_2), S^{-1}_V(y_3)\ri \\
&=&\le m, S^{-1}_V(y_1)S_V^{-2}(y_3)\ri \theta(y_2),
\end{eqnarray*}
and this concludes the proof of the formula for the action.
We now compute the coaction.
\begin{eqnarray*}
(Id_U\ot \theta)\circ \lambda_{\un {V}}(y)
&=&(\pi \otimes Id_D)\Delta (\theta(y))\\
&=&(\pi \otimes Id_D)\Delta(\gamma(S^{-1}_V(y_2))\otimes y_1)\\
&=&(\pi \otimes Id_D)((\gamma(S_V^{-1}(y_4)) \otimes y_1)
\otimes (\gamma(S_V^{-1}(y_3))\otimes y_2))\\
&=&\gamma(S_V^{-1}(y_4)y_1)\otimes \gamma(S^{-1}_V(y_3))\otimes y_2\\
&=&\gamma(S^{-1}_V(y_3)y_1)\ot \theta(y_2),
\end{eqnarray*}
for all $y\in V$, and thus the formula for the coaction is also verified.
\end{proof}

We now describe the  structure of $\un{V}$ as a Hopf algebra in the
category of left Yetter Drinfeld modules over $U^{\rm cop}$.

\begin{proposition}\label{HopfYD}
The structure of $\un{V}$ as a Hopf algebra in the category
$_{U^{\rm cop}}^{U^{\rm cop}}\cal{YD}$ is given by the formulas:
\begin{eqnarray}\label{ydmult}
x\cdot y&=&\le \gamma(y_2), S_V(x_1)x_3\ri x_2y_1;\\
\un{\Delta}(x)&=&\le \gamma(S_V(x_4)x_6),
S^{-1}_V(x_3)x_1\ri x_2\ot x_5;
\label{ydcomult} \\
\un{S}(x)&=&\le \gamma(x_4), x_1 S_V(x_3)\ri S_V(x_2).
\label{ydantipode}
\end{eqnarray}
The identity in $\un{V}$ is $\theta^{-1}(1\ot 1)=1$
and the counit is $\un{\va}=\va$.
\end{proposition}
\begin{proof} To see (\ref{ydmult}), we compute
\begin{eqnarray*}
\theta(x)\theta(y)&=&(S_U(\gamma (x_2))\ot x_1)(S_U(\gamma (y_2))\ot y_1)\\
&{{(\ref{Dmult})}\atop =}&
\le \gamma(y_2), x_3\ri \le S_U(\gamma(y_4)),
x_1\ri S_U(\gamma(y_3)\gamma(x_4))
\ot x_2y_1\\
&{{(\ref{gammagamma})}\atop =}&
\le \gamma (y_3), x_4 \ri \le S_U(\gamma(y_4)), x_1\ri
S_U(\gamma(x_3)\gamma(y_2))\ot x_2y_1\\
&=&\le \gamma (y_2), x_3\ri \le \gamma(y_3), S_V(x_1)\ri
\theta(x_2y_1)\\
&=&\le \gamma(y_2), S_V(x_1)x_3 \ri \theta(x_2y_1).
\end{eqnarray*}

Similarly, to verify (\ref{ydcomult}), we compute
\begin{eqnarray*}
\underline{\Delta}(\theta(x))&=&(\Pi \ot \Pi)(\Delta(\theta(x)))=
(\Pi \ot \Pi)(\Delta (S_U(\gamma(x_2))\ot x_1))\\
&=&\Pi (S_U(\gamma(x_4))\ot x_1)\ot \Pi(S_U(\gamma(x_3))\ot x_2))\\
&=&\Pi (S_U(\gamma(x_3))\ot x_1)\ot \theta(x_2)\\
&{{(\ref{Pi})}\atop =}&
\left(S_U(\gamma(x_5))\ot x_1\right)
\left(S^{-1}_U(\pi(S_U(\gamma(x_4))\ot x_2))\ot 1\right)
\ot \theta(x_3)\\
&{{(\ref{Dmult}, \ref{d1})}\atop =}&
\le \gamma (x_{12}), S^{-1}_V(x_4)\ri \le S_U^{-1}(\gamma(x_6)), S_V^{-1}(x_5)\ri
\le S_U^{-1}(\gamma(x_8)), x_1\ri \\
&&\times
\le \gamma(x_{10}), x_2\ri
\gamma (S_V^{-1}(x_{13})S_V(x_7)x_{11})\ot x_3\ot \theta(x_9)
\end{eqnarray*}
Now we use  the equivalent form of (\ref{vgamma})
\[
\vartheta^{-1}(S_V^{-1}(y_2))\gamma(S_V^{-1}(y_1))=
\gamma(S_V(y_2))\vartheta^{-1}(S_V^{-1}(y_1)),
\]
to replace $\le S_U^{-1}(\gamma(x_6)), S_V^{-1}(x_5)\ri \gamma(S_V(x_7))$ by
$\le S_U^{-1}(\gamma(x_7)), S_V^{-1}(x_6)\ri \gamma(S_V^{-1}(x_5))$,
and obtain
\begin{eqnarray*}
\underline{\Delta}(\theta(x))&=&
\le \gamma (x_{12}), S^{-1}_V(x_4)\ri \le S_U^{-1}(\gamma(x_7)), S_V^{-1}(x_6)\ri
\le S_U^{-1}(\gamma(x_8)), x_1\ri \\
&&\times
\le \gamma(x_{10}), x_2\ri
\gamma (S_V^{-1}(x_{13})S_V^{-1}(x_5)x_{11})\ot x_3\ot \theta(x_9).
\end{eqnarray*}
From (\ref{gammagamma}), we have
\[
\le \gamma(x_{12}), S_V^{-1}(x_4) \ri \gamma(S_V^{-1}(x_5))
\gamma (x_{11})=\le \gamma(x_{11}), S_V^{-1}
(x_5)\ri \gamma(x_{12})\gamma (S_V^{-1}(x_4)),
\]
so we can conclude that
\begin{eqnarray*}
\underline{\Delta}(\theta(x))&=&
\le  \gamma (x_{11}), S_V^{-1}(x_5)\ri \le S_U^{-1}(\gamma(x_7)), S_V^{-1}(x_6)\ri
\le S_U^{-1}(\gamma(x_8)), x_1 \ri \le \gamma(x_{10}),x_2 \ri \\
&&\times
\gamma (S_V^{-1}(x_{13})x_{12})\gamma(S_V^{-1}(x_4))\ot x_3\ot \theta(x_9)\\
&=&\le \gamma (x_{10}), S_V^{-1}(x_4)\ri \le S_U^{-1}(\gamma(x_6)), S_V^{-1}(x_5)\ri
\le S_U^{-1}(\gamma(x_7)), x_1\ri \\
&&\times
\le \gamma(x_{9}), x_2\ri \theta(x_3)\ot \theta(x_8)\\
&=& \le \gamma(S_V(x_4) x_8), S_V^{-1}(x_3)\ri \le \gamma(S_V(x_5)x_7),
x_1 \ri \theta(x_2) \ot \theta(x_6)\\
&=& \le \gamma(S_V(x_4)x_6), S^{-1}_V(x_3)x_1\ri \theta(x_2)\ot \theta(x_5).
\end{eqnarray*}
Finally, we verify (\ref{ydantipode}). From (\ref{ant}), we have that
\begin{eqnarray*}
\un{S}(\theta (x))&=&
i(\pi(\theta(x)_1))S(\theta(x)_2)\\
&=&(S_U(\gamma(x_4))\gamma(x_1)\ot 1)S(S_U(\gamma(x_3))\ot x_2)\\
&=&(S_U(\gamma(x_4))\gamma(x_1)\ot 1)(1\ot S_V(x_2))(\gamma(x_3) \ot 1)\\
&=&\le S_U^{-1}(\gamma(x_7)), S_V(x_2)\ri \le \gamma(x_5),
S_V(x_4)\ri S_U (\gamma(x_8))\gamma(x_1)\gamma(x_6)\ot S_V(x_3)\\
&=&\le \gamma(x_7), x_2\ri \le \gamma(x_5), S_V(x_4)\ri
S_U(\gamma(x_8))\gamma(x_1)\gamma(x_6)\ot S_V(x_3).
\end{eqnarray*}
By (\ref {gammagamma}) we can replace $\le \gamma(x_7),x_2 \ri
\gamma(x_1) \gamma(x_6)$ by $\le \gamma(x_6), x_1 \ri
\gamma(x_7)\gamma(x_2)$ and we obtain
\begin{eqnarray*}
\un{S}(\theta (x))&=&
\le \gamma(x_6), x_1\ri \le \gamma(x_5), S_V(x_4)\ri
S_U(\gamma(x_8))\gamma(x_7)\gamma(x_2)\ot S_V(x_3)\\
&=&\le \gamma(x_6), x_1\ri \le \gamma(x_5), S_V(x_4)\ri
\gamma(x_2)\ot S_V(x_3)\\
&=&\le \gamma(x_5), x_1S_V(x_4)\ri S_U(\gamma (S_V(x_2)))
\ot S_V(x_3) \\
&=& \le \gamma(x_4), x_1S_V(x_3)\ri \theta(S_V(x_2)).
\end{eqnarray*}
The final statement is clear.
\end{proof}

We noted in \exref{2.3} that if $V$ is quasitriangular then
$D(U^{\rm cop}, V)$ is a Hopf algebra with a projection. So for
$V=(H, R)$ a finite dimensional quasitriangular Hopf algebra and
$U=H^*$ we obtain  Drinfeld's projection \cite{dri}, and by an
analogue of Propositions \ref{leftleftYD} and \ref{HopfYD},  the
structure of $H^*$ as a braided Hopf algebra in ${}_H^H{\cal YD}$
computed by Majid in \cite{majdoubles}. (In fact, it can be proved
that this braided Hopf algebra lies in the image of a canonical
braided functor from ${}_H{\cal M}$ to ${}_H^H{\cal YD}$.)  On the
other hand, if $U = H^*$ is  quasitriangular, and $V=H$, then we
have the following.

\begin{proposition}\prlabel{3.9}
Let $H$ be a finite dimensional Hopf algebra. Then there exists a
Hopf algebra projection $\pi: D(H)\ra H^{*\rm cop}$ covering the
natural inclusion $H^{*\rm cop}\hookrightarrow D(H)$ if and only if
$H$ is CQT. Moreover, if this is the case, then $\pi$ is a
quasitriangular morphism.
\end{proposition}
\begin{proof}
Suppose $(H,\sigma)$ is CQT and let $\{ e_i,e^i \}$ be a dual basis
for $H$. Then  $H^*$ is quasitriangular with $R$-matrix
$R=\sum\limits_{i, j}\sigma(e_i, e_j)e^i\ot e^j$.  Since $D(H)=
D(H^{*\rm cop}, H)$, from \exref{2.3} the map $\pi: D(H)\ra H^{*\rm
cop}$ given by
\[
\pi(h^*\ot h)=\sum\limits_{i, j}\sigma(e_i, e_j)e^j(h)h^*(e^i\circ
S)= \sum\limits_{i}\sigma(e_i, h)h^*(e^i\circ S) ,
\]
for all $h^*\in H^*$ and $h\in H$, is a Hopf algebra morphism which
covers the inclusion $H^{*\rm cop}\hookrightarrow D(H)$.

Conversely, if such a morphism exists, then by Proposition
\ref{piiffgamma} there exists a Hopf algebra morphism $\gamma: H\ra
H^{*\rm cop}$ satisfying
\[
\le h^*_1, x_2\ri \gamma(x_1)h^*_2=\le h^*_2, x_1\ri h^*_1\gamma(x_2),
\]
for all $x\in H$ and $h^*\in H^*$, where $\le , \ri : H^*\ot H\ra k$ is
the evaluation map. It is clear that the above condition is equivalent to
\[
\le \gamma(x_1), y_1\ri x_2y_2=\le \gamma(x_2), y_2\ri y_1x_1,
\]
for all $x, y\in H$. Now, if we define
$\sigma(x, y)=\le \gamma(x), y\ri$, for all $x, y\in H$, then one can
easily check that $\sigma$ defines a CQT structure on $H$
(see also \cite[Theorem 3.3.14]{sch}).

Finally, the canonical $R$-matrix of $D(H)$ is
$R=\sum\limits_{i}(\va\ot e_i)\ot (e^i\ot 1)$. So if
$(H, \sigma)$ is CQT then the above morphism $\pi$ is quasitriangular since
\[
(\pi\ot \pi)(R)
=\sum\limits_{i, j}\sigma(e_j, S^{-1}(e_i))e^j\ot e^i
=\sum\limits_{i, j}\sigma^{-1}(e_j, e_i)e^j\ot e^i,
\]
and, from the dual version of Remarks \ref{cqtremarks} (iv), the
last term  defines a QT structure for $H^{*\rm cop}$.
\end{proof}

For $H$ finite dimensional it is not hard to see
that the categories ${}_{H^{*\rm cop}}^{H^{*\rm cop}}{\cal YD}$ and
${}_H{\cal YD}^H$ are isomorphic as braided monoidal categories (the braided
structures are the ones obtained from the left or right center
construction, see \cite{bcp}). The isomorphism is produced by
the following functor ${\cal F}$.
If $({\mf M}, \cdot ,\lambda)\in {}_{H^{*\rm cop}}^{H^{*\rm cop}}{\cal YD}$
then ${\cal F}({\mf M})={\mf M}$ becomes an object in ${}_H{\cal YD}^H$ via the structure
\[
h\bullet {\mf m}={\mf m}_{-1}(h)m_{0},\hspace{.3cm} {\mf m}\mapsto
e^i\cdot {\mf m}\ot e_i.
\]
${\cal F}$ sends a morphism to itself.

By the above identification, \prref{3.9} and Propositions
\ref{leftleftYD} and \ref{HopfYD} we obtain the following  which may
be viewed as a left version of the transmutation theory for CQT Hopf
algebras. Further details will follow in Section \ref{tt}.

\begin{corollary}\colabel{3.10}
If $H$ is a CQT Hopf algebra then $H$ has a braided Hopf algebra
structure, denoted by $\un{H}$, within ${}_H{\cal YD}^H$. Namely,
$\un{H}$ is a left-right Yetter-Drinfeld module over $H$ via
\begin{eqnarray*}
&&x\mapsto x_{(0)}\ot x_{(1)}:=x_2\ot S^{-1}(x_1)S^{-2}(x_3)\\
&&h\bullet x=\sigma(h, S^{-1}(x_1)S^{-2}(x_3))x_2=
\sigma^{-1}_{21}(S^{-1}(x_3)x_1, h)x_2,
\end{eqnarray*}
for all $h, x\in H$. $\un{H}$ is a Hopf algebra in ${}_H{\cal YD}^H$
with the same unit and counit as $H$ and
\begin{eqnarray*}
&&x\cdot y=\sigma(S(x_1)x_3, S^{-1}(y_2))x_2y_1=\sigma^{-1}_{21}
(y_2, S(x_1)x_3)x_2y_1,\\
&&\un{\Delta}(x)=\sigma(S^{-1}(x_3)x_1, S^{-1}(x_6)x_4)x_2\ot x_5=
\sigma^{-1}_{21}(S(x_4)x_6, S^{-1}(x_3)x_1)x_2\ot x_5,\\
&&\un{S}(x)=\sigma(x_1S(x_3), S^{-1}(x_4))S(x_2)=
\sigma^{-1}_{21}(x_4, x_1S(x_3))S(x_2).
\end{eqnarray*}
\end{corollary}
\begin{proof}
We only note that from the proof
of \prref{3.9} the Hopf algebra
morphism $\gamma: H\ra H^{*\rm cop}$ is given
in this case by $\gamma(h)=l_{S^{-1}(h)}$,
for all $h\in H$.
\end{proof}

If ${\cal B}$ is a bialgebra and ${\mf M}$ is a left Yetter-Drinfeld
module over ${\cal B}$,
then the map $R_{\mf M}:\ {\mf M}\ot {\mf M}\ra {\mf M}\ot {\mf M}$,
$R_{\mf M}({\mf m}\ot {\mf m}')={\mf m}'_{-1}\cd {\mf m}\ot
{\mf m}'_{0}$ is a solution in ${\rm End}({\mf M}^{\ot 3})$ of the quantum
Yang-Baxter equation
\[
R_{12}R_{13}R_{23}=R_{23}R_{13}R_{12}.
\]

Above, we have that $\un{V}$ is a left Yetter-Drinfeld module over
$U^{\rm cop}$. Then a solution $R_V\in End(V\ot V)$ to the quantum
Yang-Baxter equation is given by
\begin{eqnarray*}
R_V(x \ot y)&=&\gamma(S_V(y_3)y_1) \rhd x \ot y_2 \\
&=& \le \gamma(S_V^{-1}(y_3)y_1), S_V^{-1}(S_V^{-1}(x_3)x_1)\ri
x_2\ot y_2.
\end{eqnarray*}

%%%%%%%%%%%%%%%%%%%%%%%%%%%%%%%%%%%%%%%%%%%%%%%%%%%%%%%%%%%%%%%%%%
\section{Coquasitriangular bialgebras and generalized quantum
doubles}\label{cqtbialgs}
%%%%%%%%%%%%%%%%%%%%%%%%%%%%%%%%%%%%%%%%%%%%%%%%%%%%%%%%%%%%%%%%%%%%
\setcounter{equation}{0}
%%%%%%%%%%%%%%%%%%%%%%%%%%%%%%%%%%%%%%%
\subsection{Transmutation theory}\label{tt}
%%%%%%%%%%%%%%%%%%%%%%%%%%%%%%%%%%%%%%%%
As was mentioned in the introduction,  the braided reconstruction
theorem   associates to any CQT Hopf algebra $H$ a braided
commutative Hopf algebra $\un{\un{H}}$ in the category of right
$H$-comodules, ${\cal M}^H$. The goal of this subsection is to show
that ${\un{\un H}}$ can   be   obtained from   the structure of a
generalized quantum double with a projection.

Let $X$ and $A$ be sub-Hopf algebras of a CQT Hopf algebra  $(H,
\sigma)$. Then $\sigma$ induces a skew pairing on $A \ot X$, still
denoted by $\sigma$ and   the generalized quantum double $D(A, X)$
is defined. By (\ref{cqt4}) it follows that $XA=AX$, so $XA$ is a
sub-Hopf algebra of $H$. From \cite[3.1]{doitak}  the map
\[
\pi': D(A, X)\ra AX=XA,~~\pi' (a\ot x)=ax
\]
is a surjective Hopf algebra morphism. Although $D(H,H)$ is a Hopf
algebra with projection, we cannot, in general, make the same claim
for $D(A,X)$. Nevertheless, $XA$ and $X$ are sub-Hopf algebras of
$(H, \sigma)$, so applying the same arguments we find that
\[
\pi :D (AX, X)\ra AX=XA,~~\pi (ax\ot y)=axy
\]
is a surjective Hopf algebra morphism covering the inclusion map
$i: AX\ra D(AX, X)$ with  the  map $\gamma$ from Proposition
\ref{piiffgamma} being the inclusion of $X$ into $AX=XA$. From
Propositions \ref{leftleftYD} and \ref{HopfYD} for
$\gamma :X\hookrightarrow AX$ and
$\le, \ri=\sigma :AX\ot X\ra k$, $X$ has a braided Hopf algebra structure
in the braided monoidal category ${}_{XA}^{XA}{\cal YD}$;
this braided Hopf algebra is denoted, as usual, by
$\un{X}$. The structures are given by
\begin{eqnarray*}
&&xa\tr y=\sigma (xa, S^{-1}(S^{-1}(y_3)y_1))y_2,\\
&&\lambda_{\un{X}}(y)=S^{-1}(y_3)y_1\ot y_2,\\
&&x\cdot y=\sigma(y_2, S(x_1)x_3)x_2y_1,~~\un{1}=1,\\
&&\un{\Delta}(x)=\sigma(S(x_4)x_6, S^{-1}(x_3)x_1)x_2\ot x_5,~~\un{\va}=\va,\\
&&\un{S}(x)=\sigma(x_4, x_1S(x_3))S(x_2),
\end{eqnarray*}
for all $x, y\in X$ and $a\in A$. Next, we show that $\un{X}$ lies
in the image of a canonical braided functor from ${}^X{\cal M}$ to
${}^{XA}_{XA}{\cal YD}$, so this general context reduces to the case
$X=A=H$. We recall some background on CQT Hopf algebras and braided
monoidal categories.

For a CQT bialgebra $({\cal B}, {\varsigma})$ it is well known that its
category of left or right ${\cal B}$ comodules has a braided monoidal structure.
Namely, if ${\mf M}$ and ${\mf N}$ are two left (respectively right) ${\mf B}$
comodules then
\[
{\mf M}\ot {\mf N}\ni {\mf m}\ot {\mf n}\mapsto {\mf m}_{-1}{\mf n}_{-1}
\ot {\mf m}_0\ot {\mf n}_0\in {\cal B}\ot {\mf M}\ot {\mf N}
\]
defines the monoidal structure of ${}^{\cal B}{\cal M}$ and
\begin{equation}\label{lcombr}
c_{{\mf M}, {\mf N}}: {\mf M}\ot {\mf N}\ra {\mf N}\ot {\mf M},~~
c_{{\mf M}, {\mf N}}({\mf m}\ot {\mf n})=\varsigma ({\mf n}_{-1},
{\mf m}_{-1}) {\mf n}_0\ot {\mf m}_0
\end{equation}
defines a braiding on ${}^{\cal B}{\cal M}$, while
\begin{equation}\label{rightmonoidal}
{\mf M}\ot {\mf N}\ni {\mf m}\ot {\mf n}\mapsto
\ot {\mf m}_{(0)}\ot {\mf n}_{(0)}\ot {\mf m}_{(1)}{\mf n}_{(1)}
\in {\mf M}\ot {\mf N}\ot {\cal B}
\end{equation}
gives the monoidal structure of ${\cal M}^{\cal B}$ and
\begin{equation}\label{rcombr}
c_{{\mf M}, {\mf N}}: {\mf M}\ot {\mf N}\ra {\mf N}\ot {\mf M},~~
c_{{\mf M}, {\mf N}}({\mf m}\ot {\mf n})=\varsigma ({\mf m}_{(1)},
{\mf n}_{(1)}) {\mf n}_{(0)}\ot {\mf m}_{(0)}
\end{equation}
provides a braided structure on ${\cal M}^{\cal B}$
(see \cite{k, majfoundations} for terminology).

We denote by ${}^{({\cal B}, \varsigma)}{\cal M}$ the category of
left ${\cal B}$-comodules   with the braiding in (\ref{lcombr}).
Similarly, ${\cal M}^{({\cal B}, \varsigma)}$ is the category of
right ${\cal B}$-comodules with the braiding
defined by (\ref{rcombr}). One can easily see that
${}^{({\cal B}^{\rm cop}, \varsigma _{21})}{\cal M}\equiv
{\cal M}^{({\cal B}, \varsigma)}$,
as braided monoidal categories, where
$\varsigma _{21}=\varsigma \circ {\it tw}$ is the CQT structure of
${\cal B}^{\rm cop}$ defined in Remarks \ref{cqtremarks}.

Secondly, if $({\cal B}, \varsigma)$ is a CQT bialgebra and ${\mf
M}$ a left $\cal{B}$-comodule then ${\mf M}$ is a left
Yetter-Drinfeld module over ${\cal B}$ with the initial comodule
structure and with the ${\cal B}$-action defined by
\begin{equation}\label{modstr}
b\cd {\mf m}:=\varsigma ({\mf m}_{-1}, b){\mf m}_{0}.
\end{equation}
Thus there is  a well defined braided functor ${\mf F}_{({\cal B},
\varsigma)}: {}^{({\cal B}, \varsigma)}{\cal M}\ra {}_{\cal B}^{\cal
B}{\cal YD}$ where ${\mf F}_{({\cal B}, \varsigma)}$ sends a
morphism to itself.

\begin{lemma}\lelabel{4.0}
In the setting above, the braided Hopf algebra $\un{X}$ lies in the
image of the composite of  the canonical functor ${}^{(X,
\sigma_{21}^{-1})}{\cal M}\ra {}^{(XA, \sigma ^{-1}_{21})}{\cal M}$
and the functor ${\mf F}_{(XA, \sigma ^{-1}_{21})}$.
\end{lemma}
\begin{proof}
We have
\[
\lambda _{\un{X}}(x):=x_{-1}\ot x_{0}=S^{-1}(x_3)x_1\ot x_2.
\]
Therefore the $X$-action on $\un{X}$ can be rewritten as
\[
x\tr y=\sigma (x, S^{-1}(y_{-1}))y_{0}=
\sigma ^{-1}(x, y_{-1})x_0=\sigma ^{-1}_{21}(y_{-1}, x)y_0,
\]
and this finishes the proof.
\end{proof}

In general, if ${\cal C}$ is a braided monoidal category with
braiding $c$ then ${\cal C}^{\rm in}$ is ${\cal C}$ as a monoidal
category, but with the mirror-reversed braiding $\tilde{c}_{M,
N}=c^{-1}_{N, M}$. Note that, if $B\in {\cal C}$ is a braided Hopf
algebra with comultiplication $\un{\Delta}$ and bijective antipode
$\un{S}$ then $B^{\rm cop}$, the same object $B$, but with the
comultiplication and antipode
\[
\un{\Delta}^{\rm cop}=c^{-1}_{B, B}\circ \un{\Delta}~~{\rm and}~~
\un{S}^{\rm cop}=\un{S}^{-1 }
\]
respectively, and with the other structure morphisms the same as for
$B$, is a braided Hopf algebra  in the category ${\cal C}^{\rm in}$.
Now, since the transmutation object $\un{\un{H}}$ is a braided Hopf
algebra in ${\cal M}^H$ and not ${}^H{\cal M}$ we must  apply the
above correspondence $X\mapsto \un{X}$ to $(H^{\rm cop},
\sigma_{21})$ rather than $(H, \sigma)$ and thus obtain a braided
Hopf algebra
\[
\un{H^{\rm cop}}\in {}^{(H^{\rm cop}, \sigma^{-1})}{\cal M}\equiv
{\cal M}^{(H, \sigma_{2 1}^{-1})}\equiv {{\cal M}^{(H,
\sigma)}}^{\rm in}.
\]

We can now prove the connection between $\un{H^{\rm cop}}$ and $\un{\un{H}}$.
The structures of $\un{\un{H}}$ can be found in \cite[Example 9.4.10]{majfoundations}.

\begin{proposition}\prlabel{4.2}
Let $(H, \sigma)$ be a CQT Hopf algebra. Then $(\un{H^{\rm cop}})^{\rm cop}=\un{\un{H}}$
as braided Hopf algebras in ${\cal M}^{(H, \sigma)}$.
\end{proposition}
\begin{proof}
One can easily see that $(\un{H^{\rm cop}})^{\rm cop}$ is an object of ${\cal M}^H$
via the structure
\[
h\mapsto h_{(0)}\ot h_{(1)}=h_2\ot S(h_1)h_3,
\]
and that its algebra structure within ${\cal M}^{(H, \sigma)}$ is given by
\[
h\cdot g=\sigma (S^{-1}(h_3)h_1, g_1)h_2g_2=\sigma (S(h_1)h_3, S(g_1))h_2g_2,
\]
for all $h, g\in H$. Clearly, the unit of $(\un{H^{\rm cop}})^{\rm cop}$ is the unit of $H$.

Now, by (\ref{rcombr}) we have that $c^{-1}$, the inverse of the
braiding $c$ of ${\cal M}^{(H, \sigma_{21}^{-1})}$, is defined by
$c^{-1}_{{\mf M}, {\mf N}}({\mf n}\ot {\mf m})=\sigma({\mf n}_{(1)},
{\mf m}_{(1)}) {\mf m}_{(0)}\ot {\mf n}_{(0)}$. Therefore, the
comultiplication of $(\un{H^{\rm cop}})^{\rm cop}$ is given by
\begin{eqnarray*}
&&h\mapsto \sigma (S(h_4)h_6, S^{-1}(h_3)h_1)c^{-1}(h_5, h_2)\\
&&\hspace*{1.5cm}
=\sigma ((h_2)_{(1)}, S^{-1}((h_1)_{(1)}))c^{-1}((h_2)_{(0)}\ot (h_1)_{(0)})\\
&&\hspace*{1.5cm}
=\sigma ^{-1}((h_2)_{(2)}, (h_1)_{(2)})\sigma ((h_2)_{(1)}, (h_1)_{(1)})
(h_1)_{(0)}\ot (h_2)_{(0)}\\
&&\hspace*{1.5cm}
=h_1\ot h_2=\Delta (h),
\end{eqnarray*}
the comultiplication of $H$. The counit of $(\un{H^{\rm cop}})^{\rm cop}$ is $\va$,
the counit of $H$.

Finally, the antipode of $\un{H^{\rm cop}}$ is defined by
${\un S}(h)=\sigma (h_4S^{-1}(h_2), h_1)S^{-1}(h_3)$, so the antipode
of $(\un{H^{\rm cop}})^{\rm cop}$ is given, for all $h\in H$, by
\[
{\un S}^{-1}(h)=\sigma (S^2(h_3)S(h_1), h_4)S(h_2).
\]
Indeed, for all $h\in H$ we have
\begin{eqnarray*}
(\un{S}\circ \un{S}^{-1})(h)&=&
\sigma(h_1S(h_6), S^{-1}(h_7))\sigma (S(h_2)h_4, S(h_5))h_3\\
&{{(\ref{cqt1})}\atop =}&
\sigma(h_1, S^{-1}(h_9))\sigma(S(h_7), S^{-1}(h_8))\sigma (h_2, h_6)
\sigma (h_4, S(h_5))h_3\\
&=&\sigma(h_1, S^{-1}(h_7))v(h_4)\sigma(h_2, h_5)v^{-1}(h_6)h_3\\
&{{(\ref{p2})}\atop =}&
\sigma(h_1, S^{-1}(h_5))\sigma (h_2, S^{-2}(h_4))h_3\\
&{{(\ref{cqt2})}\atop =}&\sigma (h_1, S^{-2}(h_3)S^{-1}(h_4))h_2
=h,
\end{eqnarray*}
as needed. In a similar way we can prove that
$\un{S}^{-1}\circ \un{S}=Id_H$, the details are left to the reader.
Comparing the above structures of $(\un{H^{\rm cop}})^{\rm cop}$
with those of $\un{\un{H}}$ from \cite{majfoundations} we conclude
that$(\un{H^{\rm cop}})^{\rm cop}=\un{\un{H}}$, as braided
Hopf algebras in ${\cal M}^{(H, \sigma)}$.
\end{proof}

\begin{remark}\relabel{4.3}
From \coref{3.10} we can deduce the left version of the
transmutation theory for a CQT Hopf algebra $(H, \sigma)$ as
follows. Observe  that there is a braided functor ${\cal G}: {}^{(H,
\sigma)}{\cal M}\ra {}_H{\cal YD}^H$ defined by ${\cal G}({\mf
M})={\mf M}$ but now viewed as left-right $H$ Yetter-Drinfeld module
via
\[
h\cdot {\mf m}=\sigma(h, S^{-1}({\mf m}_{-1})){\mf m}_{0},~~
{\mf m}\mapsto \rho({\mf m}):={\mf m}_{0}\ot S^{-1}({\mf m}_{-1}).
\]
${\cal G}$ sends a morphism to itself.

Now, one can easily see that $\un{H}$ in \coref{3.10} lies in the image of the
functor ${\cal G}$. In other words we can associate to $(H, \sigma)$ a Hopf algebra
structure in ${}^{(H, \sigma)}{\cal M}$, still denoted by $\un{H}$. Note that
$\un{H}$ has the left $H$-comodule structure given by
\[
x\mapsto \l_{\un{H}}(x)=x_{-1}\ot x_{0}=S^{-1}(x_3)x_1\ot x_2,
\]
for all $x\in \un{H}$, and is a braided Hopf algebra in ${}^{(H, \sigma)}{\cal M}$
via the structure in \coref{3.10}. (In fact, this $\un{H}$ is the braided Hopf
algebra in \leref{4.0} corresponding to $X=(H, \sigma^{-1}_{21})$.) We can
conclude now that $\un{H}^{\rm op, cop}$ is the braided Hopf algebra
in ${}^{(H, \sigma)}{\cal M}$ associated to $(H, \sigma)$ through the
left transmutation theory. By $\un{H}^{\rm op, cop}$ we denote
$\un{H}$ viewed now as a Hopf algebra in ${}^{(H, \sigma)}{\cal M}$ via
\begin{eqnarray*}
&&x\diamond y=:m_{\un{H}}\circ c_{\un{H}, \un{H}}(x, y)=\sigma
(y_{-1}, x_{-1})y_0x_0
=\sigma(y_2, S(x_1)x_3)x_2y_1,\\
&&\un{\Delta}_{\un{H}^{\rm op, cop}}=c^{-1}_{\un{H}, \un{H}}\circ \un{\Delta}
=\Delta ,
\end{eqnarray*}
and the other structure morphisms equal those of  $\un{H}$.
\end{remark}
%%%%%%%%%%%%%%%%%%%%%%%%%%%%%%%%%%%%%%%%%%%%%%%%%%%%%%%%%
\subsection{The ``dual'' case}\label{dc}
%%%%%%%%%%%%%%%%%%%%%%%%%%%%%%%%%%%%%%%%%%%%%%%%%%%%%%%%%%
Let $(H, \sigma)$ be a CQT bialgebra, so that $\sigma$ is a pairing
from $ H \otimes H^{op} $ to $k$ and $  H^{\rm cop} \otimes H$ to $
k$. Then by  \resref{2.0}, we have bialgebra morphisms
\[
\phi: H^{\rm cop}\rightarrow H^0~~{\rm defined~~by~~}
\phi(h)(l)=\sigma(h,l)
\]
and
\[\psi: H^{op}\rightarrow H^0~~{\rm defined~~by~~}
\psi(h)(l)=\sigma(l,h).
\]
We denote $r_h:=\phi(h)=\sigma(h, -)$ and $H_r=Im(\phi)$.
Similarly, $l_h:=\psi(h)=\sigma(-, h)$ and $H_l = Im(\psi)$.
More generally, if $X$ and $A$ are sub-bialgebras of a
CQT bialgebra $H$ then we define
\begin{eqnarray*}
&& H_{r_A}:= \phi(A^{\rm cop}) = \{r_a = \sigma(a, -) \mid a\in A\}\\
&& H_{l_X}:= \psi(X^{op}) = \{l_x= \sigma(-, x) \mid x\in X\}  .
\end{eqnarray*}

\begin{proposition}\prlabel{2.1}
Let $(H, \sigma)$ be a CQT bialgebra (Hopf algebra) and $X, A\subseteq H$ two
sub-bialgebras (sub-Hopf algebras). Then
\begin{itemize}
\item[(i)] $H_{l_X}$ and $H_{r_A}$ are sub-bialgebras (sub-Hopf algebras)
of $H^0$.  The structure maps for $\psi(X^{op}) = H_{l_X}$ are given by
\begin{equation}
l_{xy}=l_yl_x;~~
l_1=\va;~~
\Delta(l_x)=l_{x_1}\otimes l_{x_2};~~
\va(l_x)=\va(x);~~
S(l_x)=l_{S^{-1}(x)},
\label{Hlstructure}
\end{equation}
and the structure maps for $\phi(A^{\rm cop})= H_{r_A}$ are given by
\begin{equation}
r_{ab} = r_ar_b;~~
r_1=\va;~~
\Delta(r_a)=r_{a_2}\otimes r_{a_1};~~
\va(r_a)=\va(a);~~
S(r_a)= r_{S^{-1}(a)}.
\label{Hrstructure}
\end{equation}
\item[(ii)] The bilinear form from $H_{r_A}\ot H_{l_X}$ to $ k$ defined
for all $a\in A, x\in X$, by $  r_a \ot l_x   \ra \sigma (a, x)$, is
a skew pairing between $H_{r_A} $ and $H_{l_X}$.
\item[(iii)] $H_{r_A} H_{l_X}= H_{l_X} H_{r_A}$ is a sub-bialgebra of $H^0$ and is a
sub-Hopf algebra if $A,X$ are sub-Hopf algebras of the Hopf
algebra $H$.
\end{itemize}
\end{proposition}
\begin{proof}
(i) If $X$ is a sub-bialgebra (sub-Hopf algebra) of $H$, then
$X^{op}$ is a sub-bialgebra (sub-Hopf algebra) of $H^{op}$. Since
$\sigma: H \otimes H^{op}\rightarrow k$ is a duality, by
\resref{2.0}, $\psi: X^{op}\subseteq H^{op} \rightarrow H^0$ is a
bialgebra morphism and thus $H_{l_X}$, the image of $X^{op}$ under
$\psi$, is a sub-bialgebra (sub-Hopf algebra) of $H^0$ with the
structure maps given above. The proof for $H_{r_A}$   is similar,
using the map $\phi$.

(ii) Follows directly from \resref{2.0} since $\sigma: A^{\rm cop}
\otimes X^{op}\rightarrow k$ is a skew pairing.

(iii) We refer to \cite[3.2 (b)]{doitak}. Here it is shown that the
map from $D(H^{\rm cop}, H^{\rm op})$ to $H^0$ defined by $a \ot
h\rightarrow r_al_h$ is a bialgebra map, and our statement follows.
The proof is based on the observation that
\begin{eqnarray}
&&l_yr_a=\sigma (a_1, S(y_3))\sigma (a_3, y_1)r_{a_2}l_{y_2}\label{f1}
\mbox{ or,  equivalently,  }\\
&&r_al_y=\sigma (\smi (a_3), y_1)
\sigma (a_1, y_3)l_{y_2}r_{a_2},\label{f2}
\end{eqnarray}
for any $y\in X$, $a\in A$.
\end{proof}

We now assume that $H$ is a Hopf algebra. We denote by $H_{X,A}$ the
Hopf subalgebra of $H^0$ equal to $H_{l_X}H_{r_A} = H_{r_A}H_{l_X}$.

\begin{proposition}\prlabel{1.2}
Let $(H, \sigma)$ be a CQT Hopf algebra, $X$, $A$ and ${\cal H}$
sub-Hopf algebras of $H$, and $H_{l_X}$ and $H_{r_A}$ the sub-Hopf
algebras of $H^0$ from \prref{2.1}. Let $D(H_{r_A}, H_{l_X})$ be the
generalized quantum double from the skew pairing in \prref{2.1} (ii)
induced by $\sigma$.
\begin{itemize}
\item[(i)] The map $f: D(H_{r_A}, H_{l_X})\ra H_{X, A}$, $f(r_a\ot
l_x)=r_al_x$, $a\in A$, $x\in X$, is a surjective Hopf algebra
morphism.
\item[(ii)] Then the evaluation map $\le , \ri :\ H_{X, A}\ot {\cal H}\ra k$ defined for
all $x\in X$, $a\in A$ and $h\in {\cal H}$ by
\begin{equation}
\le l_xr_a, h\ri =\sigma (h_1, x)\sigma (a, h_2) = l_x(h_1)r_a(h_2)
\end{equation}
provides a duality between the Hopf algebras $H_{X, A}$ and ${\cal
H}$, and therefore between $D(H_{r_A}, H_{l_X})$ and ${\cal H}$ via
the map $f$ in (i). In particular, $\le, \ri$ is a skew pairing on
$H_{X, A}^{\rm cop}\ot {\cal H}$ and $D(H_{r_A}, H_{l_X})^{\rm
cop}\ot {\cal H} $.
\end{itemize}
\end{proposition}
\begin{proof}
Statement (i) follows directly from \cite[2.4]{doitak} with $\alpha
= Id_{H_{r_A}}$, $\beta = Id_{H_{l_X}}$ and using (\ref{f1}).
Statement (ii) is immediate.
\end{proof}

From now on, we assume that $(H, \sigma)$ is a CQT Hopf algebra with
$X, A$ and the evaluation pairing $\le, \ri$  as above.
Consider the generalized quantum double $D(H_{X, A}^{\rm cop}, X)$.
From (\ref{Dmult}), the multiplication in $D(H_{X, A}^{\rm cop}, X)$
is given by
\begin{equation}\label{sm}
(l_xr_a\ot y)(l_{x'}r_b\ot y')=\le l_{x'_1}r_{b_3}, \smi (y_3)\ri
\le l_{x'_3}r_{b_1}, y_1\ri l_xr_al_{x'_2}r_{b_2}\ot y_2y',
\end{equation}
for all $x, x', y, y'\in X$, $a, b\in A$, and since
$D(H_{X, A}^{\rm cop}, X) = H_{X, A}^{\rm cop} \otimes X$ as a coalgebra,
the comultiplication is given by
\begin{equation}\label{sc}
\Delta (l_xr_a\ot y)=(l_{x_2}r_{a_1}\ot y_1)
\ot (l_{x_1}r_{a_2}\ot y_2).
\end{equation}
The unit is $\va \ot 1$ and the counit is defined by
$\va (l_xr_a\ot y)=\va (x)\va (a)\va (y)$.

We now prove that there is a Hopf algebra projection from
$D(H_{X, A}^{\rm cop}, X)$ to $H_{X, A}^{\rm cop}$ which covers the
canonical inclusion
$i:\ H_{X, A}^{\rm cop}\hookrightarrow D(H_{X, A}^{\rm cop}, X)$.

\begin{proposition}\prlabel{1.5}
The map $\pi :\ D(H_{X, A}^{\rm cop}, X)\ra H_{X, A}^{\rm cop}$
defined by
\begin{equation}\label{mappi}
\pi (l_xr_a\ot y)=l_xr_al_{S^{-1}(y)},
\end{equation}
for all $x, y\in X$, $a\in A$, is a Hopf algebra morphism such
that $\pi \circ i=id_{H_{X, A}^{\rm cop}}$.
\end{proposition}
\begin{proof}
Define $\gamma: X \rightarrow H^{\rm cop}_{l_X} \subseteq H^{\rm
cop}_{X, A}$ by $\gamma(x)=l_{S^{-1}(x)}$.  We show that $\gamma $
is a bialgebra morphism and that (\ref{gamma}) holds.

For $x,y \in X$, we have
\[
\gamma(xy)=l_{S^{-1}(xy)}=l_{S^{-1}(y)S^{-1}(x)}
=l_{S^{-1}(x)}l_{S^{-1}(y)}=\gamma(x)\gamma(y),
\]
so that $\gamma$ preserves multiplication. Similarly,
\[
(\gamma \otimes \gamma )(x_1\otimes x_2)=l_{S^{-1}(x_1)}\otimes l_{S^{-1}(x_2)}
=\Delta_{H^{\rm cop}_{X, A}}(\gamma(x)),
\]
so that comultiplication is preserved.
Clearly $\gamma$ preserves the unit and counit. To verify (\ref{gamma}),
we compute
\begin{eqnarray*}
\gamma(x)l_y&=&l_{S^{-1}(x)}l_y  = l_{yS^{-1}(x)}\\
&{{(\ref{cqt1})}\atop =}&
\sigma(S^{-1}(x_2),y_2)\sigma(x_1,y_3)l_{y_1S^{-1}(x_3)}\\
&{{(\ref{cqt4})}\atop =}&
\sigma(S^{-1}(x_3),y_1)\sigma(x_1,y_3)l_{ S^{-1}(x_2)y_2}\\
&=& \le l_{y_1}, S^{-1}(x_3)\ri \le l_{y_3},x_1\ri l_{y_2}\gamma(x_2).
\end{eqnarray*}
As well we have that
\begin{eqnarray*}
\gamma(x_2)r_a&=&l_{S^{-1}(x_2)}r_a \\
&{{(\ref{f1})}\atop =}&
\sigma(a_1, S(S^{-1}(x_2)))\sigma(a_3, S^{-1}(x_4))
r_{a_2}l_{S^{-1}(x_3)}\\
&=&\le r_{a_1}, x_2\ri \le r_{a_3}, S^{-1}(x_4)\ri
r_{a_2}l_{S^{-1}(x_3)}.
\end{eqnarray*}
Combining these equations, we obtain
\begin{eqnarray*}
\gamma(x)l_yr_a&=&\le l_{y_3}r_{a_1}, x_1\ri
\le l_{y_1}r_{a_3}, S^{-1}(x_3)\ri (l_{y_2}r_{a_2})\gamma(x_2)\\
&=&\le (l_yr_a)_3, x_1 \ri \le (l_yr_a)_1, S^{-1}(x_3) \ri
(l_yr_a)_2 \gamma(x_2).
\end{eqnarray*}
Since   $\pi = m_{H^0} \circ (Id\otimes \gamma)$, the statement
follows from Proposition \ref{piiffgamma}.
\end{proof}

We now apply the results of Section \ref{doubleprojections} to this
Radford biproduct.

\begin{proposition}\prlabel{3.2}
Let $H, X, A, \gamma, \pi$ be as above. Then
$D(H_{X, A}^{\rm cop}, X)\cong B \times H_{X, A}^{\rm cop}$ where $B$
is a Hopf algebra in the category
${}_{H_{X, A}^{\rm cop}}^{H_{X, A}^{\rm cop}}\cal{YD}$. In addition,
\begin{itemize}
\item[(i)]
$B=\{l_{S^{-2}(x_2)}\ot x_1\mid x\in X\}$ and is isomorphic to $X$ as a $k$-space;
\item[(ii)]
$X$ is a left $H_{X, A}^{\rm cop}$ Yetter-Drinfeld module with the structure
\begin{eqnarray}
&&l_xr_a\tr y=\le l_xr_a, \smi (y_1)S^{-2}(y_3)\ri y_2,\label{mhs}\\
&&\l _{\un{X}}(y)=l_{S^{-1}(y_1)S^{-2}(y_3)}\ot y_2.\label{chs}
\end{eqnarray}
\end{itemize}
\end{proposition}
\begin{proof}
Statement (i) follows directly from Proposition
\ref{thetabijection}, while statement (ii) follows from
Proposition \ref{leftleftYD}.
\end{proof}

From now on, $\un{X}$ will be the $k$-vector space $X$, with the
structure of Hopf algebra in the braided category
${}_{H_{X, A}^{\rm cop}}^{H_{X, A}^{\rm cop}}{\cal YD}$ induced from $B$ via the above isomorphism.
Next we compute the structure maps of $\un{X}$ in this
category.

\begin{theorem}\thlabel{1.8}
The structure of $\un {X}$ as a Hopf algebra in
${}_{H_{X, A}^{\rm cop}}^{H_{X, A}^{\rm cop}}{\cal YD}$
is given by the formulas
\begin{eqnarray}
x\circ y&=&\sigma (S(x_3)S^2(x_1), y_2)x_2y_1,\\
\un{\Delta}(x)&=&\sigma (S(x_1)x_3, S(x_4)x_6)x_2\ot x_5,\\
\un{S}(x)&=&\sigma (S^2(x_3)S(x_1), x_4)S(x_2),
\end{eqnarray}
for all $x, y\in X$. $\un{X}$ has the same unit and counit as $X
\subseteq H$.
\end{theorem}
\begin{proof}
Apply the formulas in  Proposition \ref{HopfYD} with $U=H_{X, A}$
and $V = X$.
\end{proof}

\begin{remarks}\reslabel{1.9}
(i) From Equations (\ref{v}) and (\ref{p2}) in Remarks
\ref{cqtremarks}, we obtain another formula for the antipode
$\un{S}$ of $\un{X}$. For we note that
\begin{eqnarray*}\label{sfant}
\un{S}(x)&=& \sigma(S^2(x_3),x_4)
\sigma(S(x_1), x_5) S(x_2) \\
&=&\sigma(S(x_1), x_4)v^{-1}(x_3)S(x_2) \\
&=&\sigma(S(x_1), x_4)v^{-1}(x_2) S^{-1}(x_3)\\
&=& \sigma(S(x_1), x_5)\sigma(S^2(x_2), x_3) S^{-1}(x_4)\\
&=& \sigma(S(x_1), S^{-1}(x_2)x_4) S^{-1}(x_3).
\end{eqnarray*}

(ii) The Hopf algebra isomorphism $\chi :\ \un {X}\times H_{X,
A}^{\rm cop}\ra D(H_{X, A}^{\rm cop}, X)$ is given by $\chi (x\times
l_yr_a)=(l_{S^{-2}(x_2)}\ot x_1)(l_yr_a\ot 1)$, for all $x, y\in X$
and $a\in A$.

(iii) One can easily check that, in general, $\un{X}$ is neither
quantum commutative nor quantum cocommutative as a bialgebra in
${}_{H_{X, A}^{\rm cop}}^{H_{X, A}^{\rm cop}}{\cal YD}$. But, if
 $\sigma$ restricted to $X\ot X$ gives a triangular structure on $X$ (that
is, $\sigma ^{-1}(x, y)= \sigma (y, x)$, for all $x, y\in X$) then
$\un{X}$ is quantum commutative, i.e.
\[
x\circ y=\left(x_{-1}\tr y\right)\circ x_{0},~~\forall~~x, y\in X.
\]

(iv) Any proper Hopf subalgebra $X$ of a CQT Hopf algebra $(H,
\sigma)$ can be viewed as a braided Hopf algebra in (at least) three
different left Yetter-Drinfeld module categories by applying
\prref{3.2} and \thref{1.8}  with different $A$. Specifically,
over $H_{l_X}^{\rm cop}$ (take $A=k$), over the co-opposite of
$H_{l_X}H_{r_X}=H_{r_X}H_{l_X}$ (take $A=X$), and over the
co-opposite of $H_{l_X}H_r=H_rH_{l_X}$ (take $A=H$). Note that
$H_{l_X}\subseteq H_{l_X}H_{r_X}\subseteq H_{l_X}H_{r} $.
\end{remarks}

We note that the solution to the Yang-Baxter equation  from the
Yetter Drinfeld module $\un{X}$ comes from the adjoint coaction.

For if $({\cal B}, \sigma)$ is a CQT bialgebra and ${\mf M}$ a left
$\cal{B}$-comodule then we have seen that
${\mf M}$ is a left Yetter-Drinfeld module over
${\cal B}$ with the initial comodule structure and with the ${\cal B}$-action defined by
(\ref{modstr}). In particular, we get that
$R_{\mf M}\in {\rm End}({\mf M}\ot {\mf M})$ given for all
${\mf m}, {\mf m}'\in {\mf M}$ by
\[
R_{\mf M}({\mf m}\ot{\mf m}')=\sigma ({\mf m}_{-1}, {\mf m}'_{-1}){\mf m}_{0}
\ot {\mf m}'_{0},
\]
is a solution for the quantum Yang-Baxter equation.

In the present setting, with $(H, \sigma)$ CQT, then $X$ is a left
$X$-comodule via $x \mapsto S(x_1)x_3 \ot x_2$. Then $R_{X_{\rm
ad}}\in {\rm End}(X\ot X)$ given for all $x, y\in X$ by
\[
R_{X_{\rm ad}}(x\ot y)=\sigma (S(x_1)x_3, S(y_1)y_3)x_2\ot y_2
\]
is a solution for the quantum Yang-Baxter equation. From the
discussion at the end of  Section \ref{doubleprojections}, writing
$S$ for $S_X$, we note that
\begin{eqnarray*}
R_{\un{X}}(x \ot y)&=&\le \gamma(S^{-1}(y_3)y_1),
S^{-1} (S ^{-1}(x_3)x_1)\ri x_2 \ot y_2\\
&=&\le l_{S^{-1}(S^{-1}(y_3)y_1)}, S^{-1}(S^{-1}(x_3)x_1)\ri
x_2\ot y_2 \\
&=& \sigma( S^{-1}(S^{-1}(x_3)x_1),S^{-1}(S^{-1}(y_3)y_1))
x_2\ot y_2\\
&=&\sigma(S(x_1)x_3, S(y_1)y_3)x_2 \ot y_2 \\
&=&R_{X_{\rm ad}}(x\ot y).
\end{eqnarray*}

Next, we show the connection between our $\un{X}$ and the transmutation theory.

\begin{theorem}\thlabel{unXtransm}
Let $(H, \sigma)$ be a CQT Hopf algebra and $A, X$ sub-Hopf algebras of $H$. Then there is a
braided functor $\mathbf{F}: {\cal M}^X\ra {}_{H_{X, A}^{\rm cop}}^{H_{X, A}^{\rm cop}}{\cal YD}$
such that $\mathbf{F}(\un{\un{X}})=\un{X}$.
\end{theorem}
\begin{proof}
We start by constructing the functor $\mathbf{F}$ explicitly. For
${\mf M}$   a right $X$-comodule, let  $\mathbf{F}({\mf M})={\mf M}$
as  a $k$-vector space, endowed with the following structures:
\begin{eqnarray}
&&l_xr_a\tr {\mf m}=l_xr_a(S^{-2}({\mf m}_{(1)})){\mf
m}_{(0)}~~\mbox{and}~~ \l_{{\mf M}}({\mf m})=l_{S^{-2}({\mf
m}_{(1)})}\ot {\mf m}_{(0)},\label{specstr}
\end{eqnarray}
for any $l_xr_a\in H_{X, A}$ and ${\mf m}\in {\mf M}$. From the
definitions,  $\mathbf{F}({\mf M})$ is both a left $H_{X, A}^{\rm
cop}$-module and comodule. Actually, $\mathbf{F}({\mf M})$ is a
Yetter-Drinfeld module over $H_{X, A}^{\rm cop}$. Here, using the
co-opposite of the structure maps in \prref{2.1}, we see that
relation (\ref{lyd1}) is
\begin{eqnarray*}
&&\hspace*{-1cm}
\le l_{x_1}r_{a_2}, S^{-2} ({\mf m}_{(1)})\ri
l_{x_2}r_{a_1}l_{S^{-2}({\mf m}_{(2)})}\ot {\mf m}_{(0)}\\
&&\hspace*{2cm}
=\le l_{x_2}r_{a_1}, S^{-2}({\mf m}_{(2)})\ri
l_{S^{-2}({\mf m}_{(1)})}l_{x_1}r_{a_2}\ot {\mf m}_{(0)}
\end{eqnarray*}
and this holds since
\begin{eqnarray*}
\le l_{x_1}r_{a_2},y_1 \ri l_{x_2}r_{a_1}l_{y_2}&=&\sigma(y_1, x_1)\sigma (a_2, y_2)l_{x_2}r_{a_1}l_{y_3}\\
&{{(\ref{f2})}\atop =}&\sigma (y_1, x_1)\sigma (a_4, y_2)\sigma (S^{-1}(a_3), y_3)\sigma (a_1, y_5)
l_{x_2}l_{y_4}r_{a_2}\\
&{{(\ref{cqt2})}\atop =}&\sigma(y_1, x_1)\sigma (a_1, y_3)l_{y_2x_2}r_{a_2}\\
&{{(\ref{cqt4})}\atop =}&\sigma(y_2, x_2)\sigma(a_1, y_3)l_{x_1y_1}r_{a_2}\\
&=&\le l_{x_2}r_{a_1},  y_2\ri l_{y_1}l_{x_1}r_{a_2},
\end{eqnarray*}
for all $x, y\in X$ and $a\in A$. So $\mathbf{F}$ is a well defined functor from
${\cal M}^X$ to ${}_{H_{X, A}^{\rm cop}}^{H_{X, A}^{\rm cop}}{\cal YD}$.
($\mathbf{F}$ sends a morphism to itself.)

We claim that $\mathbf{F}$ has a monoidal structure defined by $\varphi_0=Id : k\ra \mathbf{F}(k)$ and
$\varphi_2({\mf M}, {\mf N}):\ \mathbf{F}({\mf M})\ot \mathbf{F}({\mf N})\to
\mathbf{F}({\mf M}\ot {\mf N})$, the family of natural isomorphisms in
${}_{H_{X, A}^{\rm cop}}^{H_{X, A}^{\rm cop}}{\cal YD}$ given by
\[
\varphi_2({\mf M}, {\mf N})({\mf m}\ot {\mf n})=\sigma^{-1}({\mf m}_{(1)}, {\mf n}_{(1)})
{\mf m}_{(0)}\ot {\mf n}_{(0)},
\]
for all ${\mf m}\in {\mf M}\in {\cal M}^X$ and ${\mf n}\in {\mf N}\in {\cal M}^X$
(for more details see for example \cite[XI.4]{k}).
To this end, observe that $\varphi_2({\mf M}, {\mf N})$ is left
$H_{X, A}^{\rm cop}$-linear since
\begin{eqnarray*}
&& \hspace{-3cm} \varphi_2({\mf M}, {\mf N})(l_xr_a\tr ({\mf m}\ot {\mf n}))\\
&{{(\ref{lyd2})}\atop =}&
\varphi_2({\mf M}, {\mf N})(l_{x_2}r_{a_1}\tr {\mf m}\ot l_{x_1}r_{a_2}\tr {\mf n})\\
&{{(\ref{specstr})}\atop =}&  \le l_xr_a, S^{-2}({\mf n}_{(1)}{\mf
m}_{(1)})\ri
\varphi_2({\mf M}, {\mf N})({\mf m}_{(0)}\ot {\mf n}_{(0)})\\
&=& \le l_xr_a, S^{-2}({\mf n}_{(2)}{\mf m}_{(2)})\ri
\sigma^{-1}({\mf m}_{(1)}, {\mf n}_{(1)})
{\mf m}_{(0)}\ot {\mf n}_{(0)}\\
&{{(\ref{cqt4})}\atop =}& \le l_xr_a, S^{-2}({\mf m}_{(1)}{\mf
n}_{(1)})\ri \sigma^{-1}({\mf m}_{(2)}, {\mf n}_{(2)})
{\mf m}_{(0)}\ot {\mf n}_{(0)}\\
&{{(\ref{rightmonoidal}, \ref{specstr})}\atop =}&
\sigma^{-1}({\mf m}_{(1)}, {\mf n}_{(1)})l_xr_a\tr ({\mf m}_{(0)}\ot {\mf n}_{(0)})\\
&=&l_xr_a\tr \varphi_2({\mf M}, {\mf N})({\mf m}\ot {\mf n}),
\end{eqnarray*}
for all $x\in X$, $a\in A$, ${\mf m}\in {\mf M}$ and ${\mf n}\in
{\mf N}$. The fact that $\varphi_2({\mf M}, {\mf N})$ is left $H_{X,
A}^{\rm cop}$-colinear can be proved in a similar manner, the
details are left to the reader.

Clearly, $\varphi_2({\mf M}, {\mf N})$ is bijective and $\varphi_0$ makes the
left and right unit constraints in ${\cal M}^X$ and
${}_{H_{X, A}^{\rm cop}}^{H_{X, A}^{\rm cop}}{\cal YD}$ compatible. So it remains to prove that
$\varphi _2$ respects the associativity constraints of the two categories above. It is not hard to see
that this fact is equivalent to
\[
\sigma(h_1, g_1)\sigma (h_2g_2, h')=\sigma (g_1, h'_1)\sigma(h, h'_2g_2),
\]
for all $h, h', g\in H$, i.e., $\sigma$ is a $2$-cocycle,  a well
known fact which follows by applying the properties of
coquasitriangularity  or see \cite{doi}.

Moreover, we have $\mathbf{F}$ a braided functor, this means
\[
\mathbf{F}(c_{{\mf M}, {\mf N}})\circ \varphi_2({\mf M}, {\mf N})=
\varphi_2({\mf N}, {\mf M})\circ {\bf c}_{\mathbf{F}({\mf M}), \mathbf{F}({\mf N})},
\]
for all ${\mf M}, {\mf N}\in {\cal M}^X$. Indeed, on one hand, by (\ref{rcombr}) we have
\[
\mathbf{F}(c_{{\mf M}, {\mf N}})\circ \varphi_2({\mf M}, {\mf N})({\mf m}\ot {\mf n})
={\mf n}\ot {\mf m}.
\]
On the other hand, by (\ref{lyd3}) and (\ref{specstr}) we compute:
\begin{eqnarray*}
\varphi_2({\mf N}, {\mf M})\circ {\bf c}_{\mathbf{F}({\mf M}), \mathbf{F}({\mf N})}
({\mf m}\ot {\mf n})&=&\varphi_2({\mf N}, {\mf M})
(l_{S^{-2}({\mf m}_{(1)})}\tr n\ot {\mf m}_{(0)})\\
&=&\sigma(S^{-2}({\mf n}_{(1)}), S^{-2}({\mf m}_{(1)})\varphi_2({\mf N}, {\mf M})
({\mf n}_{(0)}\ot {\mf m}_{(0)})\\
&=&{\mf n}\ot {\mf m},
\end{eqnarray*}
as needed.

Finally, by \prref{4.2} and \prref{3.2} it follows that
$\mathbf{F}(\un{\un{X}})=\un{X}$, as objects
in ${}_{H_{X, A}^{\rm cop}}^{H_{X, A}^{\rm cop}}{\cal YD}$.
Furthermore, since $\mathbf{F}$ is a braided functor and $\un{\un{X}}$ is a braided Hopf algebra
in ${\cal M}^X$ it follows that $\mathbf{F}(\un{\un{X}})$ is a braided Hopf algebra in
${}_{H_{X, A}^{\rm cop}}^{H_{X, A}^{\rm cop}}{\cal YD}$ with multiplication given by
\[
\un{m}_{\mathbf{F}(\un{\un{X}})}:\ \mathbf{F}(\un{\un{X}})\ot
\mathbf{F}(\un{\un{X}}) \pile{\rTo^{\varphi_2(\un{\un{X}},
\un{\un{X}})}} \mathbf{F}(\un{\un{X}}\ot \un{\un{X}})
\pile{\rTo^{\mathbf{F}(\un{m}_{\un{\un{X}}})}}
\mathbf{F}(\un{\un{X}}),
\]
the comultiplication defined by
\[
\un{\Delta}_{\mathbf{F}({\un{\un{X}}})}:\ \mathbf{F}(\un{\un{X}})
\pile{\rTo^{\mathbf{F}(\un{\Delta}_{\un{\un{X}}})}}
\mathbf{F}(\un{\un{X}}\ot \un{\un{X}})
\pile{\rTo^{\varphi_2^{-1}(\un{\un{X}}, \un{\un{X}})}}
\mathbf{F}(\un{\un{X}})\ot \mathbf{F}(\un{\un{X}}),
\]
and the same antipode as those of $\un{\un{X}}$. More precisely, according to
the proof of \prref{4.2} the multiplication $\bullet$ is given by
\begin{eqnarray*}
x\bullet y&=&\sigma^{-1}(x_{(1)}, y_{(1)})x_{(0)}\cdot y_{(0)}\\
&=&\sigma^{-1}(x_{(1)}, S(y_1)y_3)x_{(0)}\cdot y_2\\
&=&\sigma(x_{(2)}, S^{-1}(y_4)y_1)\sigma(x_{(1)}, S(y_2))x_{(0)}y_3\\
&{{(\ref{cqt2})}\atop =}&\sigma(x_{(1)}, S^{-1}(y_2))x_{(0)}y_1=
\sigma(S(x_3)S^2(x_1), y_2)x_2y_1,
\end{eqnarray*}
for all $x, y\in X$. Similarly, we have
\begin{eqnarray*}
\un{\Delta}_{\mathbf{F}({\un{\un{X}}})}(x)&=&\sigma((x_1)_{(1)},
(x_2)_{(1)})(x_1)_{(0)}\ot
(x_2)_{(0)}\\
&=&\sigma(S(x_1)x_3, S(x_4)x_6)x_2\ot x_5,
\end{eqnarray*}
for all $x\in X$. Comparing to the structures in \thref{1.8} we conclude that
$\mathbf{F}(\un{\un{X}})=\un{X}$ as braided Hopf algebras in
${}_{H_{X, A}^{\rm cop}}^{H_{X, A}^{\rm cop}}{\cal YD}$, and this finishes the proof.
\end{proof}

We end this section with few comments.

\begin{remarks}
Let $(H, \sigma)$, $X$ and $A$ be as above.

(i) One can easily see that the map $\gamma: A\ra H_{r_A}^{\rm
cop}\subseteq H_{X, A}^{\rm cop}$ given by $\gamma(a)=r_a$ is a Hopf
algebra morphism. Moreover, for this map $\gamma$, (\ref{gamma}) reduces to
\[
r_ar_b=\sigma^{-1}(b_3, a_3)\sigma(b_1, a_1)r_{b_2}r_{a_2}
\]
which holds by (\ref{cqt4}). So on the $k$-vector space $A$ we have
a braided Hopf algebra structure, denoted by $\un{A}$,  in
${}_{H_{X, A}^{\rm cop}}^{H_{X, A}^{\rm cop}}{\cal YD}$. Namely,
$\un{A}$ is a left Yetter-Drinfeld module over $H_{X, A}^{\rm cop}$
via
\[
l_xr_a\tr b=\le l_xr_a, S^{-1}(b_1)S^{-2}(b_3)\ri b_2,~~
\lambda_{\un{A}}(a)=r_{S^{-1}(a_3)a_1}\ot a_2,
\]
and a Hopf algebra with the same unit and counit as $A$ and
\begin{eqnarray*}
&&\hspace*{5mm}
a\cdot b=\sigma (b_2, S(a_1)a_3)a_2b_1=
\sigma_{21}^{-1}(S(a_3)S^2(a_1), b_2)a_2b_1,\\
&&\hspace*{5mm}
\un{\Delta}(a)=\sigma(S(a_4)a_6, S^{-1}(a_3)a_1)a_2\ot a_5=
\sigma_{21}^{-1}(S(a_1)a_3, S(a_4)a_6)a_2\ot a_5,\\
&&\hspace*{5mm}
\un{S}(a)=\sigma(a_4, a_1S(a_3))S(a_2)=
\sigma_{21}^{-1}(S^2(a_3)S(a_1), a_4)S(a_2).
\end{eqnarray*}
Comparing with the structure of $\un{X}$ we conclude that $\un{A}$ can be obtained
from $\un{X}$ by replacing $(H, \sigma)$ with $(H, \sigma^{-1}_{21})$ and
interchanging $X$ and $A$. For this, observe that if
$\tilde{l}_a$ and $\tilde{r}_x$ are the elements of $H^0$ corresponding to
$\tilde{H}=(H, \sigma_{21}^{-1})$ then
$\tilde{r}_x=l_{S^{-1}(x)}$ and $\tilde{l}_a=r_{S(a)}$, so
$\tilde{H}_{A, X}:=\tilde{H}_{\tilde{l}_A}\tilde{H}_{\tilde{r}_X}=
H_{r_A}H_{l_X}=H_{X, A}$, as Hopf algebras.
The remaining details are now trivial.

(ii) Although it might seem  more natural to try to obtain a braided
Hopf algebra structure in ${}_{H_{X, A}^{\rm cop}}^{H_{X, A}^{\rm
cop}}{\cal YD}$ on $XA=AX$, for this we would  need a Hopf algebra
map $\pi: AX\ra H_{X, A}^{\rm cop}$ satisfying condition
(\ref{gamma}) in Proposition \ref{piiffgamma}. A candidate is $\pi
(ax)=r_al_{S^{-1}(x)}$ but, in general, it is not well defined.
(Take for example $X=A=H=({\rm SL}_q(N), \sigma)$, the CQT Hopf
algebra defined in the last section of this paper.)
\end{remarks}
%%%%%%%%%%%%%%%%%%%%%%%%%%%%%%%%%%%%%%%%%%%%%%%%%%%%%%%%%%%%%%%%%
\section{The finite dimensional case}\label{fdcase}
%%%%%%%%%%%%%%%%%%%%%%%%%%%%%%%%%%%%%%%%%%%%%%%%%%%%%%%%%%%%%%%%%%%%%
\setcounter{equation}{0} In this section we discuss how our results
thus far relate to those of Radford \cite{rad} for finite
dimensional Hopf algebras.

Throughout this section, $(H, R)$ is a finite dimensional
quasitriangular ($QT$ for short) Hopf algebra, so its dual linear
space $H^*$ has a CQT structure given by $\sigma :\ H^*\ot H^*\ra
k$, $\sigma (p, q)=p(R^1)q(R^2)$, for all $p, q\in H^*$, where
$R:=R^1\ot R^2\in H\ot H$. Thus, we can consider the sub-Hopf
algebras $H^*_l$ and $H^*_r$ of $H^{**}$ defined in the previous
section. Identifying $H^{**}$ and $H$ via the canonical Hopf algebra
isomorphism
\[
\Theta :\ H\ra H^{**},~~\Theta (h)(h^*)=h^*(h),
\]
for all $h\in H$ and $h^*\in H^*$, we will prove that the sub-Hopf
algebras $H^*_l$ and $H^*_r$ of $H^{**}$ can be identified with
the sub-Hopf algebras $R_{(l)}$ and $R_{(r)}$ of $H$ constructed in
\cite{rad}. Recall that
\[
R_{(l)}:=\{q(R^2)R^1\mid q\in H^*\}~~{\rm and}~~
R_{(r)}:=\{p(R^1)R^2\mid p\in H^*\}.
\]
Also, note that, if we write $R=\sum \limits _{i=1}^mu_i\ot v_i\in
H\ot H$, with $m$ as small as possible, then $\{u_1, \cdots , u_m\}$
is a basis for $R_{(l)}$ and $\{v_1, \cdots , v_m\}$ is a basis for
$R_{(r)}$, respectively. We extend $\{u_1, \cdots, u_m\}$ to a basis
$\{u_1, \cdots , u_m, \cdots , u_n\}$ of $H$ and denote by
$\{u^i\}_{i=\ov{1, n}}$ the dual basis of $H^*$ corresponding to
$\{u_i\}_{i=\ov{1, n}}$. Similarly, we extend $\{v_1, \cdots ,
v_m\}$ to a basis $\{v_1, \cdots , v_m, \cdots , v_n\}$ of $H$ and
denote by $\{v^i\}_{i=\ov{1, n}}$ the dual basis of $H^*$
corresponding to $\{v_i\}_{i=\ov{1, n}}$.

\begin{lemma}\lelabel{3.1}
For the above context the following statements hold:
\begin{itemize}
\item[i)] The map $\mu _l : H^*_l\ra R_{(l)}$ defined by
$\mu _l(l_q)=q(R^2)R^1$, for all $q\in H^*$, is well defined and a
Hopf algebra isomorphism. Its inverse is given by $\mu
_l^{-1}(u_i)=l_{v^i}$, for all $i=\ov{1, m}$. In particular,
$\{l_{v^i}\}_{i=\ov{1, m}}$ is a basis for $H^*_l$.
\item[ii)] The map $\mu _r: H^*_r\ra R_{(r)}$ given by
$\mu _r(r_p)=p(R^1)R^2$, for all $p\in H^*$, is well defined and a
Hopf algebra isomorphism. Its inverse is defined by $\mu
_r^{-1}(v_i)=r_{u^i}$, for all $i=\ov{1, m}$. In particular,
$\{r_{u^i}\}_{i=\ov{1, m}}$ is a basis for $H^*_r$.
\end{itemize}
\end{lemma}
\begin{proof}
We show only i), with the proof of ii) being similar. First observe
that
\[
l_q(h^*)=\sigma (h^*, q)=h^*(R^1)q(R^2)=\Theta(q(R^2)R^1)(h^*),
\]
for all $h^*\in H^*$. Thus $H^*_l=\{\Theta(q(R^2)R^1)\mid q\in
H^*\}=\Theta (R_{(l)})\cong R_{(l)}$. Clearly $\mu _l$ is precisely
the restriction and corestriction of $\Theta ^{-1}$ at
$\Theta(R_{(l)})$ and $R_{(l)}$, respectively. Hence $\mu_l$ is well
defined and a Hopf algebra isomorphism.

Finally, for all $h^*\in H^*$ we have
\[
l_{v^i}(h^*)=\sigma (h^*, v^i)=h^*(R^1)v^i(R^2) =\sum \limits
_{j=1}^mh^*(u_j)v^i(v_j),
\]
and therefore $l_{v^i}=0$, for all $i\in \{m + 1, \cdots , n\}$, and
$l_{v^i}(h^*)=h^*(u_i)=\Theta(u_i)(h^*)$, i.e. $l_{v^i}=\Theta
(u_i)$, for all $i=\ov{1, m}$. This shows that
$\mu_l^{-1}(u_i)=l_{v^i}$, for all $i=\ov{1, m}$, so the proof is
complete.
\end{proof}

From \prref{2.1}, $\sigma$ gives a   pairing on  $H_r^{* \rm cop}
\ot H^*_l$, and we   consider the generalized quantum double
$D(H^*_r, H^*_l)$. On the other hand, $H^*_r$ and $H^*_l$ are finite
dimensional Hopf algebras, so we can construct the usual Drinfeld
quantum doubles $D(H^*_l)$ and $D(H^*_r)$. To demonstrate the
connections between these Hopf algebras, we first need the following
lemma.

\begin{lemma}\lelabel{3.2}
Suppose that $\le , \ri : U\ot V\ra k$ defines a duality between the
finite dimensional Hopf algebras $U$ and $V$ such that there exists
a Hopf algebra isomorphism $\Psi :V^*\ra U$ with the property
\begin{equation}\label{f4}
\le \Psi (v^*), v\ri =v^*(v),~~\forall~~v^*\in V^*,~~v\in V.
\end{equation}
Then $D(U^{\rm cop}, V)\cong D(V)\cong D(U^{\rm op, cop})^{\rm op}$
as Hopf algebras.
\end{lemma}
\begin{proof}
As coalgebras
\begin{eqnarray*}
D(U^{\rm cop}, V)=U^{\rm cop}\ot V &{{\Psi ^{-1}\ot Id}\atop \cong}&
V^{*\rm cop}\ot V=D(V).
\end{eqnarray*}
Thus will be enough to check that $\Psi \ot Id: D(V)\ra D(U^{\rm
cop}, V)$ is an algebra morphism. For all $x^*, y^*\in V^*$ and $x,
y\in V$ we compute
\begin{eqnarray*}
&&\hspace*{-1.5cm}
(\Psi \ot Id)\left((x^*\ot x)(y^*\ot y)\right)\\
&=&y^*_1(\smi(x_3))y^*_3(x_1)\Psi(x^*y^*_2)\ot x_2y\\
&=&\le \Psi(y^*_1), \smi (x_3)\ri \le \Psi(y^*_3), x_1\ri
\Psi(x^*)\Psi(y^*_2)\ot x_2y\\
&=&\left(\Psi(x^*)\ot x\right)\left(\Psi(y^*)\ot y\right) =(\Psi\ot
Id)\left(x^*\ot x\right)(\Psi \ot Id)\left(y^*\ot y\right),
\end{eqnarray*}
as needed. Clearly, $\Psi \ot Id$ respects the units, so $\Psi \ot
Id: D(U^{\rm cop}, V)\ra D(V)$ is a Hopf algebra isomorphism.

From \cite[Theorem 3]{rad} we know that $D(U^{\rm op, cop})^{\rm
op}\cong D(U^*)$ as Hopf algebras. But $U^*\cong V^{**}\cong V$ as
Hopf algebras, and therefore $D(U^*)\cong D(V)$ as Hopf algebras. We
conclude that $D(U^{\rm op, cop})^{\rm op}\cong D(V)$ as Hopf
algebras, and this finishes the proof.
\end{proof}

\begin{proposition}\prlabel{3.3}
Let $(H, R)$ be a finite dimensional QT Hopf algebra with an
$R$-matrix $R=\sum \limits_{i=1}^mu_i\ot v_i\in H\ot H$, with $m$ as
small as possible. If $H^*_l$ and $H^*_r$ are the Hopf subalgebras
of $H^{**}\cong H$ defined above then
\[
D(H^*_r, H^*_l)\cong D(H^*_l)\cong D((H^*_r)^{\rm op})^{\rm op},
\]
as Hopf algebras.
\end{proposition}
\begin{proof}
We apply \leref{3.2} for $U=(H^*_r)^{\rm cop}$ and $V=H^*_l$. So we
shall prove that there exists a Hopf algebra morphism $\Psi :
(H^*_l)^*\ra (H^*_r)^{\rm cop}$ (or, equivalently, from
$(H^*_l)^{*\rm cop}$ to $H^*_r$) satisfying (\ref{f4}), this means
$\le \Psi (\eta), l_q\ri =\eta(q)$, for all $\eta \in (H^*_l)^*$ and
$q\in H^*$. We will use the Hopf algebra isomorphism $\xi :
R_{(l)}^{*\rm cop}\ra R_{(r)}$ from \cite[Proposition 2(c)]{rad}
defined by  $\xi (\zeta)=\zeta(R^1)R^2$, for all $\zeta \in
R_{(l)}^*$. Note that $\xi ^{-1}(v_i)=\tilde{u}^i$, for all
$i=\ov{1, m}$, where $\tilde{u}^i$ is the restriction of $u^i$ at
$R_{(l)}$.

Now, define $\Psi: (H^*_l)^{*\rm cop}\ra H^*_r$ as the composition
of the following Hopf algebra isomorphisms
\[
\Psi : (H^*_l)^{*\rm cop}
\pile{\rTo^{(\mu_l^{-1})^*}}
R_{(l)}^{*\rm cop}
\pile{\rTo^{\xi}}
R_{(r)}
\pile{\rTo^{\mu_r^{-1}}}
H^*_r.
\]
Explicitly, we have $\Psi (\eta)=\sum \limits_{i=1}^m\eta
(l_{v^i})r_{u^i}$, for all $\eta \in (H^*_l)^*$, and this allows us
to compute $ \le \Psi (\eta), l_{v^s}\ri  $  to be
\[
\sum \limits _{i=1}^m \eta (l_{v^i})\le r_{u^i}, l_{v^s}\ri = \sum
\limits_{i}^m \eta (l_{v^i})\sigma(u^i, v^s) =\sum \limits _{i,
j=1}^m \eta (l_{v^i}) u^i(u_j)v^s(v_j) =\eta(l_{v^s}),
\]
for all $s=\ov{1, m}$. Since $\{l_{v^s}\mid s=\ov{1, m}\}$ is a
basis for $H^*_l$ it follows that $\le \Psi (\eta), l_q\ri
=\eta(q)$, for all $\eta \in (H^*_l)^*$ and $q\in H^*$, so the proof
is finished.
\end{proof}

\begin{remark}\relabel{3.4}
Let $U$ and $V$ be two Hopf algebras in duality via the bilinear
form $\le , \ri$. If there is an element $\sum \limits
_{i=1}^tu_i\ot v_i\in U\ot V$ such that, for all $u\in U$ and $v\in
V$,
\begin{equation}\label{f5}
\sum \limits _{i=1}^t\le u_i, v\ri v_i=v~~{\rm and}~~ \sum \limits
_{i=1}^t\le u, v_i\ri u_i=u,
\end{equation}
then $R=\sum \limits _{i=1}^t(1\ot v_i)\ot (u_i\ot 1)$ is an
$R$-matrix for $D(U^{\rm cop}, V)$.

Let now $(H, R)$ be a finite dimensional QT Hopf algebra with
$R=\sum \limits_{i=1}^mu_i\ot v_i$, with $m$ as small as possible.
Then $\sum \limits _{i=1}^m r_{u^i}\ot l_{v^i}\in H^*_r\ot H^*_l$
satisfies the conditions in (\ref{f5}). Indeed, for all $p\in H^*$
we have
\begin{eqnarray*}
&&\hspace*{-1cm}
\sum \limits _{i=1}^m\le r_p, l_{v^i}\ri r_{u^i}=
\sum \limits _{i, j=1}^mp(u_j)v^i(v_j)r_{u^i}
=\sum \limits_{i=1}^mp(u_i)\Theta (v_i)=\Theta(p(R^1)R^2)=r_p,
\end{eqnarray*}
and in a similar manner one can verify that $\sum \limits_{i=1}^m\le
r_{u^i}, l_q\ri l_{v^i}=l_q$, for all $q\in H^*$. Therefore, ${\mf
R}=\sum \limits_{i=1}^m (1\ot l_{v^i})\ot (r_{u^i}\ot 1)$ is an
$R$-matrix for $D(H^*_r, H^*_l)$.

On the other hand, a basis on $(H^*_l)^*$ can be obtained by using
the inverse of the Hopf algebra isomorphism $\Psi$ constructed in
\prref{3.3}. Since $\Psi ^{-1}(r_{u^i})=\tilde{u}^i\circ \mu_l$ it
follows that $\{\tilde{u}^i\circ \mu _l\mid i=\ov{1, m}\}$ is a
basis of $(H^*_l)^*$. Moreover, we can easily see that it is the
basis of $(H^*_l)^*$ dual to the basis $\{l_{v^i}\mid i=\ov{1, m}\}$
of $H^*_l$, so ${\cal R}=\sum \limits _{i=1}^m(1\ot
l_{v^i})\ot (\tilde{u}^i\circ \mu_l\ot 1)$ is an $R$-matrix for
$D(H^*_l)$. Thus, we can conclude that the Hopf algebra
isomorphism $\Psi \ot Id: D(H^*_l)\ra D(H^*_r, H^*_l)$ constructed
in \prref{3.3} is actually a QT Hopf algebra isomorphism, that is,
in addition, $\left((\Psi\ot Id)\ot (\Psi \ot Id)\right)({\cal R})={\mf R}$.
The verification of the details is left to the reader.
\end{remark}

We are now able to prove that, for our context, the Hopf algebra
surjection from \prref{1.2} is in fact a QT morphism, and that it
can be deduced from the surjective QT morphism $F: D(R_{(l)})\ra
H_R:=R_{(l)}R_{(r)}=R_{(r)}R_{(l)}$ considered in \cite[Theorem
2]{rad}. Namely, $F(\zeta \ot \hbar)=\xi(\zeta)\hbar$, for all
$\zeta\in R_{(l)}^*$ and $\hbar \in R_{(l)}$, where $\xi
(\zeta)=\zeta(R^1)R^2$ is the Hopf algebra isomorphism from
$R_{(l)}^{*\rm cop}$ to $R_{(r)}$ considered in \prref{3.3}.

\begin{proposition}\prlabel{3.5}
Under the above circumstances and notations we have the following
diagram commutative.
\begin{center}
\begin{diagram}
D(H^*_r, H^*_l)&         &\rTo^{f}&          &H^*_{\sigma}:=H^*_rH^*_l=H^*_lH^*_r\\
\dTo^{\Psi ^{-1}\ot Id}_{\cong}& &        &          &\uTo_{\mu_r^{-1}\mu_l^{-1}}^{\cong}\\
D(H^*_l)&\rTo^{(\mu_l^{-1})^*\ot \mu_l}_{\cong}& D(R_{(l)})&
\rTo^{F}&H_R:=R_{(r)}R_{(l)}=R_{(l)}R_{(r)}
\end{diagram}
\end{center}
\vspace{0.5cm} Here $(\mu _l^{-1})^*$ is the transpose of $\mu
_l^{-1}$ and $\mu_r^{-1}\mu_l^{-1}$ is defined by
\[
\mu_r^{-1}\mu_l^{-1}(p(r^1)r^2q(R^2)R^1)=r_pl_q,~~ \forall~~ p, q\in
H^*,
\]
where $R=R^1\ot R^2=r^1\ot r^2$ is the $R$ matrix of $H$.
\end{proposition}
\begin{proof}
Straightforward. We only note that $p(r^1)r^2q(R^2)R^1=
p'(r^1)r^2q'(R^2)R^1$ iff
$\Theta^{-1}(r_p)\Theta^{-1}(l_q)=\Theta^{-1}(r_{p'})\Theta^{-1}(l_{q'})$,
iff $\Theta^{-1}(r_pl_q)=\Theta^{-1}(r_{p'}l_{q'})$, iff
$r_pl_q=r_{p'}l_{q'}$, so $\mu_r^{-1}\mu_l^{-1}$ is well defined.
\end{proof}

\begin{corollary}\colabel{2.6}
$H^*_{\sigma}$ is QT with $\Re=\sum \limits _{i=1}^ml_{v^i}\ot
r_{u^i}$ and in this way $f$ becomes a surjective QT Hopf algebra
morphism.
\end{corollary}

Finally, we are able to show that in the finite dimensional case the
morphism $\pi$ in \prref{1.5} arises from the particular situation
described above.

\begin{corollary}\colabel{2.8}
Let $(H, R)$ be a finite dimensional $QT$ Hopf algbera with an $R$
matrix $R=\sum \limits_{i=1}^mu_i\ot v_i$, $m$ as small as possible.
For any Hopf subalgebra $X$ of $H^*$ there exists a surjective Hopf
algebra morphism $\pi:\ D(H^{*\rm cop}_{\sigma}, X)\rightarrow
H^{*\rm cop}_{\sigma}$ covering the canonical inclusion
$i_{H^*_{\sigma}}:\ H^{*\rm cop}_{\sigma}\rightarrow D(H^{*\rm
cop}_{\sigma}, X)$.
\end{corollary}
\begin{proof}
From \coref{2.6} and \exref{2.3},   such a morphism $\pi$ exists.
Moreover, using the definition of $\pi$ in \exref{2.3}, for all $p,
q\in H^*$ and $h^*\in X$, we have
\begin{eqnarray*}
\pi (l_pr_q\ot h^*)&=&\sum \limits_{i=1}^m\le l_{u^i}, h^*\ri
l_pr_ql_{S^{* -1}(v^i)}
=\sum \limits_{i, j=1}^mu^i(u_j)h^*(v_j)l_pr_q l_{S^{* -1}(v^i)}\\
&=&\sum \limits_{i=1}^mh^*(v_i)l_pr_ql_{S^{* -1}(v^i)} =\sum
\limits_{i=1}^nh^*(v_i)l_pr_ql_{S^{* -1}(v^i)}\\& =&l_pr_ql_{S^{*
-1}(h^*)},
\end{eqnarray*}
where $S^*$ is the antipode of $ H^*$ as usual.
\end{proof}

\begin{remarks}\reslabel{5.8}
(i) Here, the sub-Hopf algebra $H^*_\sigma$ of $H^{**} \cong H$ is
the smallest sub-Hopf algebra $C$ of $H^{**}$ such that ${\mf
R}:=(\Theta\ot \Theta)(R)\in C \ot C$ and ${\mf R}$ played a crucial
role in the definition of $\pi$. Now, for $H$ CQT not necessarily 
finite dimensional with sub-Hopf 
algebras $X$ and $A$ as in \prref{1.5}, we cannot 
ensure that $H_{X,A}$ has a quasitriangular structure but we can
define the projection $\pi$ in the same way.

(ii) Applying \coref{2.8} with $X=H^*$ we see that $H^*$ has a
braided Hopf algebra structure  in the braided category of left
Yetter-Drinfeld modules over $H^{*\rm cop}_{\sigma}$. Identifying
$H^{*\rm cop}_{\sigma}$ with $H_R^{\rm cop}$ we obtain that $H^*$ is
a braided Hopf algebra in ${}_{H^{\rm cop}_R}^{H^{\rm cop}_R}{\cal
YD}$, denoted in what follows by $\un{H^*}$. From \prref{3.2} the
structure of $\un{H^*}$ as left $H^{\rm cop}_R$ Yetter-Drinfeld
module is given by
\begin{eqnarray*}
&&\hbar \cdot \v =(S^{* -1} (\v _1)S^{* -2}(\v _3))(\hbar)\v _2=
S^{-2}(\hbar _2)\rh \v \lh \smi(\hbar _1),\label{ydhs1}\\
&&\l _{\un{H^*}}(\v)= (S^{* -1} (\v _1)S^{* -2}(\v _3))(R^2)R^1\ot
\v _2 =R^1\ot R^2\cdot \v ,\label{ydhs2}
\end{eqnarray*}
for all $\hbar \in H_R$ and $\v \in H^*$, and from \thref{1.8} the
structure of $\un{H^*}$ as a bialgebra in ${}_{H^{\rm
cop}_R}^{H^{\rm cop}_R}{\cal YD}$ is the following. The
multiplication is defined by
\begin{eqnarray*}
&&\hspace*{-1cm}
\v\un{\circ}\psi
=\v _3(S(R^1_1))\v _1(S^2(R^1_2))\psi_2(R^2)\v _2\psi_1\nonumber\\
&&\hspace*{1cm} =S(R^1_1)\rh \v \lh S^2(R^1_2)\ot R^2\rh
\psi,\label{muhs}
\end{eqnarray*}
for all $\v , \psi\in H^*$, the unit is $\va$, the comultiplication
is given by
\begin{eqnarray*}
&&\hspace*{-1cm} \un{\Delta}_{\un{H^*}}(\v)=\v_1(S(R^1_1))\v
_3(R^1_2)\v_4(S(R^2_1))\v _6(R^2_2)
\v_2\ot \v_5\nonumber\\
&&\hspace*{1cm} =R^1_2\rh \v_1\lh S(R^1_1)\ot R^2_2\rh \v _2\lh
S(R^2_1),\label{cmhs}
\end{eqnarray*}
for all $\v \in H^*$, and the counit is $\Theta(1)$.
Finally, according to the part 1) of \resref{1.9},
$\un{H^*}$ is a braided Hopf algbera with antipode
\begin{eqnarray*}
&&\hspace*{-1cm}
\un{S}_{\un{H^*}}(\v)=\v_1(S(R^1))\v_4(R^2_1)
\v_2(\smi(R^2_2))\v_3\circ S^{-1}\nonumber\\
&&\hspace*{2cm}
=(R^2_1\rh \v\lh S(R^1)\smi(R^2_2))
\circ \smi .\label{ahs}
\end{eqnarray*}

Using similar arguments to those in Subsection \ref{tt} one can
easily see that $\un{H^*}$ coincides with $(\un{\un{H^{\rm
cop}}})^{\star}$, the categorical right dual of $\un{\un{H^{\rm
cop}}}$ in ${}_{H^{\rm cop}}{\cal M}$, viewed as a braided Hopf
algebra in ${}_{H^{\rm cop}_R}^{H^{\rm cop}_R}{\cal YD}$ through a
composite of two canonical braided functors. (The structures of
$\un{\un{H^{\rm cop}}}$, the associated enveloping algebra braided
group of $(H^{\rm cop}, R_{21}:=R^2\ot R^1)$, can be obtained from
\cite[Example 9.4.9]{majfoundations}, and then the braided structure
of $(\un{\un{H^{\rm cop}}})^{\star}$ can be deduced from \cite{ag}
or \cite{bcp}.)

More precisely, $\un{H^*}={\mf F}_{(H^{\rm cop}, R_{21})}
{\mf G}_{(H^{\rm cop}, R_{21})}((\un{\un{H^{\rm cop}}})^{\star})$,
as braided Hopf algebras in ${}_{H^{\rm cop}_R}^{H^{\rm cop}_R}{\cal YD}$,
where, in general, if $(H, R)$ is a QT Hopf algebra then
\begin{itemize}
\item[(i)] ${\mf F}_{(H, R)}:\ {}_H{\cal M}\ra {}_H^H{\cal M}$ is the braided
functor which acts as identity on objects and morphisms and for
${\mf M}\in {}_H{\cal M}$, ${\mf F}_H({\mf M})={\mf M}$ as $H$-module,
and has the left $H$-coaction given by $\l _M({\mf m})=R^2\ot R^1\cd {\mf m}$,
for all ${\mf m}\in {\mf M}$;
\item[(ii)] ${\mf G}_{(H, R)}:\ {}_{(H, R)}{\cal M}\ra {}_{(H_R, R)}{\cal M}$
is the functor of restriction of scalars which is braided because
the inclusion $H_R\hookrightarrow H$ is a QT Hopf algebra morphism.
\end{itemize}
The verification of all these details is left to the reader.
\end{remarks}

%%%%%%%%%%%%%%%%%%%%%%%%%%%%%%%%%%%%%%%%%%%%%%%%%%%%%%%%%%%%%%%%%%%%
\section{The Hopf algebras ${\rm SL}_q(N)$ and $U_q^{\rm ext}({\rm sl}_N)$}\selabel{slqn}
%%%%%%%%%%%%%%%%%%%%%%%%%%%%%%%%%%%%%%%%%%%%%%%%%%%%%%%%%%%%%%%%%%%%%
\setcounter{equation}{0} In this section we apply the results in
Section \ref{cqtbialgs} to the CQT Hopf algebra $H= {\rm SL}_q(N)$
introduced in \cite{takGL}.  We will show by explicit computation
that $H_\sigma \cong U_q^{\rm ext}({\rm sl}_N)$. The computation
should be compared to results of \cite[3.2]{doitak} and the
description of $U_q^{\rm ext}({\rm sl}_N)$ from \cite[Ch. 8, Theorem
33]{ks}.

Our approach is to study the image of the map $\omega$ in
\cite[3.2(b)]{doitak}, that is, $H_\sigma$, a sub-Hopf algebra of
$H^0$. If, instead, we studied the image of the map $\theta$ in
\cite[3.2(a)]{doitak}, then we would be considering a sub-bialgebra
of $(H^0)^{\rm cop}$. Constructing the isomorphism between these two
approaches  seems to be no simpler computationally than the direct
calculations we supply below.

To make this section as self-contained as possible, we first recall
the definition of the CQT bialgebra $M_q(N)$, and outline the
construction of ${\rm SL}_q(N)$.

For $V$ a  $k$-vector space with finite basis $\{e_1, \cdots ,
e_N\}$ and any   $0 \neq q\in k$,  we  associate a solution $c$ of
the Yang-Baxter equation, $c: V\ot V\ra V\ot V$, by
\[
c(e_i\ot e_j)=q^{\d _{i, j}}e_j\ot e_i + [i>j](q - q^{-1})e_i\ot
e_j,
\]
for all $1\leq i, j\leq N$, where $\d _{i, j}$ is  Kronecker's symbol and
$[i>j]$ is the Heaviside symbol, that is, $[i>j]=0$ if $i\leq j$ and
$[i>j]=1$ if $i>j$ (see for instance \cite[Proposition VIII
1.4]{k}). By the FRT construction, to any solution $c$ of the
quantum Yang-Baxter equation we can associate a CQT bialgebra,
denoted by $A(c)$, (see \cite{frt, k}), obtained by taking a
quotient of the free algebra generated by the family of
indeterminates $(x_{ij})_{1 \leq i,j \leq N}$. For the map $c$
above, $A(c)$ is denoted by $M_q(N)$ and has the following
structure, cf. \cite[Exercise 10, p. 197]{k}. As an algebra $M_q(N)$
is generated by 1 and by $(x_{i,j})_{1\leq i, j\leq N}$. (Note that
we write $x_{i,j}$ as $x_{ij}$ if the meaning is clear.)
Multiplication in $M_q(N)$ is defined by the following relations:
\begin{eqnarray}
&&x_{im}x_{in}=qx_{in}x_{im},~~\forall~~n < m,\label{mq1}\\
&&x_{jm}x_{im}=qx_{im}x_{jm},~~\forall~~i<j,\label{mq2}\\
&&x_{jn}x_{im}=x_{im}x_{jn},~~\forall~~i<j~~{\rm and}~~n<m,\label{mq3}\\
&&x_{jm}x_{in} - x_{in}x_{jm}=(q - q^{-1}) x_{im}x_{jn},~~
\forall~~i<j~~{\rm and}~~n<m.\label{mq4}
\end{eqnarray}
The coalgebra structure on the $x_{ij}$ is that of a   comatrix
coalgebra, that is
\begin{equation}\label{cmqn}
\Delta (x_{ij})=\sum \limits_{k=1}^Nx_{ik}\ot x_{kj},~~
\va(x_{ij})=\d _{i, j},
\end{equation}
for all $1\leq i, j\leq N$.

A CQT structure on $M_q(N)$ is given by the skew pairing $\sigma' :
M_q(N)\ot M_q(N)\ra k$ defined on generators $x_{im}, x_{jn}$ by
\begin{equation}\label{smqn}
\sigma' (x_{im}, x_{jn})=q^{\d_{i, j}}\d _{m ,i}\d _{n, j} + [j>i](q
- q^{-1})\d _{m, j}\d_{n ,i},
\end{equation}
and satisfying (\ref{cqt1}) to (\ref{cqt4}).  Complete details of
this construction can be found in \cite[Theorem VIII 6.4]{k} or see
\cite[(3.3)]{doitak} and \cite{doi}. This skew pairing is invertible
with inverse obtained in the formulas above by replacing $q$ by
$q^{-1}$.

For $z \in k$ such that $z^N = q^{-1}$, we define another skew
pairing $\sigma$ on $M_q(N)\ot M_q(N)$ by defining
\begin{equation}\label{nsmqn}
\sigma(1,-) = \sigma(-, 1)=\va \mbox{   and   }\sigma(x_{im},
x_{jn}) = z\sigma'(x_{im}, x_{jn})
\end{equation}
and extending. For example,
$\sigma(x_{im}x_{kl}, x_{jn})=z^2\sum\limits_{r=1}^N\sigma'(x_{im}, x_{jr})
\sigma'(x_{kl}, x_{rn})$.
Since $\sigma$ satisfies (\ref{cqt4}) when $h, h'$ are generators
$x_{ij}$, then by \cite[Lemma VIII 6.8]{k}, $\sigma$ satisfies
(\ref{cqt4}), and $\sigma$ gives $M_q(N)$ a CQT structure.
Explicitly, we have that
\begin{eqnarray}
\sigma(x_{ii},x_{ii})&=&zq;\label{zq}\\
\sigma(x_{ii},x_{jj})&=&z~~{\rm for}~~i\neq j; \label{z}\\
\sigma(x_{ij}, x_{ji})&=&z(q - q^{-1})~~{\rm if}~~i < j\label{zbracket};\\
\sigma(x_{ij},x_{st})&=&0~~{\rm otherwise.}\label{otherwise}
\end{eqnarray}
The inverse $\sigma^{-1}$ is obtained by
replacing $q$ by $q^{-1}$ and $z$ by $z^{-1}$.

The bialgebra $M_q(N)$ does not have a Hopf algebra structure, but
it possesses a remarkable grouplike central element
\[
{\rm det}_q=\sum \limits _{p\in S_N}(-q)^{-l(p)}x_{1p(1)}\cdots
x_{Np(N)},
\]
which allows us to construct ${\rm SL}_q(N):= M_q(N)/({\rm det}_q
-1)$. Here $S_N$ denotes the group of permutations of order $N$, and
$l(p)$ the number of the inversions of $p\in S_N$. As an algebra,
$SL_q(N)$ is generated by $(x_{ij})_{1\leq i, j\leq N}$ with
relations (\ref{mq1})-(\ref{mq4}) and
\begin{equation}\label{mq5}
\sum\limits _{p\in S_N}(-q)^{-l(p)}x_{1p(1)}\cdots x_{Np(N)}=1.
\end{equation}
The coalgebra structure is the comatrix coalgebra structure.  To
define the antipode $S$ of ${\rm SL}_q(N)$, denote
$X=(x_{ij})_{1\leq i, j\leq N}$ and then define $Y_{ij}$ as the
generic $(N - 1)\times (N - 1)$ matrix obtained by deleting the
$i$th row and the $j$th column of $X$. Then
\begin{eqnarray}
\hspace*{1cm}
S(x_{ij})\hspace*{-2mm}
&=&\hspace*{-2mm}
(-q)^{j - i}{\rm det}_q(Y_{ji})\nonumber   \\
\hspace*{-2mm}
&=&
\hspace*{-2mm}
(-q)^{j - i} \sum
\limits_{p\in S_{j,i}}(-q)^{-l(p)}  x_{1p(1)}\cdots
x_{j-1p(j-1)}x_{j+1p(j+1)} \cdots x_{Np(N)} ,\label{mq6}
\end{eqnarray}
where $S_{j,i}$ is the set of bijective maps from $\{1, \ldots,
j-1,j+1, \ldots, N \}$ to $\{1, \ldots,i-1,i+1, \ldots, N \}$.
Moreover, for all $1\leq i, j\leq N$, we have
\begin{equation}\label{mq7}
S^2(x_{ij})=q^{2(j -i)}x_{ij}.
\end{equation}

We include the next lemma to provide complete details of the
construction.

\begin{lemma}
For $\sigma$ the skew pairing defined by (\ref{zq}) to (\ref{otherwise}),
for all $x \in M_q(N)$,
\begin{equation}\label{sigmawelldefined}
\sigma ({\rm det}_q, x)=\sigma  (x, {\rm det}_q)=\va(x)
\end{equation}
and thus $\sigma$ is well defined on $SL_q(N)\ot SL_q(N)$.
\end{lemma}
\begin{proof}
Since ${\rm det}_q$ is a grouplike element, from (\ref{cqt1}) and
(\ref{cqt2}), it suffices to check that  (\ref{sigmawelldefined})
holds on generators. From (\ref{smqn}) it follows that
$\sigma'(x_{im}, x_{jn})=0$ unless $i\leq m $, and if $i=m$, then
$\sigma'(x_{im}, x_{jn})=q^{\delta_{i,j}}\delta_{n,j}$. Now we
compute
\begin{eqnarray*}
&&\hspace*{-1.3cm}
\sigma'({\rm det}_q, x_{ij})\\
&=&\sum \limits_{p\in S_N}(-q)^{-l(p)}\sigma' (x_{1p(1)}\cdots x_{Np(N)}, x_{ij})\\
&=& \sum \limits_{p\in S_N}(-q)^{-l(p)}\sum \limits_{k=1}^N
\sigma' (x_{1p(1)}\cdots x_{N-1p(N-1)}, x_{ik})\sigma'(x_{Np(N)},x_{kj}  )\\
&=& \sum \limits_{p\in S_N,~p(N)=N}(-q)^{-l(p)}\sum \limits_{k=1}^N
\sigma' (x_{1p(1)}\cdots x_{N-1p(N-1)}, x_{ik})\sigma'(x_{NN},x_{kj})\\
&=&q^{\d_{N, j}}\sum \limits_{p\in S_{N-1} }(-q)^{-l(p)}
\sigma'(x_{1p(1)}\cdots x_{N-1p(N-1)}, x_{ij}).
\end{eqnarray*}
By induction, it follows that
\[
\sigma'({\rm det}_q, x_{ij})= q^{\sum \limits_{k=1}^{N }\d_{ k, j}}
\d_{i, j}=q\va(x_{ij}),
\]
and hence, $\sigma ({\rm det}_q, x_{ij})=z^N q \va(x_{ij})=
\va(x_{ij})$, for all $1\leq i, j\leq N$.

In a similar manner, using the fact that $\sigma'(x_{im}, x_{jn})=0$
unless $j \geq n $, and if $j=n$, then
$\sigma'(x_{im}, x_{jn})=q^{\delta_{i,j}}\delta_{m,i}$,
one can prove that
\[
\sigma (x_{ij}, {\rm det}_q)=
z^Nq^{\sum\limits_{k=1}^N\d_{k, j}}\d_{i, j}=\va(x_{ij}).
\]
Thus $\sigma$ is well defined on the quotient $M_q(N)/(det_q -1)$
and $SL_q(N)$ has a CQT structure.
\end{proof}

We now apply the results of Section \ref{cqtbialgs} to the CQT Hopf
algebra $H = (SL_q(N), \sigma)$. For all $1\leq i, j\leq N$, let us
denote $\br_{ij}:=r_{x_{ij}}$ and $\bl_{ij}:=l_{x_{ij}}$. (Note that
we will insert commas in these subscripts only for more complicated
expressions.)

\begin{lemma}\lelabel{4.1}
Let $H= (SL_q(N), \sigma)$.
\begin{eqnarray}
&&\br_{ij}(x_{mn})=z(q - q^{-1})\d_{i, n}\d_{j, m},~~\forall~~i<j;\label{br}\\
&&\bl_{ij}(x_{mn})=z(q - q^{-1})\d_{i, n}\d_{j, m},~~\forall~~i>j;\label{bl}\\
&&\bl_{ii}(x_{mn})=\br_{ii}(x_{mn})=zq^{\delta_{i,m}}\d_{m, n};\label{blr}\\
&&\bl_{ij}=\br_{ji}=0,~~\forall~~i<j. \label{bl0}
\end{eqnarray}
\end{lemma} \begin{proof} Equations (\ref{br})-(\ref{blr}) follow directly
from (\ref{zq}) to (\ref{zbracket}). Now suppose that $i < j$.  From
(\ref{otherwise}), we have that $\bl_{ij}$ and $\br_{ji}$ are 0 on
generators $x_{mn}$. Let   $a, b\in \{x_{mn}\mid 1\leq m, n\leq
N\}$. Since $\sigma$ is a skew pairing, we have
\[
\bl_{ij}(ab)= \sigma(ab, x_{ij}) = \sum
\limits_{k=1}^N\bl_{ik}(a)\bl_{kj}(b).
\]
Since $i<j$ we cannot have both $i\geq k$ and $k\geq j$, so it
follows that the product $\bl_{ik}(a)\bl_{kj}(b)=0$, for all
$k=\ov{1, N}$. By mathematical induction it follows that
$\bl_{ij}=0$, for any $i<j$. The   statement for $\br_{ji}$ is
proved similarly.
\end{proof}
\begin{corollary}\label{grplikes} The
maps $\bl_{ii}$ and $\br_{ii}$ are  equal   grouplike elements of
$H^0$. Also, denoting $\bl_{ii}^{-1} = \bl_{S^{-1}(x_{ii})}=
S^*(\bl_{ii})$, we have \begin{equation}
\bl_{ii}^{-1}(x_{mn})=\br_{ii}^{-1}(x_{mn})=z^{-1}q^{-\delta_{i,m}}\d_{m,
n}.\label{blrinv} \end{equation}
\end{corollary}
\begin{proof} The fact that $\bl_{ii}$ and $\br_{ii}$ are grouplike
follows directly from (\ref{bl0}).  Then, since these are algebra maps equal on
generators, they are equal. The rest is immediate.
\end{proof}

The next lemma describes the commutation relations for the
generators of  $H_r$ and $H_l$.

\begin{lemma}\label{commutation}
In $H_r$ the commutation relations for the generators $\br_{ij}$
are given by:
\begin{eqnarray}
&&\br_{im}\br_{in}=q\br_{in}\br_{im}, ~~\forall~~ n <m;~~
\br_{jm}\br_{im}=q\br_{im}\br_{jm},~~\forall~~i<j;\label{hr1}\\
&& \br_{jn}\br_{im}=\br_{im}\br_{jn};~~ \br_{jm}\br_{in} -
\br_{in}\br_{jm} =(q - q^{-1})\br_{im}\br_{jn},~~
\forall~~i<j,~~n<m.\label{hr2}
\end{eqnarray}
In $H_l$ the commutation relations for the generators are given by:
\begin{eqnarray}
&&\bl_{in}\bl_{im}=q\bl_{im}\bl_{in},~~\forall~~n<m~~; ~~
\bl_{im}\bl_{jm}=q\bl_{jm}\bl_{im},~~\forall~~i<j;\label{hl1}\\
&& \bl_{im}\bl_{jn} = \bl_{jn}\bl_{im}; ~~ \bl_{in}\bl_{jm} -
\bl_{jm}\bl_{in} =(q - q^{-1})\bl_{jn}\bl_{im},~~ \forall~~i<j~~{\rm
and}~~n<m.\label{hl2}
\end{eqnarray}
As well, for all $i,j$, we have that
\begin{equation}\label{hl3'}
{\bf l}_{11 }{\bf l}_{22}\cdots {\bf l}_{NN}=
{\bf r}_{11 }{\bf r}_{22}\cdots {\bf r}_{NN}=\va
~~{\rm and}~~
\br_{ii}\br_{jj} = \bl_{ii}\bl_{jj}=\bl_{jj}\bl_{ii} =
\br_{jj}\br_{ii}.
\end{equation}
\end{lemma}
\begin{proof}
Relations (\ref{hr1}) to (\ref{hl2}) follow from (\ref{mq1}) to
(\ref{mq4}) together with (\ref{Hlstructure}), (\ref{Hrstructure}).
The first relation in (\ref{hl3'}) follows from (\ref{mq5}),
(\ref{bl0}) and Corollary \ref{grplikes}. By the second relation  in
(\ref{hl2}) we have that $\bl_{ii}\bl_{jj} - \bl_{jj}\bl_{ii} =(q -
q^{-1})\bl_{ji}\bl_{ij} $ if $i<j$, and $\bl_{jj}\bl_{ii} -
\bl_{ii}\bl_{jj}=(q - q^{-1})\bl_{ij}\bl_{ji} $  if $j<i$. By
(\ref{bl0}), these are both 0. Therefore,
$\bl_{ii}\bl_{jj}=\bl_{jj}\bl_{ii}$, for all $1\leq i, j\leq N$.
\end{proof}

Note that (\ref{hr2}) and (\ref{hl2}) imply that
\begin{eqnarray}\label{mge} \nonumber && \br_{jm}\br_{in}= \br_{in}\br_{jm}
~~ \mbox{if   $i<j$, $n<m$ and either  $i>m$ or $j>n$};\\
&& \bl_{in}\bl_{jm}=\bl_{jm}\bl_{in}~~ \mbox{if   $i<j$, $n<m$ and
either  $i<m$ or $j<n$}.
\end{eqnarray}

We now describe another set of algebra generators for $H_l$ and
$H_r$.

\begin{proposition}\prlabel{4.3}
Let $1\leq i\leq N$ and $1\leq s, t\leq N - 1$ and define:
\begin{eqnarray*}
&&\hk_i:=\bl_{ii}, \hspace{.2cm}
\hk_i^{-1}:=\bl_{ii}^{-1};\\
&& E_s:=\bl_{s+1,s+1}^{-1}\bl_{s+1,s} = q \bl_{s+1,s}\bl_{s+1,s+1}^{-1};\\
&& F_s:=(q - q^{-1})^{-2}\br_{ss}^{-1}\br_{s,s+1}= (q -
q^{-1})^{-2}q \br_{s,s+1}\br_{ss}^{-1}.
\end{eqnarray*}
As an algebra $H_l$ is generated by $\hk_i$, $\hk_i^{-1}$, and the $E_s$,
with the following relations.
\begin{eqnarray}
&&\hk_i\hk_j=\hk_j\hk_i,~~
\hk_i\hk_i^{-1}=\hk_i^{-1}\hk_i=\va,~~\hk_1\hk_2\cdots \hk_N=\va ;\label{hl1s}\\
&&\hk_iE_t\hk_i^{-1}=q^{\d_{i, t} - \d_{i, t+1}}E_t;\label{hl2s}\\
&&E_tE_s=E_sE_t,~~\mbox{if $\mid s - t\mid >1$};\label{hls3}\\
&&E_s^2E_t - (q + q^{-1})E_sE_tE_s + E_tE_s^2=0, ~~\mbox{if $\mid s
- t\mid =1$}.\label{hls4}
\end{eqnarray}

Similarly, $H_r$ is generated as an algebra by the $\hk_i$, $\hk_i^{-1}$,
and the $F_s$, with relations (\ref{hl1s}) and
\begin{eqnarray}
&&\hk_iF_t\hk_i^{-1}=q^{\d_{i, t + 1} - \d_{i, t}}F_t;\label{hr1s}\\
&&F_tF_s=F_sF_t,~~\mbox{if $\mid s - t\mid >1$};\label{hr2s}\\
&&F_s^2F_t - (q + q^{-1})F_sF_tF_s + F_tF_s^2=0,~~ \mbox{if $\mid s
- t\mid =1$}.\label{hr3s}
\end{eqnarray}
\end{proposition}
\begin{proof}
We prove only the statements for   $H_l$; the proofs for  $H_r$
are similar. Note that the map $\hk_i\mapsto \hk_i^{-1}$,
$E_s\mapsto F_s$, is a well defined   algebra isomorphism from $H_l$
to $H_r$.

Since $\{\bl_{ij}\mid i\geq j\}$ is a set of algebra generators for
$H_l$, it suffices to prove that for any $1\leq i\leq N$ and $1\leq
j\leq N-i$ the element $\bl_{i+j,i}$ can be written as a linear
combination of products of the elements $\hk_1^{\pm 1}, \cdots ,
\hk_N^{\pm 1}, E_1, \cdots , E_{N-1}$ to see that these elements
generate.  For $1 \leq i \leq N-1$, we have $\bl_{i+1,i}=
\hk_{i+1}E_i = q^{-1}E_i\hk_{i+1} $. Suppose that $\bl_{i+j,i}$ can
be written as a linear combination of products of the elements
$\hk_1^{\pm 1}, \cdots , \hk_N^{\pm 1}, E_1, \cdots , E_{N-1}$. By
the second relation in (\ref{hl2}),
\begin{equation*} \bl_{i+j,i}\bl_{i+j+1,i+j} -
\bl_{i+j+1,i+j}\bl_{i+j,i} =(q - q^{-1})\bl_{i+j+1,i}\bl_{i+j,i+j}
\end{equation*}
for all $1\leq j\leq N-i-1$. Now we compute
\begin{eqnarray*}
&& \hspace*{-1.3cm}
\bl_{i+j+1,i}\\
&=&\hspace*{-2mm}
(q - q^{-1})^{-1}[\bl_{i+j,i}\bl_{i+j+1,i+j}\bl_{i+j,i+j}^{-1} -
\bl_{i+j+1,i+j}\bl_{i+j,i}\bl_{i+j,i+j}^{-1}]\\
&{{(\ref{hl1}, \ref{hl2})}\atop =}&\hspace*{-2mm}
(q - q^{-1})^{-1}q^{-1}
[\bl_{i+j,i}E_{i+j}\bl_{i+j+1,i+j+1}\bl_{i+j,i+j}^{-1} -
\bl_{i+j+1,i+j}
\bl_{i+j,i+j}^{-1}\bl_{i+j,i}]\\
&{{(\ref{hl1})}\atop =}&\hspace*{-2mm}
(q^2 - 1)^{-1}
[\bl_{i+j,i}E_{i+j}\hk_{i+j+1}\hk_{i+j}^{-1} -
q^{-1}E_{i+j}\hk_{i+j+1} \hk_{i+j}^{-1}\bl_{i+j,i}],
\end{eqnarray*}
and the statement follows by mathematical induction.

Now we show that (\ref{hl1s}) to (\ref{hls4}) hold. Relations
(\ref{hl1s}) are immediate from (\ref{hl3'}). In order to prove
(\ref{hl2s}) observe that
\[
\hk_iE_t\hk_i^{-1}=\bl_{ii}\bl_{t+1,t+1}^{-1}\bl_{t+1,t}\bl_{ii}^{-1}
=\bl_{t+1,t+1}^{-1}\bl_{ii}\bl_{t+1,t}\bl_{ii}^{-1}.
\]
We now verify (\ref{hl2s}) case by case.  \begin{itemize}
\item[ Case I:]
If  $i < t $ or $ i>t+1$, by (\ref{hl2}) we have
$\bl_{ii}\bl_{t+1,t}=\bl_{t+1,t}\bl_{ii}$, and therefore
\[
\hk_iE_t\hk_i^{-1}=\bl_{t+1,t+1}^{-1}\bl_{t+1,t}=E_t=q^{\d _{i, t} -
\d _{i, t + 1}}E_t.
\]
\item[ Case II:]  If $i=t$, by (\ref{hl1}) we have
$\bl_{tt}\bl_{t+1,t}=q\bl_{t+1,t}\bl_{tt}$, so that
\[
\hk_tE_t\hk_t^{-1}=q\bl_{t+1,t+1}^{-1}\bl_{t+1,t}=qE_t=q^{\d_{t, t}
- \d_{t, t+1}}E_t.
\]
\item[ Case III:] If $i=t+1 $,   again using (\ref{hl1}) we have
$\bl_{t+1,t}\bl_{t+1,t+1}=q\bl_{t+1,t+1}\bl_{t+1,t}$, and
\[
\hk_{t+1}E_t\hk_{t+1}^{-1}=\bl_{t+1,t}\bl_{t+1,t+1}^{-1}=q^{-1}E_t=
q^{\d_{t+1, t} - \d_{t+1, t+1}}E_t.
\]
\end{itemize}
Now to verify (\ref{hls3}), we first note that if $t<s - 1$,
(\ref{mge}) implies that the pairs   $\bl_{t+1,t}$ and
$\bl_{s+1,s+1}$, $\bl_{t+1,t}$ and $\bl_{s+1,s}$, and $\bl_{s+1,s}$
and $\bl_{t+1,t+1}$ commute and thus
\[
E_tE_s=\bl_{t+1,t+1}^{-1}\bl_{t+1,t}\bl_{s+1,s+1}^{-1}\bl_{s+1,s}=E_sE_t.
\]
If $t>s+1$ then $s<t-1$, so interchanging $s$ and $t$ in the
argument above we obtain $E_sE_t=E_tE_s$.

It remains to prove the relation  in (\ref{hls4}). We first compute
\begin{eqnarray*}
E_sE_{s+1}&=&\bl_{s+1,s+1}^{-1}\bl_{s+1,s}\bl_{s+2,s+2}^{-1}\bl_{s+2,s+1}\\
&{{(\ref{mge})}\atop =}&\bl_{s+1,s+1}^{-1}\bl_{s+2,s+2}^{-1}\bl_{s+1,s}\bl_{s+2,s+1}\\
&{{(\ref{hl2})}\atop =}&\bl_{s+2,s+2}^{-1}\bl_{s+1,s+1}^{-1}
[\bl_{s+2,s+1}\bl_{s+1,s} + (q - q^{-1})\bl_{s+2,s}\bl_{s+1,s+1}]\\
&{{(\ref{hl1}, \ref{hl2})}\atop =}&
q^{-1}\bl_{s+2,s+2}^{-1}\bl_{s+2,s+1}\bl_{s+1,s+1}^{-1}\bl_{s+1,s} +
(q - q^{-1})\bl_{s+2,s+2}^{-1}\bl_{s+2,s}.
\end{eqnarray*}
In other words, we have proved that
\begin{equation}\label{n1}
E_sE_{s+1}=q^{-1}E_{s+1}E_s + (q -
q^{-1})\bl_{s+2,s+2}^{-1}\bl_{s+2,s},
\end{equation}
for all $1\leq s\leq N-2$. Clearly, this is equivalent to
\begin{equation}\label{en1}
E_{s+1}E_s=qE_sE_{s+1} - q(q - q^{-1})\bl_{s+2,s+2}^{-1}\bl_{s+2,s},
\end{equation}
for all $1\leq s\leq N-2$. Using the two relations above we obtain
\begin{eqnarray*}
E_s^2E_{s+1}&=&q^{-1}E_sE_{s+1}E_s + (q - q^{-1})E_s
\bl_{s+2,s+2}^{-1}\bl_{s+2,s},\\
E_{s+1}E_s^2&=&qE_sE_{s+1}E_s - q(q -
q^{-1})\bl_{s+2,s+2}^{-1}\bl_{s+2,s}E_s,
\end{eqnarray*}
from which we compute
\begin{eqnarray*}
&&\hspace*{-1cm}
E_s^2E_{s+1} - (q^{-1} + q)E_sE_{s+1}E_s + E_{s+1}E_s^2\\
&&\hspace*{1cm} =(q - q^{-1})[E_s\bl_{s+2,s+2}^{-1}\bl_{s+2,s} -
q\bl_{s+2,s+2}^{-1}\bl_{s+2,s}E_s].
\end{eqnarray*}
Now, since
\begin{eqnarray*}
E_s\bl_{s+2,s+2}^{-1}\bl_{s+2,s}&=&
\bl_{s+1,s+1}^{-1}\bl_{s+1,s}\bl_{s+2,s+2}^{-1}\bl_{s+2,s}\\
&{{(\ref{mge})}\atop =}&\bl_{s+1,s+1}^{-1}\bl_{s+2,s+2}^{-1}\bl_{s+1,s}\bl_{s+2,s}\\
&{{(\ref{hl1})}\atop =}&q\bl_{s+2,s+2}^{-1}\bl_{s+1,s+1}^{-1}\bl_{s+2,s}\bl_{s+1,s}\\
&{{(\ref{hl2})}\atop
=}&q\bl_{s+2,s+2}^{-1}\bl_{s+2,s}\bl_{s+1,s+1}^{-1}\bl_{s+1,s}
=q\bl_{s+2,s+2}^{-1}\bl_{s+2,s}E_s,
\end{eqnarray*}
it follows that $E_s^2E_{s+1} - (q^{-1} + q)E_sE_{s+1}E_s +
E_{s+1}E_s^2=0$, for all $1\leq s\leq N-2$.

Similarly, again using (\ref{n1}) and (\ref{en1}) we have
\begin{eqnarray*}
E_{s+1}^2E_s&=&qE_{s+1}E_sE_{s+1} - q(q - q^{-1})E_{s+1}\bl_{s+2,s+2}^{-1}\bl_{s+2,s},\\
E_sE_{s+1}^2&=&q^{-1}E_{s+1}E_sE_{s+1} + (q -
q^{-1})\bl_{s+2,s+2}^{-1}\bl_{s+2,s}E_{s+1},
\end{eqnarray*}
and therefore
\begin{eqnarray*}
&&\hspace*{-1cm}
E_{s+1}^2E_s - (q + q^{-1})E_{s+1}E_sE_{s+1} + E_sE_{s+1}^2\\
&&\hspace*{1cm} =(q - q^{-1})[\bl_{s+2,s+2}^{-1}\bl_{s+2,s}E_{s+1} -
qE_{s+1}\bl_{s+2,s+2}^{-1}\bl_{s+2,s}].
\end{eqnarray*}
Now, we have
\begin{eqnarray*}
qE_{s+1}\bl_{s+2,s+2}^{-1}\bl_{s+2,s}
&{{(\ref{hl1})}\atop =}&q^2\bl_{s+2,s+2}^{-1}\bl_{s+2,s+1}\bl_{s+2,s}\bl_{s+2,s+2}^{-1}\\
&{{(\ref{hl1})}\atop =}&q\bl_{s+2,s+2}^{-1}\bl_{s+2,s}\bl_{s+2,s+1}\bl_{s+2,s+2}^{-1}\\
&{{(\ref{hl1})}\atop
=}&\bl_{s+2,s+2}^{-1}\bl_{s+2,s}\bl_{s+2,s+2}^{-1}\bl_{s+2,s+1}
=\bl_{s+2,s+2}^{-1}\bl_{s+2,s}E_{s+1},
\end{eqnarray*}
and therefore $E_{s+1}^2E_s - (q + q^{-1})E_{s+1}E_sE_{s+1} +
E_sE_{s+1}^2=0$, for all $1\leq s\leq N-2$, and this finishes the
proof.
\end{proof}

Next we fix some notation:
\begin{definition}\label{braidedcomm} For any $x$ and $y$, let $[x, y]_q:=qxy - yx$.
\begin{itemize}
\item[(i)] For $1 \leq i \leq j \leq N$, define  $ {\cal F}_{i, i}:= 1$,
${\cal F}_{i, i + 1}:= q(1-q^{-2})^2 F_i$, and
\[
{\cal F}_{i, j}:=q(1-q^{-2})^{j-i+1}[F_{j-1}, [F_{j-2},\cdots
[F_{i+1}, F_i]_q\cdots ]_q]_q,  ~~\mbox{for}~~j \geq i+2 ,
\]
i.e., the ${\cal F}_{i, j}$  are defined inductively by
${\cal F}_{i, j+1} = (1-q^{-2})[F_j, {\cal F}_{i, j}]_q$.
\item[(ii)] Similarly, define $ {\cal E}_{j,j}:= 1$,
$ {\cal E}_{j + 1, j }:= q^{-1}E_j$ and,
\[
{\cal E}_{i, j}:=q^{-1}(q^2 -1)^{j-i}[\cdots [[E_j, E_{j+1}]_q,
E_{j+2}]_q, \cdots , E_{i-1}]_q, ~~\mbox{for}~~ i \geq j+2,
\] i.e., the ${\cal E}_{i,j}$ are defined inductively by ${\cal
E}_{i+1,j} = (q^2 -1)^{-1}[{\cal E}_{i,j}, E_i]_q$.
\end{itemize}
\end{definition}
\vspace{1mm}
Now, using \prref{4.3},
$\bl_{ij}$ and $\br_{ij}$ can be expressed in terms of the new
generators $\hk_i^{\pm 1}$, $E_s$  and $F_s$ in the spirit of Serre
relations.
\begin{corollary}\colabel{4.3}
For $1 \leq i \leq j \leq N$, we have that $\bl_{ji}={\cal E}_{j,i}\hk_j$,
and $\br_{ij}={\cal F}_{i,j}\hk_i$.
\end{corollary}
\begin{proof}
We prove  the statement for $\br_{ij}$ by
mathematical induction on $j$; the similar proof for $\bl_{ji}$ is
left to the reader. If $j=i$, the statement is clear. If $j=i+1$,
then the statement follows from Lemma \ref{commutation} and the
definitions. Now assume that  $\br_{ij} = {\cal F}_{ij}\hk_i$ for
$\br_{ij}$ with $j \geq i+1$ and we show that it holds for
$\br_{i,j+1}$. By (\ref{hr2}), we have
\[
\br_{j,j+1}\br_{i,j} -
\br_{i,j}\br_{j,j+1} = (q-q^{-1})\br_{i,j+1}\br_{j,j},
\]
and thus
\begin{eqnarray*}
\br_{i,j+1}&=&(q-q^{-1})^{-1}[\br_{j,j+1} \br_{i,j}\br_{j,j}^{-1} -
\br_{i,j}\br_{j,j+1}\br_{j,j}^{-1}]\\
&{{(\ref{hr1})}\atop
=}&(q-q^{-1})^{-1}[q\br_{j,j+1}\br_{j,j}^{-1}\br_{i,j} -
q^{-1}\br_{i,j}\br_{j,j}^{-1}\br_{j,j+1}]\\
&{{(\ref{hr1})}\atop =}&(q-q^{-1})[F_{j}\br_{i,j} -
q^{-1}\br_{i,j}F_{j}]\\
& = &(1-q^{-2})[F_{j}, \br_{i,j}]_q\\
&=& (1 - q^{-1})[F_j, {\cal F}_{ij} \hk_i]_q \mbox{by the
induction assumption}\\
&{{(\ref{hr1s})}\atop=}&
(1 -q^{-2})[F_j, \cal{F}_{ij}]_q \hk_i\\
&=&(1 - q^{-2}) {\cal F}_{i,j+1} \hk_i,
\end{eqnarray*}
and the proof is complete.
\end{proof}

The next proposition describes the comultiplication, counit and antipode for
$H_l$ and $H_r$.

\begin{proposition}\prlabel{4.4}
The coalgebra structure $\Delta$, $\va$ and the antipode $S$ for
$H_l$ are determined by
\begin{eqnarray}
&&\Delta(\hk_i^{\pm 1})=\hk_i^{\pm 1}\ot\hk_i^{\pm 1},~~
\va(\hk_i^{\pm 1})=1,\label{hls5}\\
&&\Delta(E_s)=\va\ot E_s + E_s\ot \hk_{s+1}^{-1}\hk_s,~~
\va(E_s)=0,\label{hls6}\\
&&S(\hk_i^{\pm 1})=\hk_i^{\mp
1},~~S(E_s)=-E_s\hk_s^{-1}\hk_{s+1},\label{hls7}
\end{eqnarray}
for $1\leq i\leq N$ and $1\leq s\leq N-1$. Similarly,  the
coalgebra structure $\Delta$, $\va$ and the antipode $S$ for $H_r$
are determined by (\ref{hls5}), the first equality in (\ref{hls7})
and
\begin{eqnarray}\label{hr4s}
&&\Delta (F_s)=F_s\ot \va + \hk_{s+1}\hk_s^{-1}\ot F_s,~~\va(F_s)=0,\\
\label{hr5s} && S(F_s)=-\hk_s\hk_{s+1}^{-1}F_s.
\end{eqnarray}
\end{proposition}
\begin{proof}
We give the details for $E_s$ and
leave those for $F_s$ to the reader.
Also, we note that the algebra isomorphism $\hk_i\mapsto \hk_i^{-1}$,
$E_s\mapsto F_s$, defined in the proof of \prref{4.3},
is actually a Hopf algebra isomorphism between $H_l$ and $H_r^{\rm cop}$.

Now, we compute:
\begin{eqnarray*}
\Delta(E_s)&=&\Delta(\bl_{s+1,s+1}^{-1}\bl_{s+1,s})=(\bl_{s+1,s+1}^{-1}\ot
\bl_{s+1,s+1}^{-1})
\sum \limits _{k=s}^{s+1}\bl_{s+1,k}\ot \bl_{ks}\\
&=&(\bl_{s+1,s+1}^{-1}\ot \bl_{s+1,s+1}^{-1})(\bl_{s+1,s}\ot
\bl_{ss} + \bl_{s+1,s+1}\ot \bl_{s+1,s})\\
&=&E_s\ot \hk_{s+1}^{-1}\hk_s + \va \ot E_s,
\end{eqnarray*}
as needed. Also,
$\va(E_s)=\va(\bl_{s+1,s+1}^{-1}\bl_{s+1,s})=\d_{s+1, s}=0$, for all
$1\leq s\leq N-1$. The formula for $S(E_s)$ is then immediate.
\end{proof}

\begin{remark}\relabel{4.5}
For   $1\leq i\leq N-1$ we now define $K_i:=\hk_{i+1}^{-1}\hk_i
=\hk_i \hk_{i+1}^{-1} $, and then define $H'_l$ to be the subalgebra
of $H_l$ generated by $\{K_i^{\pm 1}, E_i\mid 1\leq i\leq N - 1\}$.
From \prref{4.3} one can easily see that the relations between the
algebra generators of $H'_l$ are (\ref{hls3}), (\ref{hls4}) and
\begin{equation}\label{hlp1}
K_i^{\pm 1}K_i^{\mp 1}=\va ,~~K_iK_j=K_jK_i,~~
K_iE_jK_i^{-1}=q^{a_{ij}}E_j,
\end{equation}
for all $1\leq i, j\leq N-1$, where $a_{ii}=2$, $a_{ij}=-1$ if
$\mid i - j\mid =1$, and $a_{ij}=0$ if $\mid i - j\mid >1$. Actually,
$H'_l$ is a Hopf subalgebra of $H_l$. The induced Hopf algebra
structure is given by
\begin{eqnarray*}
&&\Delta(K_i^{\pm 1})=K_i^{\pm 1}\ot K_i^{\pm 1},~~\va(K_i^{\pm 1})=1,\\
&&\Delta(E_i)=\va \ot E_i + E_i\ot K_i,~~\va(E_i)=0,\\
&&S(K_i^{\pm 1})=K_i^{\mp 1},~~S(E_i)=-E_iK_i^{-1},
\end{eqnarray*}
for all $1\leq i\leq N-1$. Similarly, define $H'_r$ to be the
subalgebra of $H_r$ generated by $\{K_i^{\pm}, F_i\mid 1\leq i\leq N
- 1\}$. Then its algebra generators satisfy the relations in
(\ref{hr2s}), (\ref{hr3s}), the first two relations in (\ref{hlp1})
and
\begin{equation}\label{hrp1}
K_iF_jK_i^{-1}=q^{-a_{ij}}F_j,~~\forall~~1\leq i, j\leq N-1,
\end{equation}
where $a_{ij}$ are the scalars defined above. Moreover, $H'_r$ is a
Hopf subalgebra of $H_r$. The elements $K_i^{\pm}$ are grouplike
elements and
\begin{equation}
\Delta(F_i)=K_i^{-1}\ot F_i + F_i\ot
\va,~~\va(F_i)=0,~~S(F_i)=-K_iF_i,
\end{equation}
for all $1\leq i\leq N-1$. Also, the map $K_i\mapsto K_i^{-1}$,
$E_i\mapsto F_i$, is well defined and a Hopf algebra
isomorphism from $H'_l$ to $H'^{\rm cop}_r$.
\end{remark}

We can now describe the Hopf algebra structure of
$H_{\sigma}:=H_lH_r$.

\begin{proposition}\prlabel{4.6}
The Hopf algebra $H_{\sigma}:=H_lH_r$ can be presented as follows.
It is the algebra generated by $\{\hk_i^{\pm 1}\mid 1\leq i \leq
N\}\cup \{E_s, F_s\mid 1\leq s\leq N-1\}$ with relations
(\ref{hl1s})-(\ref{hr3s}) and
\begin{equation}\label{hsigma}
E_sF_t - F_tE_s=  \d_{s, t}(q - q^{-1})^{-1}[\hk_s\hk_{s+1}^{-1} -
\hk_s^{-1}\hk_{s+1} ],
\end{equation}
for all $1\leq s, t\leq N-1$. The Hopf algebra structure of
$H_{\sigma}$ is given by the relations in \prref{4.4}. In other
words, $H_{\sigma}\equiv U_q^{\rm ext}({\rm sl}_N)$, the extended
Hopf algebra of $U_q({\rm sl}_N)$ defined in \cite[Section 8.5.3]{ks}.
\end{proposition}
\begin{proof}
It remains only to prove  (\ref{hsigma}). It suffices to show first
that (\ref{hsigma}) holds when applied to $x_{mn}$ for all $1\leq m,
n\leq N$,  and then to show that if the maps in (\ref{hsigma}) are
equal on $a$ and $b$ then they are also equal on $ab$.

The relations from \leref{4.1} will be key in the following
computations. Also we will sometimes use the notation $K_s =
\hk_s\hk_{s+1}^{-1}$ from \reref{4.5}.

Recall that  $\hk_i(x_{mn})=\bl_{ii}(x_{mn})=
zq^{\d_{i, m}}\d_{m, n}$ by (\ref{blr}) and $\hk_i^{-1}(x_{mn})=
\bl_{ii}^{-1}(x_{mn}) = z^{-1} q^{ -\d_{i, m}}\d_{m, n}$ by (\ref{blrinv}).
Thus,
\begin{eqnarray*}
&&\hspace*{-1cm}E_s(x_{mn})=\bl_{s+1,s+1}^{-1}\bl_{s+1,s}(x_{mn})
=\sum\limits_{k=1}^N\bl_{s+1,s+1}^{-1}(x_{mk})
\bl_{s+1,s}(x_{kn})\\
&&\hspace*{1cm} =q^{-\d_{m, s+1}}(q - q^{-1})\d_{s+1, n}\d_{s, m}=(q
- q^{-1})\d_{s+1, n}\d_{s, m}.
\end{eqnarray*}
Similarly, we compute  $F_t(x_{mn})=(q - q^{-1})^{-1} \d_{t,
n}\d_{t+1, m}$, and therefore
\begin{eqnarray*}
E_sF_t(x_{mn})&=&\sum\limits_{k=1}^NE_s(x_{mk})F_t(x_{kn})
=(q - q^{-1})\d_{s, m}F_t(x_{s+1,n})\\
&=&\d_{s, m}\d_{t, n}\d_{s, t}=\d_{s, m}\d_{s, t}\d_{m, n},
\end{eqnarray*}
and
\begin{eqnarray*}
F_tE_s(x_{mn})&=&\sum\limits_{k=1}^nF_t(x_{mk})E_s(x_{kn})
=(q - q^{-1})^{-1}\d_{t+1, m}E_s(x_{tn})\\
&=&\d_{t+1, m}\d_{s+1, n}\d_{s, t}=\d_{s+1, m}\d_{s, t}\d_{m, n}.
\end{eqnarray*}
We conclude that
\begin{equation}\label{hsigma2}
(E_sF_t - F_tE_s)(x_{mn})=\d_{s, t}[\d_{s, m} - \d_{s+1, m}]\d_{m, n}.
\end{equation}
On the other hand, we have
\begin{eqnarray*}
\hk_s\hk_{s+1}^{-1}(x_{mn})&=&
\sum\limits_{k=1}^N\hk_s(x_{mk})\hk_{s+1}^{-1}(x_{kn})\\
&=& z  q^{  \d_{s, m}}\d_{m,k} z^{-1} q^{ -\d_{s+1, m}}\d_{k,n}
=q^{\d_{s, m} - \d_{s+1, m}}\d_{m, n}.
\end{eqnarray*}
Similarly, $\hk_s^{-1}\hk_{s+1} (x_{mn})=q^{\d_{s+1, m} - \d_{s,
m}}\d_{m,n}$. Thus
\begin{equation*}
[\hk_s\hk_{s+1}^{-1} - \hk_s^{-1}\hk_{s+1}](x_{mn})\\
= [q^{\d_{s, m} - \d_{s+1, m}} - q^{\d_{s+1, m} - \d_{s, m}}]\d_{m, n}.
\end{equation*}
Since
\[
q^{\d_{s, m} - \d_{s+1, m}} - q^{\d_{s+1, m} -
\d_{s, m}}=\left\{\begin{array}{ll}
0&,\mbox{if $m\not\in \{s, s+1\}$}\\
q - q^{-1}&,\mbox{if $m=s$}\\
-(q - q^{-1})&, \mbox{if $m=s+1$}
\end{array}\right. ,
\]
it follows that $\d_{s, t}(q - q^{-1})^{-1}[K_s - K_s^{-1}](x_{mn})=
\d_{s, t}[\d_{s, m} - \d_{s+1, m}]\d_{m, n}$. Together with
(\ref{hsigma2}) this proves  that (\ref{hsigma}) holds on
generators.

\par Since the $K_s$ are grouplike, the right hand side of
(\ref{hsigma}) is $(K_s, K_s^{-1})$-primitive. Also the left hand
side is $(K_s, K_t^{-1})$-primitive.  This follows from the fact
that the matrix $(a_{ij})_{1\leq i, j\leq N}$ defined in \reref{4.5}
is symmetric, and thus
\[
K_t^{-1}E_s\ot F_tK_s=q^{-a_{ts}}E_sK_t^{-1}\ot
q^{a_{st}}K_sF_t=E_sK_t^{-1}\ot K_sF_t.
\]

If $s \neq t$, then both sides of  (\ref{hsigma}) are 0 on
generators $a, b \in SL_q(N)$, and thus, being skew-primitive,  on
$ab$. Otherwise, both sides are $(K_s, K_s^{-1})$-primitive and thus
equal on the product $ab$. The statement then follows by induction.
\end{proof}

From \prref{4.4} and \prref{4.6} we conclude that the Hopf algebras
$H_l$ and $H_r$ are precisely $U_q(\mf{b}_-)$ and $U_q(\mf{b}_+)$,
the Borel-like Hopf algebras associated to $U^{\rm ext}_q({\rm
sl}_N)$. The fact that $U_q(\mf{b}_+)^{\rm cop}$ and $U_q(\mf{b}_-)$
are Hopf algebras in duality is well known but using the general
theory above we are now able to present a more conceptual proof.

\begin{corollary}\colabel{4.7}
The pairing $\le , \ri :U_q(\mf{b}_+)\ot U_q(\mf{b}_-)\ra k$
defined by
\begin{eqnarray*}
&&\le \hk_i, \hk_j \ri =\le \hk_i^{-1}, \hk_j^{-1}\ri = z q^{\d_{i,
j}  },~~ \le \hk_i, \hk_j^{-1}\ri
=\le \hk_i^{-1}, \hk_j\ri = z^{-1}q^{ - \d_{i, j}},\\
&&\le \hk_i, E_s\ri =\le \hk_i^{-1}, E_s\ri =0,~~
\le F_s, \hk_j\ri =\le F_s, \hk_j^{-1}\ri =0,\\
&&\le F_s, E_t\ri = (q - q^{-1})^{-1}\d_{s, t},
\end{eqnarray*}
for all $1\leq i, j\leq N$ and $1\leq s, t\leq N - 1$, defines a
duality between $U_q(\mf{b}_+)^{\rm cop}$ and $U_q(\mf{b}_-)$ as in
\prref{2.1}.

Also, $U_q(\mf{b}_-)$ is self dual, in the sense that there is a
duality between $U_q(\mf{b}_-)$ and itself  given for all $1\leq i,
j\leq N$ and $1\leq s, t\leq N-1$ by
\begin{eqnarray*}
&&\le \hk_i, \hk_j \ri =\le \hk_i^{-1}, \hk_j^{-1}\ri = z^{-1}q^{  -
\d_{i, j}},~~
\le \hk_i, \hk_j^{-1}\ri =\le \hk_i^{-1}, \hk_j\ri =zq^{\d_{i, j} },\\
&&\le \hk_i, E_s\ri =\le \hk_i^{-1}, E_s\ri =
\le E_s, \hk_i\ri =\le E_s, \hk_i^{-1}\ri =0,\\
&&\le E_s, E_t\ri = (q - q^{-1})^{-1}\d_{s, t}.
\end{eqnarray*}
\end{corollary}
\begin{proof}
The first set of equations follows  from \prref{2.1} and from
(\ref{zq})- (\ref{otherwise}).   We verify only the third equation
of the first set.  Note first that $\le \br_{ss},  \bl_{ij} \ri = 0
$ unless $i=j$ by (\ref{zq})- (\ref{otherwise}) so that $\le
\br_{ss}^{-1}, E_t \ri = 0$.  Now we compute
\begin{eqnarray*}
\le F_s, E_t\ri
&=&(q - q^{-1})^{-2}\le \br^{-1}_{ss}\br_{s,s+1}, E_t\ri \\
&{{(\ref{hls6})}\atop =}&(q - q^{-1})^{-2}[\le \br_{ss}^{-1}, \va
\ri \le \br_{s,s+1}, E_t\ri + \le \br_{ss}^{-1}, E_t\ri \le
\br_{s,s+1},
\hk_{t+1}^{-1}\hk_t\ri ]\\
&=&(q - q^{-1})^{-2}\le \br_{s,s+1}, \bl_{t+1,t+1}^{-1}\bl_{t+1,t}\ri \\
&=&(q - q^{-1})^{-2}[\le \br_{ss}, \bl_{t+1,t+1}^{-1}\ri
\le \br_{s,s+1}, \bl_{t+1,t}\ri \\
&&\hspace*{3cm}
+ \le \br_{s,s+1}, \bl_{t+1,t+1}^{-1}\ri \le \br_{s+1,s+1}, \bl_{t+1,t}\ri ]\\
&=&(q - q^{-1})^{-2}\hk_{t+1}^{-1}(x_{ss})\bl_{t+1,t}(x_{s,s+1})\\
&{{(\ref{bl}, \ref{blr})}\atop =}& (q - q^{-1})^{-2}z^{-1} q^{  -
\d_{s, t+1}}z  (q - q^{-1})\d_{s, t} =(q - q^{-1})^{-1}\d_{s, t},
\end{eqnarray*}
as we claimed. The second set of equations follows from the
isomorphism between $H_r^{\rm cop}$ and $H_l$ described at the
beginning of the proof of \prref{4.4}.
\end{proof}

\begin{corollary}\colabel{4.8}
The Hopf algebra $U_q^{\rm ext}({\rm sl}_N)$ is isomorphic to a
factor of the generalized quantum double $D(U_q(\mf{b}_+),
U_q(\mf{b}_-))\equiv D(U_q(\mf{b}_-)^{\rm cop}, U_q(\mf{b}_-))$.
\end{corollary}
\begin{proof}
It is an immediate consequence of \coref{4.7} and \prref{1.2}.
\end{proof}

\begin{remark}\relabel{4.9}
Further to \reref{4.5}, define $H'_{\sigma}=H'_lH'_r =
H'_{\sigma}=H'_rH'_l$,  a Hopf subalgebra of $H_{\sigma}\equiv
U_q^{\rm ext}({\rm sl}_N)$. As expected, $H'_{\sigma}\equiv U_q({\rm
sl}_N)$, the Hopf algebra with algebra generators $\{K_i^{\pm 1},
E_i, F_i\mid 1\leq i\leq N-1\}$ and relations from \reref{4.5} as
well as
\[
E_iF_j - F_jE_i=\d_{ij}\frac{K_i - K_i^{-1}}{q - q^{-1}},~~
\forall~~1\leq i, j\leq N -1,
\]
from  \prref{4.6}. The comultiplication, the counit and the antipode are
defined in \reref{4.5}. (For more details see \cite[VI.7\& VII.9]{k}.)

Note that $H'_r$ and $H'_l$ are exactly the Borel-like Hopf algebras
associated to $U_q({\rm sl}_N)$ and the situation is similar to that
of  \coref{4.7}, i.e., there is a duality on $H'^{\rm cop}_r \ot
H'_l$  given by
\begin{eqnarray*}
&&\le K_i, K_j\ri =\le K_i^{-1}, K_j^{-1}\ri =q^{a_{ij}},~~
\le K_i^{-1}, K_j\ri =\le K_i, K_j^{-1}\ri =q^{-a_{ij}},\\
&&\le K_i, E_s\ri =\le K_i^{-1}, E_s\ri=\le F_s, K_i\ri=\le F_s, K_i^{-1}\ri =0,\\
&&\le F_s, E_t\ri = (q - q^{-1})^{-1}\d_{s, t},
\end{eqnarray*}
for all $1\leq i, j\leq N$ and $1\leq s, t\leq N - 1$. Moreover,
since $H'^{\rm cop}_r$ is isomorphic to $H'_l$ , then $H'_l$ is self
dual and   the generalized quantum double $D(H'_r, H'_l)\equiv
D((H'_l) ^{\rm cop},   H'_l  )$. Furthermore, one can easily see
that the associated generalized quantum double $D(H'_r, H'_l)$ can
be identified as a sub-Hopf algebra of $D(H_r, H_l)$; hence we have
a Hopf algebra morphism from $D(H'_r, H'_l)$ to $H_{\sigma}$. The
image of this morphism is $H'_{\sigma}$, so $U_q({\rm sl}_N)$ is
also a factor of a generalized quantum double.
\end{remark}

\begin{corollary}\colabel{4.10}
The evaluation pairings on    $U_q^{\rm ext}({\rm sl}_N) \ot {\rm
SL}_q(N)$ and $ U_q({\rm sl}_N) \ot {\rm SL}_q(N)$,  are given by
\begin{eqnarray*}
&&\le \hk_i, x_{mn}\ri = z q^{  \d_{i, m}}\d_{m, n},~~
\le \hk_i^{-1}, x_{mn}\ri = z^{-1} q^{ - \d_{i, m}}\d_{m, n},\\
&&\le E_s, x_{mn}\ri =(q - q^{-1})\d_{s+1, n}\d_{s, m},~~ \le F_s,
x_{mn}\ri =(q - q^{-1})^{-1}\d_{s, n}\d_{s+1, m};
\end{eqnarray*}
and
\begin{eqnarray*}
&&\le K_i, x_{mn}\ri =q^{\d_{i, m} - \d_{i+1, m}}\d_{m, n},~~
\le K_i^{-1}, x_{mn}\ri =q^{\d_{i+1, m} - \d_{i, m}}\d_{m, n},\\
&&\le E_i, x_{mn}\ri =(q - q^{-1})\d_{i+1, n}\d_{i, m},~~ \le F_i,
x_{mn}\ri =(q - q^{-1})^{-1}\d_{i, n}\d_{i+1, m},
\end{eqnarray*}
respectively.
\end{corollary}
\begin{proof}
Both $U_q^{\rm ext}({\rm sl}_N)$ and $U_q({\rm sl}_N)$ are Hopf
subalgebras of ${\rm SL}_q(N)^0$, so  the evaluation map gives dual
pairings.  The formulas for the pairings come directly from
\leref{4.1} or from the proof of \prref{4.6}.
\end{proof}
From \thref{1.8}, ${\rm SL}_q(N)$ is a braided Hopf algebra in the
braided category of left Yetter-Drinfeld modules over $U_q^{\rm
ext}({\rm sl}_N)^{\rm cop}$. We  end this paper by computing the
structures of the braided Hopf algebra $\un{{\rm SL}_q(N)}$. First
we need the following.

\begin{lemma}\lelabel{4.10}
For any $1\leq s\leq N-1$ and $1\leq m, n\leq N$ we have
\begin{eqnarray}
\bl_{s+1,s}(S(x_{mn}))&=&
- z^{-1}(q - q^{-1})\d_{s+1, n}\d_{s, m},\label{hsigma5}\\
\br_{s,s+1}(S(x_{mn}))&=&- z^{-1}q^{-2  }(q - q^{-1})
\d_{s, n}\d_{s+1, m}.\label{hsigma6}
\end{eqnarray}
\end{lemma}
\begin{proof}
As usual, we only prove (\ref{hsigma5}); (\ref{hsigma6}) can be
proved similarly. Recall that   $\sigma^{-1}$ is obtained by
replacing $q$ by $q^{-1}$ and $z$ by $z^{-1}$ in (\ref{smqn}) and
(\ref{nsmqn}). We compute
\begin{eqnarray*}
\bl_{s+1,s}(S(x_{mn}))&=&\sigma (S(x_{mn}), x_{s+1,s})
=\sigma^{-1}(x_{mn}, x_{s+1,s})\\
&=& z^{-1}[ q^{-\delta_{m, s+1}}\delta_{m, n}\delta_{s, s+1}
+  [s+1>m](q^{-1} - q)\d_{s+1, n}\d_{s, m}]\\
&=&- z^{-1} (q - q^{-1})\d_{s+1, n}\d_{s, m},
\end{eqnarray*}
since $[s+1>m]\d_{s, m}=\d_{s, m}$, so the proof is complete.
\end{proof}

Now we can describe concretely the left Yetter-Drinfeld module
structure of $\un{{\rm SL}_q(N)}$ over $U_q^{\rm ext}({\rm
sl}_N)^{\rm cop}$.

\begin{proposition}\prlabel{6.14}
$\un{{\rm SL}_q(N)}$ is a left Yetter-Drinfeld module over $U_q^{\rm
ext}({\rm sl}_N)^{\rm cop}$ via the structure
\begin{eqnarray*}
&&\hspace*{-4mm} \hk_i\tr x_{mn}=q^{\d_{i, n} - \d_{i, m}}x_{mn},~~
\hk_i^{-1}\tr x_{mn}=q^{\d_{i, m} - \d_{i, n}}x_{mn},\\
&&\hspace*{-4mm}
E_s\tr x_{mn}=(1 - q^{-2})[q^{-1}\d_{s+1, n}x_{ms}
- q^{\d_{s, n} - \d_{s+1, n}}\d_{s, m}x_{s+1n}],\\
&&\hspace*{-4mm} F_s\tr x_{mn}=q(1 -q^{-2})^{-1}[q^{\d_{s, m} -
\d_{s+1, m}}\d_{s, n}x_{ms+1} -
q^{-1}\d_{s+1, m}x_{sn}],\\
&&\hspace*{-4mm} x_{mn}\longmapsto \sum \limits_{j=n}^Nq^{2(j -
n)}{\cal E}_{j, n}   \hk_j \left (\hk_m^{-1}\ot x_{mj} \right. \\
&&\hspace*{1.2cm}
+ \sum\limits_{ 1\leq i\leq m - 1}  \sum\limits_{ p\in {\cal B}_{i,
m}} \left.  (-q)^{-l(p) +    (m-i)}    {\cal E}_{m, p(m)}\cdots
{\cal E}_{i+1, p(i+1)} \hk_i^{-1}\ot x_{ij} \right ),
\end{eqnarray*}
where the ${\cal E}_{i,j}$ are from Definition \ref{braidedcomm}.
Also for $i\leq m-1$, ${\cal B}_{i, m}$ denotes the set of bijective
maps $p: \{i+1, \cdots, m\}\ra \{i,\cdots , m-1\}$ such that
$p(k)\leq k$, for all $i+1\leq k\leq m$.
\end{proposition}
\begin{proof}
We apply \prref{3.2} to this setting. By (\ref{mhs}), and recalling
again from (\ref{blr}) and (\ref{blrinv}) that  $\hk_i^{\pm
1}(x_{mn})=(zq^{ \d_{i, m}})^{\pm1}\d_{m, n}$   for   $1\leq i, m,
n\leq N$, we have
\begin{eqnarray*}
\hk_i^{\pm 1}\tr x_{mn}&=& \sum\limits_{j, k=1}^N\hk_i^{\pm
1}(\smi(x_{mj})S^{-2}(x_{kn}))x_{jk} =\sum\limits_{j,
k=1}^N\hk_i^{\mp 1}(x_{mj})\hk_i^{\pm 1}(x_{kn})x_{jk}\\ & =&
%\sum\limits_{j,k=1}^N z^{-1}q^{-\d_{im}}\d_{mj} z
%q^{\d_{ik}}\d_{kn}x_{jk}
q^{\pm(\d_{i, n} - \d_{i, m})}x_{mn}.
\end{eqnarray*}

Next, using (\ref{hsigma5}), we compute
\begin{eqnarray*}
E_s(S(x_{mn}))&=&\le \bl_{s+1s+1}^{-1}\bl_{s+1s},
S(x_{mn})\ri=\sum\limits_{i =1}^N
\hk_{s+1}(x_{in})\bl_{s+1s}(S(x_{mi}))\\
&  = & z q^{\d_{s+1, n} }\bl_{s+1s}(S(x_{mn}))=
-q^{\d_{s+1, n}}(q - q^{-1})\d_{s, m}\d_{s+1, n}\\
&=&-q(q - q^{-1})\d_{s, m}\d_{s+1, n}.
\end{eqnarray*}
Now, we use  the above computation together with the fact that
$E_s(x_{mn})=(q - q^{-1})\d_{s+1, n}\d_{s, m}$ and $K_s(x_{mn}) =
q^{\d_{s+1,m} - \d_{s,m}}\d_{s,m}$ from the proof of \prref{4.6}, to
compute the action of $E_s$ on $x_{mn}$.
\begin{eqnarray*}
&&\hspace*{-1.7cm}
E_s\tr x_{mn}\\
&{{(\ref{mhs})}\atop =}&\hspace*{-2mm}
\sum\limits_{i, j=1}^NE_s(\smi(x_{mi})S^{-2}(x_{jn}))x_{ij}\\
&{{(\ref{hls6}, \ref{mq7})}\atop =}&\hspace*{-2mm}
\sum \limits_{i, j=1}^N [q^{2(j - n)}\d_{m, i}E_s(x_{jn}) + q^{2(m - i) +
2(j - n)}E_s(S(x_{mi}))\hk_{s+1}^{-1}\hk_s(x_{jn})]x_{ij}\\
&=&\hspace*{-2mm}
q^{-1}(q - q^{-1})[q^{-1}\d_{s+1, n}x_{ms} - \d_{s, m}
\sum\limits_{j=1}^Nq^{2(j - n)}K_s(x_{jn})x_{s+1j}]\\
&=&\hspace*{-2mm}
(1 - q^{-2})[q^{-1}\d_{s+1, n}x_{ms} - q^{\d_{s, n} -
\d_{s+1, n}}\d_{s, m}x_{s+1n}],
\end{eqnarray*}
as   claimed. Similarly, one can compute $F_s\tr x_{mn}$; the
verification of the details is left to the reader.

Finally, from (\ref{chs}) and (\ref{mq7}), we see that
the coaction of $U_q^{\rm ext}({\rm sl}_N)^{\rm cop}$ on ${\rm
SL}_q(N)$ is defined by
\[
x_{mn}\mapsto \sum \limits_{i,
j=1}^N\bl_{S^{-1}(x_{mi})S^{-2}(x_{jn})}\ot x_{ij}
 = \sum\limits_{1\leq i, j \leq N}q^{2(j - n) + 2(m -
i)}\bl_{jn}\bl_{S(x_{mi})}\ot x_{ij}.
\]
From (\ref{mq6}) we have
\begin{eqnarray*}
\bl_{S(x_{mi})}&=&(-q)^{i - m} \sum \limits_{p\in S_{i,
m}}(-q)^{-l(p)}\bl_{x_{1,p(1)}\cdots
x_{i-1,p(i-1)}x_{i+1,p(i+1)}\cdots x_{N,p(N)}}\\
&=&(-q)^{i - m}\sum \limits_{p\in S_{i, m}}(-q)^{-l(p)}
\bl_{N,p(N)}\cdots \bl_{i+1,p(i+1)}\bl_{i-1,p(i-1)}\cdots
\bl_{1,p(1)},
\end{eqnarray*}
and this forces $k\geq p(k)$, for any $  k\not=i $. Since $p$ is
bijective,   non-zero summands occur only when
\begin{itemize}
\item[(i)]   $i=m$ and  $p(k)=k$, for all $  k\not=i $;
\item[(ii)]  $i\leq m-1$ and
$p(1)=1, \cdots , p(i-1)=i-1$, $\{p(i+1), \cdots , p(m)\}=\{i,
\cdots , m-1\}$, $p(m+1)=m+1, \cdots , p(N)=N$, and $p(k)\leq k$,
for any $k\in \{i+1, \cdots , m\}$.
\end{itemize}
In other words, we have proved that $\bl_{S(x_{mm})}=\bl_{mm}^{-1}$, and
that for $1\leq i\leq m - 1$,
\begin{eqnarray*}
&&\hspace*{-3mm}
\bl_{S(x_{mi})}\\
&&=(-q)^{i-m}\sum\limits_{p\in {\cal B}_{i,
m}}(-q)^{-l(p)} \bl_{NN}\cdots \bl_{m+1,m+1}\bl_{m,p(m)}\cdots
\bl_{i+1,p(i+1)}\bl_{i-1,i-1}\cdots \bl_{11}.
\end{eqnarray*}
From (\ref{mge}), for  $i+1\leq k\leq m$, $m+1\leq t\leq N$ and
$1\leq s\leq i-1$ we have
\begin{equation}\label{hsigma8}
\bl_{kp(k)}\bl_{tt}=\bl_{tt}\bl_{kp(k)}~~{\rm and}~~
\bl_{ss}\bl_{kp(k)}=\bl_{kp(k)}\bl_{ss},
\end{equation}
so that, using (\ref{hl3'}), we have
\[
\bl_{S(x_{mi})}=(-q)^{i - m}\sum\limits_{p\in {\cal B}_{i, m}}
(-q)^{-l(p)}\bl_{mp(m)}\cdots \bl_{i+1p(i+1)} \bl_{ii}^{-1}\cdots
\bl_{mm}^{-1},
\]
for $1\leq i\leq m - 1$. Now, by \coref{4.3} we have
$\bl_{kp(k)}={\cal E}_{k, p(k)}\hk_k$, for all $k\in \{i+1, \cdots ,
m\}$. In particular, if $i=m - 1$ we then have
\[
\bl_{S(x_{mi})}=-q^{-1}\bl_{mp(m)}\hk_{m-1}^{-1}\hk_m^{-1}=
-q^{-1}{\cal E}_{m, p(m)}\hk_i^{-1}.
\]
Now let $1\leq i\leq m - 2$. From (\ref{hl2s}) it follows that
$\hk_sE_t=E_t\hk_s$, for all $s>t+1$, so from the definition of
${\cal E}_{i, j}$ we conclude that $\hk_{m - t}{\cal E}_{k,
p(k)}={\cal E}_{k, p(k)}\hk_{m - t}$, for all $0\leq t\leq m - i -
2$ and $i + 1\leq k\leq m - t - 1$. Consequently,
\[
\bl_{S(x_{mi})}=(-q)^{i - m}\sum\limits_{p\in {\cal B}_{i, m}}
(-q)^{-l(p)}{\cal E}_{m, p(m)}\cdots {\cal E}_{i+1,
p(i+1)}\hk_i^{-1},
\]
for all $1\leq i\leq m - 1$.  Again using \coref{4.3} for $\bl_{jn}$
and the   formulas above for $\bl_{S(x_{mi})}$ we obtain  the
coaction in the statement of the theorem.
\end{proof}

We now apply the results of \thref{1.8}   to obtain the
multiplication, comultiplication, unit, counit and antipode for
$\un{{\rm SL}_q(N)}$.

\begin{theorem}\thlabel{6.15}
The structure of $\un{{\rm SL}_q(N)}$ as a braided Hopf algebra  in
the category of left $U_q^{\rm ext}({\rm sl}_N)^{\rm cop}$
Yetter-Drinfeld modules is the following. The unit and counit are as
in  ${\rm SL}_q(N)$. The multiplication is defined by
\begin{eqnarray*}
&&\hspace*{-5mm}
x_{im}\circ x_{jn}=q^{\d_{i, n} - \d_{m, n}}x_{im}x_{jn}
+ (q - q^{-1})\d_{i, n}\sum\limits_{s>n}q^{2(s - i) - \d_{s, m}}x_{sm}x_{js}\\
&&\hspace*{-3mm} - [m>n]q^{\d_{i, n}}(q - q^{-1})x_{in}x_{jm} -
[m>n+1](q - q^{-1})^2\d_{i, n}\sum\limits_{n<s<m}q^{2(s -
i)}x_{ss}x_{jm}.
\end{eqnarray*}
Comultiplication is given by
\begin{eqnarray*}
&&\hspace*{-5mm}
\un{\Delta}(x_{im})=q^{-2N - 1}\{q^{- \d_{i, m}}\sum\limits_s
q^{2s + \d_{s, i} + \d_{s, m}}
x_{is}\ot x_{sm} + (q - q^{-1})[x_{im}\ot x_{+}^{>m}\\
&& - \d_{i, m}\sum\limits_{s; t>i}q^{2(s + \d_{s, i})}
x_{ts}\ot x_{st} + x_{+}^{>i}\ot x_{im}] + (q - q^{-1})^2
[\sum\limits_{s>i, m}q^{2s - 1}x_{sm}\ot x_{is}\\
&& - [i>m]q^{2i + 1}\sum\limits_{s>i}x_{sm}\ot x_{is}
- [m>i]q^{2m + 1}\sum\limits_{s>m}x_{sm}\ot x_{is}]\\
&& - (q - q^{-1})^3\sum \limits_{t>s>i, m}q^{2s}x_{tm}\ot x_{it}\},
\end{eqnarray*}
where we denoted $x_{+}^{>k}:=\sum\limits_{s>k}q^{2s}x_{ss}$. The
antipode is determined by
\begin{eqnarray*}
\un{S}(x_{im})=q^{2N+1}
[q^{-2m - \d_{i, m}}S(x_{im}) - (q - q^{-1})\d_{i, m}S(x_{-}^{>m})],
\end{eqnarray*}
where $x_{-}^{>m}:=\sum\limits_{s>m}q^{-2s}x_{ss}$ and $S$ is the antipode of
${\rm SL}_q(N)$.
\end{theorem}
\begin{proof}
From \leref{4.1} we have
\begin{eqnarray*}
\bl_{ij}(x_{mn})&=&zq^{ \d_{m, i}}\d_{m, n}\d_{i, j} + [i>j]z(q -
q^{-1})
\d_{n, i}\d_{m, j},\label{blsec}\\
\bl^{-}_{ij}(x_{mn})&=& z^{-1}q^{  -\d_{m, i}}\d_{m, n}\d_{i, j} -
[i>j]z^{-1} (q - q^{-1})\d_{n, i}\d_{m, j},\label{blmin}
\end{eqnarray*}
where $\bl^{-}_{ij}$ denotes the map $\sigma^{-1}(-,
x_{ij})=\bl_{S^{-1}(x_{ij})}$.

Now, the multiplication of $\un{{\rm SL}_q(N)}$  from \thref{1.8}
can be rewritten as
\[
x\circ y=\sigma^{-1}(x_3, y_2)\sigma(S^2(x_1), y_3)x_2y_1,
\]
and together with (\ref{mq7}) this allows us to compute
\begin{eqnarray*}
x_{im}\circ x_{jn}&=&\sum\limits_{s, t, u, v=1}^N
q^{2(s - i)}\sigma^{-1}(x_{tm}, x_{uv})
\sigma (x_{is}, x_{vn})x_{st}x_{ju}\\
&=&\sum\limits_{s, t, u, v=1}^N
q^{2(s - i)}\bl^{-}_{uv}(x_{tm})\bl_{vn}(x_{is})x_{st}x_{ju}\\
&=&\sum\limits_{s, t, u, v=1}^Nq^{2(s - i)}
\left(q^{-\d_{t, u}}\d_{t, m}\d_{u, v} - [u>v](q - q^{-1})\d_{m, u}\d_{v, t}\right)\\
&&\hspace*{5mm}\times
\left(q^{\d_{i, v}}\d_{i, s}\d_{v, n} + [v>n] (q - q^{-1})\d_{s, v}\d_{i, n}\right)x_{st}x_{ju}.
\end{eqnarray*}
Splitting the sum above into  four separate sums, we obtain the
 formula in the statement of the theorem.

The computation of $\un{\Delta}(x_{im})$ is much more complicated. Firstly,
observe that the comultiplication in \thref{1.8} can be rewritten as
\begin{equation}\label{comultunslqn}
\un{\Delta}(x)=\sigma(x_1, x_6)\sigma^{-1}(x_2, x_8)\sigma(x_4, x_9)v(x_5)x_3\ot x_7,
\end{equation}
where $v$ is the map    $v(y) = \sigma(y_1,S(y_2))$ as in (\ref{v}).
Next, we compute
\begin{eqnarray*}
v(x_{im})&=&\sum\limits_{k=1}^N
\sigma(x_{ik}, S(x_{km}))
=\sum\limits_{k=1}^N
q^{2(m - k)}\bl^{-}_{km}(x_{ik})\\
&=&\sum\limits_{k=1}^N z^{-1}q^{2(m - k)  }(q^{-\d_{i, k}}\d_{i,
k}\d_{k, m}
- [k>m](q - q^{-1})\d_{i, m})\\
&=&z^{-1} [q^{-\d_{i, m}} - q(1 - q^{-2})\sum\limits_{k>m}(q^{-2})^{k - m}]\d_{i, m}\\
&=&z^{-1}[q^{-1} - q^{-1}(1 - (q^{-2})^{N - m})]\d_{i, m} =z^{-1}q^{
- 2(N - m) - 1}\d_{i, m},
\end{eqnarray*}
for all $1\leq i, m\leq N$. Therefore, by (\ref{comultunslqn}) we have
\begin{eqnarray*}
\un{\Delta}(x_{im})&=&\sum\limits_{a, \cdots , h=1}^N
\bl_{ef}(x_{ia})\bl^{-}_{gh}(x_{ab})
\bl_{hm}(x_{cd})v(x_{de})x_{bc}\ot x_{fg}\\
&=&\sum\limits_{a, b, c, d, f, g, h=1}^N
q^{-2(N - d) - 1}\left(q^{\d_{i, d}}\d_{i, a}\d_{d, f} +
[d>f](q - q^{-1})\d_{d, a}\d_{i, f}\right)\\
&&\hspace*{1cm}
\times \left(q^{-\d_{a, g}}\d_{a, b}\d_{g, h} - [g>h](q - q^{-1})\d_{g, b}\d_{h, a}\right)\\
&&\hspace*{1cm}
\times \left(q^{\d_{c, h}}\d_{c, d}\d_{h, m} + [h>m](q - q^{-1})\d_{h, d}\d_{c, m}\right)
x_{bc}\ot x_{fg}.
\end{eqnarray*}
Again, splitting the above sum into  eight separate sums, a tedious
but straightforward computation yields the formula for
$\un{\Delta}(x_{im})$ in the statement of the theorem.

Finally, the equation for $\un{S}(x_{ij})$ follows easily by the
general formula in \thref{1.8} and by computations similar to the
ones above.
\end{proof}
%%%%%%%%%%%%%%%%%%%%%%%%%%%%%%%%%%%%%%%%%%%%%%%%%%%%%%%%%%%%

\end{document}